\documentclass[11pt, a4paper]{article}
\usepackage[utf8]{inputenc}
\usepackage{amssymb,amsmath,dsfont}
\usepackage{graphicx}
\usepackage{wrapfig}
\usepackage{theorem} 
\usepackage[margin=2.5cm]{geometry}
\usepackage[dvipsnames,usenames]{color}
\usepackage{epstopdf}
\usepackage[font=small,labelfont=bf]{caption}

\usepackage{float}

\usepackage[colorlinks=true, linkcolor=blue, urlcolor=magenta, citecolor=blue]{hyperref}
{\theorembodyfont{\rmfamily}\newtheorem{test}{Test}}
\newtheorem{lemma}{Lemma}[section]
\newtheorem{theorem}[lemma]{Theorem}
\newtheorem{prop}[lemma]{Proposition}
\newtheorem{coro}[lemma]{Corollary}
{\theorembodyfont{\rmfamily}\newtheorem{remark}[lemma]{Remark}}

\title{Inverse problem of reconstruction of degenerate diffusion coefficient in a parabolic equation}
\author{Piermarco Cannarsa, Anna Doubova, Masahiro Yamamoto}
\date{\today}

\newcommand{\ep}{\varepsilon}
\newcommand{\www}{\widetilde}
\newcommand{\ppp}{\partial}

\newcommand{\ppppa}{\partial_x(x^{\alpha}a(x)\partial_x}

\newcommand{\ppppp}{\partial_x(p(x)\partial_x}

\newcommand{\la}{\lambda}
\newcommand{\va}{\varphi}
\newcommand{\R}{\mathds{R}}

\newcommand{\N}{\mathds{N}} 

\newcommand{\wwww}{e^{s\va}}
\newcommand{\weight}{e^{2s\varphi}}

\newcommand{\exex}{e^{2se^{\la\delta^2}}}

\begin{document}

\maketitle 

\begin{abstract}
We consider the inverse problem of identification of degenerate diffusion coefficient of the form $x^\alpha a(x)$ in a one dimensional parabolic equation by some extra data.  We first prove by energy methods the uniqueness and Lipschitz stability results for the identification of a constant coefficient $a$ and the power $\alpha$ by knowing an interior data at some time. 

On the other hand, we obtain the uniqueness result for the identification of a general diffusion coefficients $a(x)$ and also the power $\alpha$  form a boundary data on one side of the space interval. The proof is based on global Carleman estimates for a hyperbolic problem and an inversion of the integral transform similar to the Laplace transform. Finally, the theoretical results are satisfactory verified by numerically experiments. 
\end{abstract}

\section{Introduction}


\bigskip
In this paper we will consider several inverse problems for the following degenerate parabolic equation with two type of degeneracy in the diffusion coefficient. Mainly, weak degeneracy, with $0< \alpha < 1$
\begin{equation}\label{eq.dw}
\left\{\begin{array}{ll}
\partial_t u- \partial_x(x^{\alpha} a(x) \partial_x u) = f(x,t), & (x,t)\in(0,\ell)\times(0,T),  \\[1mm]
u(0,t) = 0, \quad u(\ell,t) = 0, &  t\in (0,T),  \\[1mm]
u(x,0) = u_0(x),  & x\in(0,\ell),
\end{array}\right.
\end{equation}
and strong degeneracy, with $1\leq \alpha <2$   
\begin{equation}\label{eq.sd}
\left\{\begin{array}{ll}
\partial_t u- \partial_x(x^{\alpha} a(x) \partial_x u) = f(x,t), & (x,t)\in(0,\ell)\times(0,T), \\[1mm]
(x^\alpha\partial_x  u(x,t)|_{x=0}  = 0, \quad u(\ell,t) = 0, &  t\in (0,T), \\[1mm] 
u(x,0) = u_0(x),  & x\in(0,\ell),
\end{array}\right.
\end{equation}
where $T>0$, $\ell >0$, $u_0\in L^2(0,\ell)$, $u_0\neq 0$,  $f\in L^2((0,\ell)\times(0,T))$ are given, $a\in C^1([0,\ell])$ with 
$0<\underline{a} \leq a(x)\leq \overline{a}<+\infty$ for all $x\in(0,\ell)$, with $\underline{a}$ and $\overline{a}$ a positive constants. 

\medskip
Our goal is to determine or estimate  the diffusion coefficient for the considered degenerate problems. Although we can discuss more general cases, we concentrate here on the one-dimensional linear  equation. The following three types of inverse problems of the identification of the unknown diffusion coefficient will be considered and analyzed. 

\begin{description}

\item[Linear case:] given $\alpha=1$, we will deal with the inverse problem of identification of a constant coefficient, $a(x)\equiv a$, where $a$ is unknown positive constant, form given suitable internal data at some time  (see Section~\ref{sec.lincase}). 
\item[Power-like case:] given $a(x) = 1$, the inverse problem will be the identification of the power $\alpha\in (0,2)$ from suitable internal observations at some time (see Section~\ref{sec.power}). 
\item[General case:]  the inverse problem of the identification of a general diffusion coefficient $x^\alpha a(x)$ form one boundary observation (see Section~\ref{sec.gencase}). 
\end{description}

\medskip
In recent years, degenerate parabolic equations have received increasing attention, in view of the important related theoretical analysis and practical applications, as for example climate science (see Sellers~\cite{Sellers}, D\'{\i}az~\cite{diaz}, Huang~\cite{JiHuang}),  populations genetics Ethier~\cite{Ethier}, vision Citti and Manfredini~\cite{Citti}, financial mathematics~Black and Scholes~\cite{BS} and fluid dynamic Oleinik and Samokhin~\cite{Oleinik}. See also Cannarsa, Martinez and Vancostenoble~\cite{CMV} and references therein. 

\medskip
Inverse problems are the kind of problems known as not well-posed in the sense of Hadamard (see \cite{Hadamard}). This means that, either a solution does not exist, or it is not unique and/or small errors in known data lead to large errors in the computation of the solution(s). 
In particular, the first theoretical issues that we analyze in our paper are uniqueness and stability.  The main techniques used here are usual in control and parameter identification theory: tools from real and complex analysis, energy methods and Hardy's inequality, Laplace and other similar integral transforms, and Carleman inequalities.

\medskip
The literature related to  inverse problems for degenerate parabolic PDEs, despite their fundamental importance and practical applications, is rather recent and scarce. For example, the inverse  source problem was considered in Tort~\cite{Tort1}, Cannarsa, Tort, and Yamamoto~\cite{CTY}, Cannarsa, Martinez and Vancostenoble~\cite{CMV}, Deng and all~\cite{Deng} and Hussein and all~\cite{HusseinLesnic}, where numerical reconstruction was also considered.  The inverse problem of recovering the first-order coefficient of degenerate parabolic equations is analyzed in~Deng and Yang \cite{Deng2} and Kamynin~\cite{Kamynin}. An inverse diffusion problem was considered in~\cite{Tort2}, where a constant diffusion coefficient was recovered from mesurement of second order derivatives by means of Carleman estimates provided in Cannarsa, Tort and Yamamoto~\cite{CTY}. 

As for inverse problems of determining spatially 
varying coefficients such as conductivities and source terms for 
non-degenerate parabolic equations, we note that 
there has been a substantial amount of works.  Here we are restricted to 
refer to chapters from Isakov \cite{Is}, and Yamamoto \cite{Ya}.
See also the references therein.  In particular,
as for the same type of inverse problems for one-dimensional regular 
parabolic equations, we refer to Murayama \cite{Mura},  see also
Pierce \cite{Pier}, Suzuki and Murayama \cite{SuMu}.

On the other hand, most of the techniques that have been developed to treat non-degenerate parabolic equations cease to be applicable in the degenerate case.
For instance, for strongly degenerate operators the trace of the co-normal derivate must vanish on the part of the boundary where ellipticity fails, yielding no useful measurement\footnote{Although in this paper we consider equations that degenerate on a strict subset of the boundary, our approach applies unaltered  when degeneracy occurs on the whole boundary, like in the Budyko-Sellers model~\cite{Sellers}.}.

\medskip
The second main issue of this paper is the reconstruction method for our inverse problems. The goal is to compute approximations of the diffusion coefficient as solutions to the inverse problems using the observation data. 
To our best knowledge, there are no works combining the theoretical studies on the uniqueness, stability 
and the numerical reconstruction of the diffusion coefficient for a degenerate parabolic equation. 

An efficient technique to achieve this, as shown below, relies on a reformulation of the search of the diffusion coefficient as 
an extremal problem. This is classical nowadays and has been applied in a lot of situations; see for instance Lavrentiev and all ~\cite{Lavrentiev}, Samarskii and  Vabishchevich~\cite{Samarskii}, and Vogel~\cite{Vogel}.

%
%
 \bigskip
 The paper is organized  as follows. In Section~\ref{sec:weelpos}, we consider the well-posedness of the direct problem \eqref{eq.dw} and \eqref{eq.sd}. In Section~\ref{sec.lincase} we prove, using the energy methods, the Lipschitz stability and the uniqueness result for the linear case corresponding to the constant coefficient $a$. Section~\ref{sec.power} will be devoted to the power-like case of the identification of the power $\alpha$ in  \eqref{eq.dw} and \eqref{eq.sd}. Again, the Lipschitz stability and the uniqueness will be achieved by the energy methods. In Section~\ref{sec.gencase} a general case of the reconstruction of the unknown general diffusion coefficient  $x^\alpha a(x)$ will be treated by means of the Carleman estimates, obtaining the uniqueness result for both $\alpha$ and $a(x)$. 
Finally, in Section~\ref{sec.numerics1}, in order to illustrate theoretical results obtained in the previous sections, we perform satisfactory numerical experiments corresponding to the considered inverse problems.  
 

 \section{Preliminaries}\label{sec:weelpos}

\medskip
In this section we will consider the well-posedness of the direct problem related to \eqref{eq.dw} and \eqref{eq.sd}. 
We begin with recalling the natural functions spaces where the  problem can be set. For any $\ell>0$,  let us call $H=L^2(0,\ell)$. For $\alpha\in(0,2)$ we will consider the following function spaces: 
\begin{equation}\label{eq.spaces}
H_\alpha^1(0,\ell) =
\left\{\begin{array}{l}
 \{ u\in H \,  |\,    \displaystyle \int_0^\ell x^\alpha|\partial_x u|^2\, dx<\infty,\,  u(0) = u(\ell)= 0\},  \quad 0<\alpha<1,
 \\[4mm]
 \{ u\in H \,  |\,    \displaystyle \int_0^\ell x^\alpha|\partial_x u|^2\, dx<\infty,\,   u(\ell) =  0\} ,\quad  1\le \alpha<2
 \end{array}\right.
\end{equation}
and 
\begin{equation*}
H_\alpha^2(0,\ell) =
 \{ u\in H^1_\alpha(0,\ell) \,  |\,  x\mapsto x^\alpha \partial_x u\in H^1(0,\ell)  \},  \quad 0<\alpha<2.
\end{equation*}

\begin{remark}[Neumann boundary condition]\label{Neumann boundary condition}
For $u\in H^2_\alpha(0,\ell)$ with $1\le \alpha < 2$, we have $x^\alpha \partial_x u|_{x=0}=0$. Indeed, if
$x^\alpha\partial_x u(x) \to L$ when $x\to 0$, then $x^\alpha|\partial_x u(x) |^2 \sim L^2/x^\alpha$ and, therefore $L=0$ otherwise  $u\notin H^1_\alpha(0,\ell)$.
\end{remark}
We now adapt the classical Poincar\'e inequality to the above weighted spaces.
\begin{lemma}[Poincare's inequality] \label{lemma1} For all $\varepsilon>0$, there exists a  constant $C_p=C_p(\varepsilon,\ell)>0$ such that for all $\alpha\in [0,2-\varepsilon]$, the following holds: 
\begin{equation}\label{PI}
\int_0^\ell |u|^2\, dx \leq C_p\int_0^\ell x^\alpha |\partial_x u |^2\, dx \quad \forall\, u\in H_\alpha^1(0,\ell).
\end{equation}
\end{lemma}
\textbf{Proof of Lemma~\ref{lemma1}--} Let $\ell>0$ be given. We will distinguish two cases. 
\begin{enumerate}
\item Assume $\alpha=1$, for all $x\in (0,\ell)$, we have 
 \begin{equation*}
 \begin{array}{lll}
 |u(x) | &\le &\displaystyle \int_x^\ell |\partial_x u| \, dy \le \Big(\int_x^\ell y |\partial _x u(y)|^2\, dy \Big)^{1/2}
 \Big(\int_x^\ell \dfrac{1}{y} \, dy\Big)^{1/2} 
 \\[4mm]
 &\le &
 \displaystyle \Big(\log \frac{1}{\ell}-\log \dfrac{1}{x}\Big)^{1/2}  \Big(\int_x^\ell y |\partial _x u(y)|^2\, dy \Big)^{1/2}.
 \end{array}
\end{equation*}
Therefore, we obtain 
 \begin{equation*}
 \int_0^\ell |u(x) |^2\, dx   \leq \Big(\int_0^\ell \log \dfrac{1}{x}\, dx\Big)   \int_0^\ell x |\partial _x u(x)|^2\, dx = 
 \ell(1-\log \ell) \int_0^\ell x |\partial _x u(x)|^2\, dx.
\end{equation*}


\item Assume now that $\alpha \neq 1$. We can write 
\begin{equation*}
\begin{array}{lll}
|u(x) | &\leq & \left | \displaystyle\int_x^\ell y^{\alpha/2} u'(y) \, \dfrac{dy}{y^{\alpha/2}}\right |
\leq  \left ( \displaystyle\int_x^\ell  \dfrac{dy}{y^{\alpha/2}} \right)^{1/2}
\displaystyle\left(\int_0^\ell  y^{\alpha} |u'(y)|^2 \, dy\right)^{1/2} 
\\[4mm]
 &\leq & \left(\dfrac{\ell^{1-\alpha}-x^{1-\alpha}}{1-\alpha} \right)^{1/2}
 \displaystyle\left(\int_0^\ell  y^{\alpha} |u'(y)|^2 \, dy\right)^{1/2} 
 \end{array}
\end{equation*}
Therefore, we obtain 
\begin{equation*}
\begin{array}{lll}
\displaystyle\int_0^\ell |u(x) |^2\, dx &= &
\left(\displaystyle\int_0^\ell  \dfrac{\ell^{1-\alpha}-x^{1-\alpha}}{1-\alpha} \, dx \right)
 \displaystyle\int_0^\ell  x^{\alpha} |\partial_x u(x)|^2 \, dx  
 \\[4mm]
&= & \dfrac{\ell^{2-\alpha}}{2-\alpha}
 \displaystyle\int_0^\ell  x^{\alpha} |\partial_x u(x)|^2 \, dx =C_p  \displaystyle\int_0^\ell  x^{\alpha} |\partial_x u(x)|^2 \, dx,
\end{array}
\end{equation*}
here we have used that $\alpha<2$. This ends the proof. 

\hfill $\square$
\end{enumerate}

\bigskip
For given constants $0<\underline{a}<\overline{a}$,   we assume that 
 \begin{equation}
\label{eq:SA}
a\in C^1([0,\ell]),\qquad 0<\underline{a} \leq a(x)\leq \overline{a}<+\infty\quad\forall x\in(0,\ell)
\end{equation}
Problems \eqref{eq.dw} and \eqref{eq.sd} can be recast  in the abstract form 
\begin{equation}
\label{eq:HCP}
\begin{cases}
 u'(t)=Au(t)
 &
 t\ge 0
 \vspace{.1cm}
 \\
 u(0)=u_0
\end{cases}
\end{equation}
by introducing the linear operator $A:D(A)\subset H \to H$ defined by
\begin{equation}
\label{eq:A}
D(A) = H_\alpha^2(0,\ell), \qquad Au = \partial_x(x^\alpha a(x) \partial_x u)\quad\forall u\in D(A).
\end{equation}
Observe that the degenerate coefficient $\widetilde a(x)=x^\alpha a(x)$ is of class $C^1([0,\ell])$ and satisfies
\begin{equation}\label{hypa}
\lim_{x\to 0}\dfrac{ x\widetilde a'(x)}{\widetilde a(x)} =\lim_{x\to 0}\dfrac{x^{\alpha+1}a'(x)+\alpha x^\alpha a(x)}{x^\alpha a(x)} = \alpha \in (0,2). 
\end{equation}
Therefore, we can appeal to~\cite{MarVan} to obtain a well-posedness result for  \eqref{eq:HCP} as well as for the nonhomogeneous problem
\begin{equation}
\label{eq:NHCP}
\begin{cases}
 u'(t)=Au(t)+f(t)
 &
 t\ge 0
 \vspace{.1cm}
 \\
 u(0)=u_0
\end{cases}
\end{equation}

So,  one can prove the following, where we recall that $D((-A)^{1/2})=H_\alpha^1(0,\ell)$ for the operator $A$ defined in \eqref{eq:A}.

\begin{theorem}\label{th2.3}
Assume \eqref{eq:SA}. Then,  $A$ is the infinitesimal generator of a strongly continuous semigroup of contractions, $e^{tA}$. Moreover, $e^{tA}$ is analytic.
Therefore, for any $u_0\in H$,  the solution $u(x,t) = (e^{tA}u_0)(x)$ of problem  \eqref{eq:HCP} satisfies:

\begin{itemize}
\item[(i)] $u\in C([0,\infty);H)$,
\item[(ii)] $t\mapsto u(\cdot, t)$ is analytic as a map $(0,+\infty) \longrightarrow H$,
\item[(iiI)] $u(t) \in \bigcap_{n\geq 1} D(A^n)$ \quad $\forall\, t>0$. 
\end{itemize}
Furthermore, if $u_0\in D((-A)^{1/2})$, then for any $f\in L^2(0,T;H)$ the mild solution
\begin{equation*}
u(t):=e^{tA}u_0+\int_0^te^{(t-s)A}f(s)ds\qquad(t\in[0,T])
\end{equation*}
of problem \eqref{eq:NHCP} belongs to $H^1(0,T;H)\cap C([0,T];D((-A)^{1/2}))\cap L^2(0,T;D(A))$ and satisfies the equation in \eqref{eq:NHCP} for a.e. $t\in[0,T]$.
\end{theorem}

Consequently, the equation in \eqref{eq.dw} and \eqref{eq.sd} is satisfied in classical sense on $(0,\ell)\times(0,+\infty)$, as well as boundary conditions, taking into account  Remark~\ref{Neumann boundary condition}. As for the initial condition, we recall that $u_0\in L^2(0,\ell)$ implies that $u\in C^0([0,+\infty);L^2(0,\ell))$ and 
$u_0\in H^1_\alpha(0,\ell)$ implies that $u\in C^0([0,+\infty);H^1_\alpha(0,\ell))$. Hereafter, we will assume solutions as smooth as required.

\section{Reconstruction by energy methods}

\subsection{The linear case: uniqueness and Lipschitz stability}\label{sec.lincase}

Let us consider the linear case in which $a(x) \equiv a$, where $a$ is a constant such that  $0<\underline{a}\le a \le \overline{a}$  and $\alpha = 1$. This corresponds to
the problem \eqref{eq.sd} with strong degeneracy, that is written as follows: 

\begin{equation}\label{dpe1}
\left\{\begin{array}{ll}
\partial_t u - a\, \partial_x( x\partial_x u)= f(x,t), & x\in(0,\ell), \quad t\in (0,T), \\[1mm]
x\, \partial_x u(x,t)|_{x=0} =0 \quad   u(\ell,t) = 0,  &  t\in (0,T),  \\[1mm]
u(x,0) = u_0(x),  & x\in(0,\ell).
\end{array}\right.
\end{equation}
In this section we will analyze the uniqueness and Lipschitz stability for the inverse problem of the identification of the constant coefficient $a$ in \eqref{dpe1} from measurements $\partial_t u(x,t_0)$ and $x\partial_x u(x,t_0)$, for some $t_0\in (0,T]$ and for all $x\in (0,\ell)$, as is indicated in Theorem~\ref{th:linearcase}. In Section~\ref{sec:lincase1} we will perform numerical reconstructions of $a$ related to this inverse problem.  

\medskip 
On the other hand, let us note that it is also possible to consider other observations, for example, using the boundary data $\partial_x u(\ell,t)$ for all $t\in (0,T)$. The related inverse problem will be the following: given $\eta=\partial_x u(\ell,\cdot)$, find $a$ such that the corresponding solution $u^a$ of \eqref{dpe1} satisfies $\partial_x u(\ell,t) = \eta(t)$, for all  $t\in (0,T)$. 
However,  uniqueness and stability for this second inverse problem are, to our knowledge, open questions that we will consider in a forthcoming paper.  Nonetheless,  we have satisfactory numerical approaches for this second inverse problem, that we will present in the Section~\ref{sec.numerics2}.

\bigskip
Let $u_i$, $i=1,2$ be the solution of \eqref{dpe1} with $a=a_i$, that is to say
\begin{equation}\label{dpe1i}
\left\{\begin{array}{ll}
\partial_t u_i - a_i\, \partial_x( x\partial_x u_i)= f(x,t), & x\in(0,\ell), \quad t\in (0,T), \\[1mm]
x\, \partial_x u_i(x,t)|_{x=0} =0 \quad   u_i(\ell,t) = 0,  &  t\in (0,T),  \\[1mm]
u_i(x,0) = u_0(x),  & x\in(0,\ell).
\end{array}\right.
\end{equation}

The main result of this section is the following: 
\begin{theorem}\label{th:linearcase}
Assume that $a_1,\, a_2 \in C^1([0,\ell])$, $0<\underline{a} \leq a_1,a_2\leq \overline{a}<+\infty$ for all $x\in(0,\ell)$, $u_0\in L^2(0,\ell)$, $u_0\neq 0$ and $f\in L^2((0,\ell)\times(0,T))$. Let $u_1,u_2$ be solutions of \eqref{dpe1i} satisfying the following hypothesis: 
for some $t_0\in (0,T]$ and $\mu>0$ such that 
 \begin{equation}\label{LB}
\int_0^\ell x |\partial_x u_i(x,t_0)|^2\, dx \ge  \mu, \quad i=1,2. 
 \end{equation}
Then there exists a positive constant $C=C(\underline{a},\mu)$ such that 
 \begin{equation}\label{LS}
|a_2-a_1| \le C \left (\int_0^\ell \big( |\partial_t u_1(x,t_0)-\partial_t u_2(x,t_0)|^2+ x|\partial_x u_2(x,t_0)-\partial_x u_1(x,t_0) |^2\big)\, dx\right)^{1/2}. 
 \end{equation}
\end{theorem}

\begin{remark}
Notice that, as we have mentioned before, we can obtain similar uniqueness result  for more general equations containing zero order therms.  
\end{remark}

\begin{remark}
Recall that, since $e^{tA}$ is analytic, $u_i(t,\cdot)\in H_1^2(0,\ell)$ for a.e. $t>0$, $i=1,2$.  Moreover, assuming in addition that $f\equiv 0$, we have the following.
\begin{enumerate} 
\item  By backward uniqueness we know that, if $u_0\neq 0$, then the solution $u$ of \eqref{dpe1} satisfies
\begin{equation}\label{eq:pos}
\displaystyle\int_0^\ell x |\partial_x u(x,t)|^2 \, dx >0  \quad \forall\,  t\in [0,T].  
\end{equation}
This explains the role played by assumption \eqref{LB} above.
\item Since $A$ is dissipative, we have that 
\[
\displaystyle\int_0^\ell x |\partial_x u_i(x,t)|^2 \, dx \ge \int_0^\ell x |\partial_x u_i(x,T)|^2\, dx \quad \forall\,  t\in [0,T].  
\]
Therefore, 
\[
\mu = \min_{i=1,2} \displaystyle\int_0^\ell x |\partial_x u_i(x,T)|^2 \, dx
\]
is a possible choice of $\mu$ for every $t_0\in (0,T]$. 
\item Another way to construct $\mu>0$ satisfying \eqref{LB} is to take $u_0\in D(A)$, observing that 
\begin{equation*}
\displaystyle\int_0^\ell x |\partial_x u(t,x)|^2\, dx 
\ge 
\displaystyle\int_0^\ell x |\partial_x u_0(x)|^2\, dx  
- 2t\displaystyle\int_0^\ell  |(x\partial_x u_0)_x|^2\, dx. 
\end{equation*}
Therefore,  \eqref{LB} holds with $\mu = \dfrac{1}{2} \displaystyle\int_0^\ell x |\partial_x u_0(x)|^2\, dx $ for all 
\[
0< t_0 < \dfrac{1}{4} \, 
\dfrac{\displaystyle\int_0^\ell x |\partial_x u_0(x)|^2\, dx}{\displaystyle\int_0^\ell  |(x\partial_x u_0)_x|^2\, dx}\,.
\]
\end{enumerate}

\end{remark}

\medskip
\noindent
%
\textbf{Proof of Theorem~\ref{th:linearcase}} -- We assume that $a_2\ge a_1$ and set $w = u_1-u_2$. Then, we can write  
 \begin{equation*}
\partial_t w - a_1\, \partial_x( x\partial_x w)= (a_1-a_2) \partial_x( x\partial_x u_2).
 \end{equation*}
Multiplying by $u_2$ and integrating by parts this equality, we get 
 \begin{equation*}
 \begin{array}{l}
\displaystyle (a_2-a_1)\int_0^\ell x|\partial_x u_2|^2\, dx 
=
 \displaystyle\int_0^\ell u_2 \partial_t w \, dx 
 + a_1\int_0^\ell x\partial_x w \partial_x u_2\, dx  \\[4mm]
\le 
\displaystyle \left(\int_0^\ell |u_2|^2\, dx  \right)^{1/2} \left(\int_0^1 |\partial_t w|^2\, dx  \right)^{1/2} 
+ 
a_1  \left(\int_0^\ell x |\partial_x w|^2\, dx  \right)^{1/2} \left(\int_0^\ell x |\partial_x u_2|^2\, dx  \right)^{1/2}.
 \end{array}
 \end{equation*}
So, by \eqref{PI}, we have
 \begin{equation}\label{eq:ineq}
 \begin{array}{l}
 \displaystyle |a_2-a_1|\int_0^\ell x|\partial_x u_2|^2\, dx  
 \\[4mm]
\le 
\displaystyle \left(\sqrt{C_p} \left(\int_0^\ell |\partial_t w|^2\, dx  \right)^{1/2} + 
a_1  \left(\int_0^\ell x |\partial_x w|^2\, dx  \right)^{1/2}\right) \left(\int_0^\ell x |\partial_x u_2|^2\, dx  \right)^{1/2},
 \end{array}
 \end{equation}
 which in turn yields the the conclusion in view of \eqref{LB} and ends the proof of this theorem. 
 \hfill $\square$
 
 \bigskip

 It is clear that from the proof of this theorem we can deduce the following uniqueness result.
\begin{coro}
Let us assume that $a_1,\, a_2\ge \underline{a}>0$, $f\equiv 0$, and $u_0\in L^2(0,\ell)$, $u_0\neq 0$. If $u_1,u_2\in C^0([0,T];L^2(0,\ell))$ are solutions of \eqref{dpe1i} such that for some $t_0\in (0,T]$ we have 
$\partial_t u_1(x,t_0)=\partial_t u_2(x,t_0)$ and $ \partial_x u_2(x,t_0)=\partial_x u_1(x,t_0)$ for all $x\in(0,\ell)$, then 
$a_1=a_2$. 
\end{coro}
Indeed, it is enough to realize that we do not need  hypothesis \eqref{LB}  for uniqueness, because we do not need the constant in \eqref{eq:ineq} to be uniform. Therefore, form \eqref{eq:ineq} and thanks to the inequality \eqref{eq:pos}, we obtain uniqueness. 

\subsection{The power-like case: uniqueness and Lipschitz stability}\label{sec.power}

In this section we will consider the inverse problem of identification of the power $\alpha$ in the degenerate parabolic problems \eqref{eq.dw} and \eqref{eq.sd}, assuming for simplicity that $a(x)\equiv 1$ and $f\equiv 0$. That is to say,  that we will consider  
the following problems: with weak degeneracy  $(0\leq \alpha < 1)$   
\begin{equation}\label{eq.CP.w}
\left\{\begin{array}{ll}
\partial_t u- \partial_x(x^{\alpha}  \partial_x u) = 0, & (x,t)\in(0,\ell)\times(0,T),  \\[1mm]
u(0,t) = 0, \quad u(\ell,t) = 0, &  t\in (0,T),  \\[1mm]
u(x,0) = u_0(x),  & x\in(0,\ell)
\end{array}\right.
\end{equation}
and with strong degeneracy  $(1\leq \alpha <2)$   
\begin{equation}\label{eq.CP}
\left\{\begin{array}{ll}
\partial_t u- \partial_x(x^{\alpha}  \partial_x u) = 0, & (x,t)\in(0,\ell)\times(0,T), \\[1mm]
x^\alpha\partial_x  u(x,t)|_{x=0}  = 0, \quad u(\ell,t) = 0, &  t\in (0,T), \\[1mm] 
u(x,0) = u_0(x),  & x\in(0,\ell).
\end{array}\right.
\end{equation}

Similarly to what we have done in the Section~\ref{sec.lincase}, we will analyze  uniqueness and stability  for the inverse problem of identification of the power $\alpha$ in the above problems from measurements $\partial_t u(x,t_0)$ and $x\partial_x u(x,t_0)$, for some $t_0\in (0,T]$ and for all $x\in (0,\ell)$. In Section~\ref{sec:numpower1} we will perform the numerical reconstruction of $\alpha$ for the problems  \eqref{eq.CP.w} and \eqref{eq.CP}.

\medskip 

As we have noted in the previous section, is also possible to consider  the inverse problem of identification of the power $\alpha$ with only one boundary observation, for example  $\partial_x u(\ell,t)$ for all $t\in (0,T)$. However, the uniqueness and stability related to this second inverse problem is, to our knowledge, an open questions. Nevertheless,  we have satisfactory numerical approaches for this second inverse problem, that we will expose in the Section~\ref{sec:powenum2}.


In this section, we analyze  uniqueness and stability for the power $\alpha$ in \eqref{eq.CP} and \eqref{eq.CP.w}.  Let us start with a uniqueness result.

\begin{theorem}\label{th:uniqueness2}
Let assume that  $0<\ell \le1$, $u_0\in L^2(0,\ell)$, $u_0\neq 0$ and  $u_i$, $i=1,2$, be the solution of  \eqref{eq.CP} (or  \eqref{eq.CP.w})  corresponding to $\alpha_i$ such that for some $t_0\in (0,T]$ we have  $\partial_t u_1(x,t_0)= \partial_t u_2(x,t_0)$  and  $\partial_x u_1(x,t_0)= \partial_x u_2(x,t_0)$ for all $x\in(0,\ell)$, then 
$\alpha_1=\alpha_2$. 
\end{theorem}
%
\textbf{Proof of Theorem} -- We can assume, without loss of generality, that $\alpha_1< \alpha_2$. Setting
 \begin{equation*}
w(x,t ) = u_2(x,t) - u_1(x,t),
 \end{equation*}
we have
 \begin{equation}\label{eq:1}
\partial_t w -\partial_x(x^{\alpha_2}\partial_x w) = \partial_x((x^{\alpha_2} -x^{\alpha_1}) \partial_x u_1).
 \end{equation}
 Then, multiplying this equality by $u_1$ and integrating by parts, we obtain 
 \begin{equation}\label{eq.5}
\displaystyle\int_0^\ell(x^{\alpha_1} - x^{\alpha_2})|\partial_x u_1|^2 \, dx 
= 
\int_0^\ell (u_1\partial_t w +x^{\alpha_2}\partial_x w\, \partial_x u_1)\, dx.
 \end{equation}
By hypothesis,  we have that there exists $t_0\in (0,T]$ such that $\partial_t w(\cdot,t_0) = \partial_x w(\cdot,t_0) =0$. Then, from \eqref{eq.5} we have
\[
\displaystyle\int_0^\ell (x^{\alpha_1}-x^{\alpha_2})|\partial_x u_1(t_0,x)|^2\, dx\le 0.
\]
Since $x^{\alpha_1}>x^{\alpha_2}$ for all $x\in (0,\ell)$, the above inequality yields 
\[
\partial_x u_1(\cdot,t_0) = 0\quad \text{a.e. on } (0,\ell). 
\]
By backward uniqueness, this implies that $u_0=0$ which is excluded by the hypothesis. 
This ends the proof.  \hfill $\square$

\begin{remark} In the above uniqueness result, the observability assumption $\partial_x u_1(\cdot,t_0) = \partial_x u_2(\cdot,t_0) $ can be replaced by $u_1(\cdot,t_0) = u_2(\cdot,t_0) $. 
 
\end{remark}
\begin{theorem}[Lipschitz stability]\label{th.power}
Let us assume that  $0<\ell<1$, $0\le \alpha_1,\alpha_2 \le 2-\varepsilon$, for some $\varepsilon \in (0,2)$, and $u_0\in L^2(0,\ell)$, $u_0\neq 0$. Let $u_1$ and $u_2$ be the solutions of  \eqref{eq.CP} (or  \eqref{eq.CP.w}), corresponding to $\alpha_i (i=1,2)$, satisfying the following hypothesis: for some $t_0\in(0,T]$ and $\mu>0$
\begin{equation}\label{eq.h4}
\int_0^\ell x^{\alpha_i}|\partial_x u_i(x,t_0)|^2 \, dx > \mu. 
\end{equation}
Then there exists a constant $C=C(\varepsilon,\ell)>0$ such that 
 \begin{equation}\label{LSA}
|\alpha_2-\alpha_1| \le C \left (\int_0^\ell \Big( |\partial_t (u_1(x,t_0)- u_2(x,t_0))|^2
+ x^{\alpha_1 \vee\alpha_2 }|\partial_x (u_1(x,t_0)-u_2(x,t_0)) |^2\Big)\, dx \right)^{1/2},
 \end{equation}
 where $ \alpha_1 \vee \alpha_2=\max\{ \alpha_1,\alpha_2\}$.
\end{theorem}


\begin{remark}
The hypothesis  $0<\ell<1$ in the Theorem~\ref{th.power} is essential in the proof. 
\end{remark}

\noindent
\textbf{Proof of Theorem~\ref{th.power} --} We assume, without loss of generality, that $\alpha_1< \alpha_2$. We notice that for $0\le \alpha_1< \alpha_2$, we have $H_{\alpha_1}^1(0,\ell)  \hookrightarrow H^1_{\alpha_2} (0,\ell)$, where the spaces $H^1_{\alpha_i} (0,\ell)$ are given by  \eqref{eq.spaces}.

As before, setting 
 \begin{equation*}
w(x,t ) = u_2(x,t) - u_1(x,t),
 \end{equation*}
we have \eqref{eq:1} and \eqref{eq.5}. Now, we observe that 
  \begin{equation*}
\displaystyle\int_0^\ell(x^{\alpha_1} - x^{\alpha_2})|\partial_x u_1|^2 \, dx 
= 
\int_0^\ell (1-x^{\alpha_2-\alpha_1} )x^{\alpha_1} |\partial_x u_1|^2 \, dx 
\geq 
(1-\ell^{\alpha_2-\alpha_1}) \int_0^\ell x^{\alpha_1} |\partial_x u_1|^2 \, dx.
 \end{equation*}
On the other hand, we  have
 
 \begin{equation} \label{eq.ast}
1-\ell^{\alpha_2-\alpha_1} = \displaystyle\int_\ell^1 \dfrac{d}{d\lambda} \lambda^{\alpha_2-\alpha_1} \, d\lambda
=
(\alpha_2-\alpha_1) \displaystyle\int_\ell^1 \lambda^{\alpha_2-\alpha_1-1} \, d\lambda
\geq 
(\alpha_2-\alpha_1) \ell(1-\ell ).
 \end{equation}
Here, we have used that
 
 \begin{equation*} 
  \int_\ell^1 \lambda^{\alpha_2-\alpha_1-1}\, d\lambda \geq \ell (1-\ell),
 \end{equation*}
 since if $\alpha_2-\alpha_1-1\leq 0$, we have $\lambda^{\alpha_2-\alpha_1-1} \geq 1 > \ell$ and if 
 $0<\alpha_2-\alpha_1-1<1$, we have $\lambda^{\alpha_2-\alpha_1-1} \geq \ell^{\alpha_2-\alpha_1-1} > \ell$. 
\par\noindent
Therefore,  we conclude that 
 
  \begin{equation} \label{eq.6}
(\alpha_2-\alpha_1)\ell (1-\ell) \displaystyle\int_0^\ell x^{\alpha_1} |\partial_x u_1|^2\, dx
\leq
\displaystyle\int_0^\ell(x^{\alpha_1} - x^{\alpha_2})|\partial_x u_1|^2 \, dx.
 \end{equation}
Next, by Poincare's inequality \eqref{PI}, for the right-hand side of \eqref{eq.5} we have that  
 
\begin{equation}\label{eq.7}
\begin{array}{l}
\displaystyle\int_0^\ell (u_1\partial_t w +x^{\alpha_2}\partial_x w\, \partial_x u_1)\, dx 
\\[4mm]
\le  
\left( \displaystyle\int_0^\ell |u_1|^2\, dx\right)^{1/2} \left( \displaystyle\int_0^\ell |\partial_t w |^2\, dx\right)^{1/2} \\[4mm]
\phantom{+} +
\left( \displaystyle\int_0^\ell x^{\alpha_2}|\partial_x u_1|^2\, dx\right)^{1/2}
\left( \displaystyle\int_0^\ell x^{\alpha_2}|\partial_x w|^2\, dx\right)^{1/2}
\\[4mm]
 \le 
\sqrt{C_p}\left( \displaystyle\int_0^\ell x^{\alpha_1}|\partial_x u_1|^2\, dx\right)^{1/2} 
\left( \displaystyle\int_0^\ell |\partial_t w |^2\, dx\right)^{1/2} \
\\[4mm]
\phantom{+}  +
\ell^{(\alpha_2-\alpha_1)/2} \left(\displaystyle\int_0^\ell x^{\alpha_1} |\partial_x u_1 |^2\, dx\right)^{1/2}\left(\displaystyle\int_0^\ell x^{\alpha_2} |\partial_x w |^2\, dx\right)^{1/2}.
\end{array}
 \end{equation}
 The conclusion follows from \eqref{eq.6} and \eqref{eq.7}. This ends the proof. \hfill $\square$

\section{Reconstruction by Carleman inequalities: a general case}\label{sec.gencase}

We consider now the identification of a spatially varying diffusion 
coefficient $x^{\alpha}a(x)$ in a degenerate parabolic equation \eqref{eq.dw}:
\begin{equation}\label{1.1}
\left\{\begin{array}{ll}
\partial_t u- \partial_x(x^{\alpha} a(x) \partial_x u) = 0, & (x,t)\in(0,1)\times(0,T),  \\[1mm]
u(0,t) = 0, \quad u(1,t) = 0, &  t\in (0,T),  \\[1mm]
u(x,0) = u_0(x),  & x\in(0,1),
\end{array}\right.
\end{equation}
where we take $\ell = 1$ for simplicity.
We assume 
\[
0 < \alpha < 1,
\]
which corresponds to the weak degeneracy.  We can discuss also the case
$1<\alpha<2$ but we postpone for keeping the article concise. 

For suitably given initial value $u_0$ and coefficient 
$a$ which are stated below, there exists a unique solution to 
\eqref{1.1} by Theorem~\ref{th2.3} in Section~\ref{sec:weelpos} and by $u_{a,\alpha}$, 
we denote the solution. 

\medskip

The inverse problem that we will consider in this section is to recover
$a(x)$, $0<x<1$ and $\alpha \in (0,1)$ from the boundary measurements 
$\partial_x u_{a,\alpha}(1,t)$, provided that unknown $a$ and initial
value $u_0$ belong to some admissible sets introduced below. 

\bigskip
%

\subsection{Uniqueness result}

Henceforth, we denote by $'$ the  derivative with respect to 
one single variable, for example, $u'_0(x) = \frac{d u_0}{dx}$, 
$\eta'(t) = \frac{d \eta}{dt}$,\dots, since there is no fear of confusion. 

For statement of our main results on the uniqueness in determining
diffusion coefficients, we need to introduce admissible sets where 
unknown coefficients vary.  Moreover, we need some regularity for the
solutions to (25), which requires compatibility conditions involving
unknown coefficients and initial values.  Thus the descriptions of
the admissible sets are concerned with both coefficients and 
initial values.

First we set $W^{1,1}(0,T) := \{ a\in L^1(0,T);\, 
a'\in L^1(0,T)\}$ and 
\begin{align*}
& \mathcal{A}_0 := \{ (a, \alpha, u_0)\in C^3([0,1]) \times 
(0,1) \times C^4([0,1]);\, \,
u_0, \, (x^{\alpha}a(x)u_0')' \in W^{1,1}(0,1), \\
& (x^{\alpha}a(x))^{1/2}u_0', \, \, (x^{\alpha}a(x))^{1/2} 
(x^{\alpha}a(x)u_0')''  \in L^2(0,1), \,\,
 x^{\alpha}a(x)u_0' \in H^1(0,1), \\
& x^{\alpha}a(x)(x^{\alpha}a(x)u_0')'' \in H^1(0,1),
\end{align*}
\begin{equation}\label{26}
u_0(\xi) = 0, \quad (x^{\alpha}au_0')'(\xi) = 0 \quad 
\mbox{for $\xi=0,\,1$} \}. 
\end{equation}

\begin{remark}
We assume $(a,\alpha,u_0)\in \mathcal{A}_0$ in order to conclude
the regularity for $\www{u}_{a,\alpha}$ that we will introduce in Lemma~\ref{Lemma1}.
\end{remark}

We arbitrarily fix constants $\delta_0 > 0$, $\delta_1>0$, $\ep_0>0$, 
$\ep_1>0$, and a function $a_0 \in C^2([0,1])$ satisfying 
$a_0>0$ on $[0,1]$.
In terms of $\mathcal{A}_0$, we define an admissible set of 
coefficients $x^{\alpha}a(x)$ and initial values $u_0(x)$:
\begin{equation}\label{27}
\begin{array}{ll}
\mathcal{A}(\ep_0,\ep_1,\delta_0,\delta_1)
 := & \{ (a,\alpha, u_0) \in \mathcal{A}_0:\, 
a\ge \delta_0  \text{ on }  [0,1], \, a(x) = a_0(x)  
\mbox{ for $1-\ep_0\le x \le 1$}, 
\\[2mm]
& \vert u_0'(x)\vert \ge \delta_1  \text{ for } \ep_1 \le x \le 1\}.
\end{array}
\end{equation}
By $a=a_0$ on $[1-\ep_0,1]$, we notice that the non-degenerate 
component $a(x)$ of a 
diffusion coefficient is assumed to be given near $x=1$, and our 
concerns is the determination in an interior interval inside 
$(0,1)$.
\\

The main uniqueness result for the inverse problem of determining 
a diffusion coefficient $x^{\alpha}a(x)$ 
is the following:

\begin{theorem}\label{th.uniqueness}
Let $u_{a,\alpha}$ and 
$u_{b,\beta}$ be the solutions of \eqref{1.1} corresponding to 
$(a, \alpha,u_0)$, $(b, \beta,u_0) \in \mathcal{A}(\ep_0,\ep_1,\delta_0,
\delta_1)$, respectively.  
If  
\begin{equation}\label{28}
\ppp_xu_{a,\alpha}(1,t) = \ppp_xu_{b,\beta}(1,t), \quad 0<t<t_0, 
\end{equation}
with some $t_0>0$, then
\[
x^{\alpha}a(x) = x^{\beta}b(x), \quad \ep_1 \le x \le 1.
\]
The constant $t_0>0$ can be chosen independently of $\ep_0,\ep_1,\delta_0$
and $\delta_1$.
\end{theorem}
\bigskip 

The theorem asserts the uniqueness in determining of $x^{\alpha}a(x)$ for 
$\ep_1 \le x \le 1$ within an admissible set 
$\mathcal{A}(\ep_0,\ep_1,\delta_0, \delta_1)$ of coefficients and 
initial value.

\medskip

The characterization of the admissible set 
$\mathcal{A}(\ep_0,\ep_1,\delta_0, \delta_1)$ is not given 
separately for $x^{\alpha}a(x)$ and $u_0(x)$.  However, we can provide 
a subset of it as follows.  We arbitrarily fix $u_0 \in C^4([0,1])$ 
and positive constants $\ep_1, \ep_2, \delta_1$ such that 
\begin{equation}\label{29}
u_0(x) = 0 \quad \mbox{for $0\le x\le \ep_2$}, \qquad
\vert u_0'(x)\vert \ge \delta_1 \quad \mbox{for $\ep_1\le x\le 1$},
\quad u_0(1) = u_0''(1) = 0.
\end{equation}

Whenever $0<\alpha<1$ and a function $a$ satisfies
\begin{equation}\label{30}
\left\{ \begin{array}{rl}
& a \in C^3([0,1]), \quad a\ge \delta_0  \text{ on }  [0,1], \quad
a(x) = a_0(x) \quad \mbox{for $1-\ep_0\le x \le 1$}, \\[2mm]
& (x^{\alpha}a)'(1) = \alpha a(1) + a'(1) = 0,
\end{array}\right.
\end{equation}
we see that $(a,\alpha,u_0) \in \mathcal{A}(\ep_0,\ep_1,\delta_0,\delta_1)$.
Thus, (\ref{30}) is a description of an admissible set of unknown coefficients,
which is written independently of $u_0$.

\bigskip
\begin{coro}
Let $u_0$ satisfy (\ref{29}).
If $(a,\alpha)$ and $(b,\beta)$ satisfy (\ref{30}) and  
$$
\ppp_xu_{a,\alpha}(1,t) = \ppp_xu_{b,\beta}(1,t), \quad 0<t<t_0
$$
with some $t_0>0$, then 
\[
x^{\alpha}a(x) = x^{\beta}b(x), \quad \ep_1 \le x \le 1.
\]
\end{coro}

\bigskip

In Theorem~\ref{th.uniqueness}, we cannot identify the non-degenerate 
component $a(x)$ separately from the degeneracy power $\alpha$. 
However, under a stricter condition on an initial value $u_0$, we can 
uniquely determine both $a(x)$ and $\alpha$.  More precisely, we define 
another admissible set $\mathcal{A}_1(\ep_0,\delta_0)$ as a subset of 
$\mathcal{A}(\ep_0,\ep_1,\delta_0,\delta_1)$ by

\begin{equation}\label{31}
\begin{array}{ll}
\mathcal{A}_1(\ep_0,\delta_0)
:=  &  \{ (a,\alpha,u_0)\in \mathcal{A}_0;\, 
a\ge \delta_0  \text{ on }  [0,1], 
 \\[2mm]
& a(x) = a_0(x) \quad 
\mbox{for $1-\ep_0\le x \le 1$}, \quad
\vert u_0'(x)\vert > 0  \quad\text{for $0<x\le 1$} \}.
\end{array}
\end{equation}

We can write
$$
\mathcal{A}_1(\ep_0,\delta_0)
= \bigcap_{\ep_1>0, \delta_1>0} \mathcal{A}(\ep_0,\ep_1,\delta_0,\delta_1).
$$

Then we have

\begin{theorem}\label{th.uniqueness2}
Assume that $u_{a,\alpha}$ and $u_{b,\beta}$ be the solutions of \eqref{1.1} 
corresponding to $(a, \alpha, u_0)$, $(b, \beta, u_0) \in 
\mathcal{A}_1(\ep_0,\delta_0)$, respectively.  If 
%
$$
\ppp_xu_{a,\alpha}(1,t) = \ppp_xu_{b,\beta}(1,t), \quad 0<t<t_0, 
$$
with some $t_0 \in (0,T)$, then
\[
a(x) = b(x), \quad 0 \le x \le 1 \quad \text{and} \quad \alpha=\beta.
\]
\end{theorem}
\bigskip

\noindent
\textbf{Proof of Theorem~\ref{th.uniqueness2}.}
The proof of Theorem~\ref{th.uniqueness2} follows directly from Theorem~\ref{th.uniqueness}. Indeed, let $(a, \alpha, u_0), (b, \beta, u_0) \in 
\mathcal{A}_1(\ep_0,\delta_0)$.
Then, for any $\ep_1>0$, we can choose $\delta_1>0$ 
such that $\vert u_0(x)\vert \ge \delta_1$ for $\ep_1\le x\le 1$.
Therefore, applying Theorem~\ref{th.uniqueness}, we see that 
$x^{\alpha}a(x) = x^{\beta}b(x)$ for $\ep_1\le x\le 1$.
Here we notice that in Theorem~\ref{th.uniqueness} the constant $t_0>0$ is independent of 
$\ep_1$ and $\delta_1$.

Since $\ep_1 > 0$ is arbitrary and $a, b\in C([0,1])$, we obtain
\[
x^{\alpha}a(x) = x^{\beta}b(x), \quad 0\le x\le 1.
\]
Assume $\alpha > \beta$.  Then $x^{\alpha-\beta}a(x) = b(x)$ for 
$0 < x \le 1$.  Letting $x\downarrow 0$, we reach $b(0) = 0$.
This contradicts $b(0) \ge \delta_0 > 0$.
Similarly, $\alpha < \beta$ is impossible.  Hence $\alpha=\beta$, and so
$a(x) = b(x)$ for $0\le x\le 1$.  Thus Theorem~\ref{th.uniqueness2} is 
proved. 
\hfill $\square$.
\bigskip

Our proof of Theorem~\ref{th.uniqueness} is based on Carleman estimates, and the methodology 
requires some positivity related to the initial value, and in our case, 
$\vert u_0'(x)\vert$ must be positive in the interval $[\ep_1,1]$ where
we want to determine coefficients.  It may be practical that we formulate 
an inverse problem with the zero initial condition and 
some non-zero boundary condition.  However, in that case, our proof does not 
work.

As is described in Section~\ref{sec.rhp}, the proof of Theorem~\ref{th.uniqueness} is involved with an inversion of the integral transform which is similar to the Laplace transform.
Therefore,  the stability rate is expected to be a quite weak, although 
we can try to consider stability for our inverse problem which should be 
conditional under some a priori boundedness condition.

On the other hand, we can consider an inverse problem for the system where
$u(x,0) = u_0(x)$, $x\in (0,1)$ in \eqref{1.1} is replaced by 
$u(x,t_0) = u_0(x)$, $x\in (0,1)$ with fixed $t_0>0$.
In this case, we can prove conditional stability of Lipschitz or H\"older type
(see Imanuvilov and Yamamoto~\cite{IY0}, and Yamamoto~\cite{Ya}).
We remark that in Theorem~\ref{th.uniqueness} the inverse problem corresponds 
to  the classical initial boundary value problem.

\medskip

\subsection{Proof of Theorem~\ref{th.uniqueness}}

This section is devoted to the proof of Theorem~\ref{th.uniqueness}.
The proof is based on Carleman estimate and such a 
methodology can be found in Bukhgeim and Klibanov \cite{BK}, and then 
Bellassoued and Yamamoto \cite{BY2017}, Huang, Imanuvilov and Yamamoto 
\cite{HIY}, Imanuvilov and Yamamoto \cite {IY0}, \cite{IY1}, \cite{IY2}, 
and Klibanov \cite{Kl2}.  However, the methodology itself for a parabolic equation  needs
data $u(\cdot,t_0)$ at $t_0>0$, not $u(\cdot,0)$ for a parabolic type of 
equations.
Therefore, we first transform our inverse problem for parabolic equation
to an inverse hyperbolic problem by an integral transform in $t$ and then 
we establish the uniqueness for the transformed inverse hyperbolic problem.
Thus, the proof of Theorem~\ref{th.uniqueness} is divided into the 
following three steps.

\medskip\noindent
{\bf Section~\ref{sec.rhp}: First Step.} We transform the inverse parabolic problem to an inverse hyperbolic 
problem.

\medskip\noindent
{\bf Section~\ref{sec.ufhp}: Second Step.} We apply the method from \cite{BK} 
to a
linearized inverse hyperbolic problem and establish the uniqueness.
For a general coefficient $a(x)>0$ of the principal term, 
a key Carleman estimate can be 
proved locally only in $x$, and we have to argue starting 
at $x=1$ and we can prove the uniqueness in a small $x$-interval.  Then, we 
repeat the argument for the uniqueness until the $x$-interval of the 
uniqueness can reach $x=\ep_1$.

\medskip\noindent
{\bf Section~\ref{sec.proofth}: Thirds Step.} We translate the uniqueness for 
the 
linearized inverse hyperbolic problem to the original inverse parabolic 
problem and complete the proof of Theorem~\ref{th.uniqueness}.

\subsubsection{First Step: Reduction to hyperbolic problem}\label{sec.rhp}

We introduce an integral transform called the Reznitskaya transform, and we 
prove
Lemma~\ref{Lemma4} which connects the inverse parabolic problem to an inverse 
hyperbolic problem. For  this transform, we can refer to Romanov~\cite{Ro}. 
Notice that
Klibanov~\cite{Kl2}, pp. 594-595 have used the transform to prove the 
uniqueness for an inverse problem for non-degenerate parabolic equations.
In our article, the regularity of the solution to \eqref{1.1} requires 
careful justification in applying the transform, which is made through
Lemmas~\ref{Lemma1}-\ref{Lemma3}. 

\bigskip

Let us first consider a system of the second order in $t$. More precisely, 
let $\www{u}_{a,\alpha}
= \www{u}_{a,\alpha}(x,t)$ satisfy 
\begin{equation}\label{2.1}
\left\{ \begin{array}{ll}
 \ppp_t^2\www{u}_{a,\alpha}= \ppppa \www{u}_{a,\alpha}), &  0<x<1, \quad t>0, 
\\[2mm]
\www{u}_{a,\alpha}(0,t) = \www{u}_{a,\alpha}(1,t) = 0, & t>0, \\[2mm]
\www{u}_{a,\alpha}(x,0) = u_0(x), \quad \ppp_t\www{u}_{a,\alpha}(x,0) = 0, 
& 0<x<1.
\end{array}\right.                  
\end{equation}

In view of \eqref{29} and \eqref{30}, by Alabau-Boussouira, Cannarsa and 
Leugering \cite{FCL}, we can prove the following result. 

\begin{lemma}\label{Lemma1}
Let $(a,\alpha,u_0) \in \mathcal{A}(\ep_0,\ep_1, \delta_0, \delta_1)$,
$\mu \in (0,1)$ 
and $T>0$ be arbitrary. There exists a unique solution to \eqref{2.1} 
such that 
\begin{equation}\label{2.2}
\www{u}_{a,\alpha}\in C^3([0,\infty); H^1(\mu,1)) \cap C^2([0,\infty); 
H^2(\mu,1)) \cap C^1([0,\infty); H^3(\mu,1))                  
\end{equation}
and
\begin{equation}\label{2.3}
\ppp_t\ppp_x^2 \www{u}_{a,\alpha} \in L^{\infty}((\mu,1) \times (0,T)).
\end{equation}
Moreover, the following regularity properties hold: 
\begin{equation}\label{2.4}
\ppp_t^2\www{u}_{a,\alpha} \in L^{\infty}(0,\infty;L^2(\mu,1)), \quad
   \ppp_t\www{u}_{a,\alpha} \in L^{\infty}(0,\infty;H^1(\mu,1)), \quad 
\www{u}_{a,\alpha} \in L^{\infty}(0,\infty;H^2(\mu,1))                       
\end{equation}
and for arbitrarily given $M>0$, 
there exists a constant $C=C(M,\mu,T)>0$ such that 
\begin{align*}
\Vert \ppp_t\ppp_x^2\www{u}_{a,\alpha}\Vert_{L^{\infty}(\mu,1)\times 
(0,T))} 
& + \sum_{k,\ell=0}^1 \Vert \ppp_t^k\ppp_{x,t}^{\ell}\www{u}_{a,\alpha}\Vert
_{C([0,T];H^2(\mu,1))}
\\
&+ \sum_{k=0}^1 \Vert \ppp_t^k\www{u}_{a,\alpha}\Vert
_{H^2((\mu,1)\times (0,T))}
\le C(M,\mu,T)
\end{align*}
if $\Vert a\Vert_{C^2([0,1])} + \Vert u_0\Vert_{C^4([0,1])}
\le M$.
\end{lemma}
\bigskip

\noindent
\textbf{Proof of Lemma~\ref{Lemma1}.}
The proof is an application of a usual lift-up argument of the regularity
to the result given in \cite{FCL}. 
We set
\[
A_{a,\alpha}u_0(x) := (x^{\alpha}a(x)u_0')' = \frac{d}{dx}\left(x^{\alpha}a(x)
\frac{du_0}{dx}(x) \right).
\]
Then, $A_{a,\alpha}u_0$ satisfies \eqref{26} and we can write 
\[\left\{ \begin{array}{l}
A_{a,\alpha}u_0 \in W^{1,1}(0,T), \\[2mm]
x^{\alpha/2}\sqrt{a(x)}(A_{a,\alpha}u_0)'
= (x^{\alpha}a(x))^{1/2}(x^{\alpha}a(x)u_0')'' \in L^2(0,1), \\[2mm]
A_{a,\alpha}u_0(0) = A_{a,\alpha}u_0(1) = 0,\\[2mm]
x^{\alpha}a(x)(x^{\alpha}a(x)u_0')'' \in H^1(0,1).
\end{array}\right.
\]
Let us consider the following hyperbolic problem with $A_{a,\alpha}u_0$ as an initial condition: 
\[
\left\{ \begin{array}{ll}
\ppp_t^2v = \ppppa v), &  0<x<1, \, t>0, \\[2mm]
v(0,t) = v(1,t) = 0, &  t>0, \\[2mm]
v(x,0) = A_{a,\alpha}u_0(x), \quad \ppp_tv(x,0) =  0, & 0<x<1.
\end{array}\right.
\]
By \cite{FCL}, there exists a unique solution $v$ to this problem  
which verifies 
\[
v \in C^2([0,\infty);L^2(0,1))
\cap C^1([0,\infty);H^1(\mu,1))\cap C^0([0,\infty);H^2(\mu,1))
\]
for arbitrary $\mu\in(0,1)$.

We can verify that $A_{a,\alpha}u$ satisfies \eqref{1.1} in the sense of
distributions.  The uniqueness of solution to \eqref{1.1} yields 
\[
v = A_{a,\alpha}u \in C^1([0,\infty);H^1(\mu, 1)),
\]
and so
\[
u \in C^1([0,\infty);H^3(\mu, 1)).
\]
Thus the proof of Lemma~\ref{Lemma1} is complete. \hfill
$\square$
\bigskip


Next, we introduce an integral transform called the Reznitskaya transform 
(see~\cite{Ro}, pp. 213--215) as follows:
\begin{equation} 
\label{2.5BIS}
(K\eta)(t) := \int^{\infty}_0 \eta(\tau)G(t,\tau) d\tau, \quad
\eta \in L^{\infty}(0,\infty),  
\end{equation}
where
\[
G(t,\tau) = \frac{1}{\sqrt{\pi t}} e^{-\frac{\tau^2}{4t}}, \quad t, \tau > 0.
\]

We note that 
\begin{equation}\label{2.5}
\ppp_{\tau}^kG(t,\cdot) \in L^{\infty}(0,\infty) \cap
L^1(0,\infty), \quad 
\ppp_t^kG(\cdot,\tau) \in L^{\infty}(0,\infty)\quad 
\mbox{for all $t, \tau>0$ and $k\in \N$.}
\end{equation}

Then, the following result can be proved. 
\begin{lemma}\label{Lemma2}
Let $\eta \in W^{2,\infty}(0,\infty)$ satisfy $\eta'(0) = 0$.  Then
\[
\frac{d}{dt}(K\eta)(t) = \left(K \frac{d^2\eta}{d\tau^2}\right)(t), 
\quad t>0.
\]
\end{lemma}
\bigskip
\textbf{Proof of Lemma~\ref{Lemma2}.} The proof is given in~\cite{Ro}, 
pp. 213--214, and for completeness we will here present the proof. 
Since $\ppp_tG(t,\tau) = \ppp_{\tau}^2G(t,\tau)$ for $t>0$ and $\tau>0$, 
and $\eta \in W^{2,\infty}(0,\infty)$, we have
\begin{align*}
 \frac{d}{dt}(K\eta)(t) & = \int^{\infty}_0 \eta(\tau)\ppp_tG(t,\tau) d\tau
= \int^{\infty}_0 \eta(\tau)\ppp_{\tau}^2G(t,\tau) d\tau
\\[2mm]
& = \Big[ \eta(\tau)\ppp_{\tau}G(t,\tau) - \eta'(\tau)G(t,\tau)\Big]^{\tau=\infty}
_{\tau=0} + \int^{\infty}_0 \eta''(\tau)G(t,\tau) d\tau, \quad t>0.
\end{align*}
Using that $\eta'(0) = 0$ and $\ppp_{\tau}G(t,0) = 0$ for $t>0$, we obtain
\[
\frac{d}{dt}(K\eta)(t) = \int^{\infty}_0 \eta''(\tau)G(t,\tau) d\tau
= \left(K\frac{d^2\eta}{dt^2}\right)(t), \quad t>0.
\]
The proof of Lemma~\ref{Lemma2} is complete.\hfill $\square$

\bigskip

We will also need the following result. 
\begin{lemma}\label{Lemma3}

(i) Let $\www{u}_{a,\alpha}$ satisfy \eqref{2.1} with 
$(a,\alpha,u_0) \in \mathcal{A}(\ep_0,\ep_1,\delta_0,\delta_1)$, 
and let us set 
\begin{equation}\label{38}
V_{a,\alpha}(x,t) := (K\www{u}_{a,\alpha}(x,\cdot))(t)
= \int^{\infty}_0 \www{u}_{a,\alpha}(x,\tau)G(t,\tau) d\tau, 
\quad 0<x<1, \, t>0.
\end{equation}
Then 
\begin{equation}\label{2.6} 
\ppp_tV_{a,\alpha}(x,t) = \ppppa V_{a,\alpha}(x,t)) 
\quad \mbox{for almost all $x \in (0,1)$ 
and $t>0$}                                   
\end{equation}
and
\begin{equation}\label{2.7} 
\left\{ \begin{array}{l}
V_{a,\alpha}(0,t) = V_{a,\alpha}(1,t) = 0, \quad t>0, \\[2mm]
V_{a,\alpha}(x,0) = u_0(x), \quad \mbox{for almost all $x \in (0,1)$}. 
\end{array}\right.
\end{equation}
(ii) Let $u_{a,\alpha}$ be the solution to \eqref{1.1} where 
$(a,\alpha,u_0) \in \mathcal{A}(\ep_0,\ep_1,\delta_0,\delta_1)$. 
Then 
\begin{equation}\label{2.8} 
u_{a,\alpha}(x,t) = \int^{\infty}_0 \www{u}_{a,\alpha}(x,\tau)G(t,\tau) d\tau, 
\quad 0<x<1, \, t>0.                
\end{equation}
\end{lemma}

\bigskip\noindent
\textbf{Proof of Lemma~\ref{Lemma3}.} 
Part \textit{(ii)} follows readily from part \textit{(i)} by the uniqueness of 
the 
solution $u_{a,\alpha}$ (see Theorem~\ref{th2.3} in Section~\ref{sec:weelpos}).
Therefore, it suffices to prove part \textit{(i)}. 
Henceforth, we set 
\[
(g,h) := \int^1_{\ep_2} g(x)h(x)\, dx,
\]
with arbitrary $\ep_2>0$.
For $\theta \in C^{\infty}_0(\ep_2,1)$, we set 
\[
\www{w}(t) = (\www{u}_{a,\alpha}(\cdot,t), \theta) 
= \displaystyle\int^1_{\ep_2} \www{u}_{a,\alpha}(x,t)\theta(x)\, dx 
\quad \text{and} \quad   w(t) = (V_{a,\alpha}(\cdot,t), \theta) \quad 
\text{for } t>0.
\]
Then, the initial condition in \eqref{2.1} and the regularity \eqref{2.4} 
yield
\[
\www{w}\in W^{2,\infty}(0,\infty), \quad \ppp_t\www{w}(0) = 0
\]
and
\[
w(t) = (V_{a,\alpha}(\cdot,t), \theta) 
= \int_{\varepsilon_2}^1 V_{a,\alpha}(x,t)\theta(x)\,dx
= \int^{\infty}_0 \www{w}(\tau)G(t,\tau) d\tau
= (K\www{w})(t) , \quad t>0.
\]
Therefore, Lemma~\ref{Lemma2} implies
\[
\frac{dw}{dt}(t) = \frac{d}{dt}(K\www{w}) (t) 
= 
\left(K\frac{d^2\www{w} }{d\tau^2}\right) (t)
=  \int^{\infty}_0 \frac{d^2}{d\tau^2}\www{w}(\tau)
G(t,\tau) d\tau, 
\]
that is,
\begin{align*}
\frac{d}{dt}(V_{a,\alpha}(\cdot,t), \theta) 
& = \int^{\infty}_0 \frac{d^2}{d\tau^2}(\www{u}_{a,\alpha}(\cdot,\tau), 
\theta) G(t,\tau) d\tau
= \int^{\infty}_0 \left( \int^1_{\ep_2} \frac{\ppp^2 
\www{u}_{a,\alpha}(x,\tau)}
{\ppp\tau^2}\theta(x) dx \right) G(t,\tau) \,d\tau    
\\[2mm]
& = \int^1_{\ep_2} \left( \int^{\infty}_0 \frac{\ppp^2\www{u}_{a,\alpha}}
{\ppp \tau^2}(x,\tau) G(t,\tau) d\tau\right) \theta(x)\,  dx
\\[2mm]
& = \left( \left(K \frac{\ppp^2\www{u}_{a,\alpha}}{\ppp \tau^2}\right)
(\cdot,t), \, \theta\right).
\end{align*}
Therefore,
\begin{equation}\label{2.9} 
\left( \frac{\ppp V_{a,\alpha}}{\ppp t}(\cdot,t), \theta\right)
= \left( \left(K\frac{\ppp^2\www{u}_{a,\alpha}}{\ppp \tau^2}\right), 
\, \theta\right)
\end{equation}
for all $\theta \in C^{\infty}_0(\ep_2,1)$.

Recall that 
\[
A_{a,\alpha}w(x,t) = \ppp_x( x^{\alpha}a(x) \ppp_xw(x,t)), \quad 
0<x<1, \, t>0.
\]
By \eqref{2.2}, \eqref{2.4} and \eqref{2.5}, we can justify the 
exchange of the operation by $A_{a,\alpha}$ and $\int^\infty_0 \cdots
d\tau$ in (\ref{38}),
and obtain
\[
A_{a,\alpha}V_{a,\alpha}(x,t) = \int^{\infty}_0 (A_{a,\alpha}
\www{u}_{a,\alpha})(x,\tau) G(t,\tau) \, d\tau
= K(A_{a,\alpha}\www{u}_{a,\alpha})(x,t), \quad \ep_2 < x < 1, \, t>0.
\]
Hence, using \eqref{2.9} and since $\www{u}_{a,\alpha}$ is the solution of 
\eqref{2.1}, we have 
\[
\left( \frac{\ppp V_{a,\alpha}}{\ppp t}(\cdot,t) 
- A_{a,\alpha}V_{a,\alpha}(\cdot, t), \, \theta\right)
= \left( K\left( \frac{\ppp^2\www{u}_{a,\alpha}}{\ppp \tau^2} 
- A_{a,\alpha}\www{u}_{a,\alpha}\right)
(\cdot,t), \, \theta\right) = 0
\]
for all $\theta \in C^{\infty}_0(\ep_2,1)$ and $t>0$.  
Consequently, since $\ep_2>0$ is  arbitrary, we reach 
\[
\frac{\ppp V_{a,\alpha}}{\ppp t}(x,t) = \frac{\ppp}{\ppp x}\left( 
x^{\alpha}a(x)\frac{\ppp V_{a,\alpha}}{\ppp x}(x,t) \right)
\]
for almost all $x \in (0,1)$ and $t>0$.  Therefore, \eqref{2.6} holds.

Using $\www{u}_{a,\alpha}(0,t) = \www{u}_{a,\alpha}(1,t) = 0$ for $t>0$, 
we obtain that the first condition in 
\eqref{2.7} is directly seen.  

By the change of variables $\xi = \frac{\tau}{2\sqrt{t}}$, we have 
\begin{equation}\label{431}
V_{a,\alpha}(x,t) = \int^{\infty}_0 \www{u}_{a,\alpha}(x,\tau)G(t,\tau) d\tau 
= \frac{2}{\sqrt{\pi}} \int^{\infty}_0 \www{u}_{a,\alpha}(x,2\xi\sqrt{t})
e^{-\xi^2} d\xi, \quad t>0.
\end{equation}
Since $\www{u}_{a,\alpha} \in C([0,\infty); C([\ep_2,1]))$, and using that 
$\ppp_t\www{u}_{a,\alpha} \in L^{\infty}(0,\infty;H^1(\ep_2,1))$ in 
\eqref{2.4}, we have
\[
\lim_{t\to 0} \www{u}_{a,\alpha}(x,2\xi\sqrt{t}) = \www{u}_{a,\alpha}(x,0)
\]
for almost all $x \in (\ep_2,1)$ and $\xi>0$, and  
by $\www{u}_{a,\alpha} \in L^{\infty}(0,\infty;H^1(\ep_2,1))$ in \eqref{2.4}, 
we apply the Lebesgue convergence theorem to \eqref{431} and obtain
\[
V_{a,\alpha}(x,0) = \frac{2}{\sqrt{\pi}} \int^{\infty}_0 
\www{u}_{a,\alpha}(x,0)e^{-\xi^2} d\xi 
= \www{u}_{a,\alpha}(x,0)\frac{2}{\sqrt{\pi}} \int^{\infty}_0 
e^{-\xi^2} d\xi = u_0(x),
\]
for almost all $x \in (\ep_2,1)$, which is the second condition 
in \eqref{2.7}. Since $\ep_2>0$ is arbitrary, the proof of 
Lemma~\ref{Lemma3} is complete. \hfill $\square$


\bigskip

Moreover, we prove
\begin{lemma}\label{Lemma4} 
Let us assume that $\ppp_xu_{a,\alpha}(1,t) = \ppp_xu_{b,\beta}(1,t)$
 for $0 < t < t_0$.  Then, 
\[
\ppp_x\www{u}_{a,\alpha}(1,t) = \ppp_x\www{u}_{b,\beta}(1,t), \quad 0<t<\infty.
\]
\end{lemma}

\bigskip
\noindent
\textbf{Proof of Lemma~\ref{Lemma4}.} 
Since $u_{a,\alpha}(x,\cdot)$ is analytic in $t>0$ with fixed $x \in [\ep_2,1]$,
by Theorem~\ref{th2.3} (ii) in Section~\ref{sec:weelpos}, we see that  $\ppp_xu_{a,\alpha}(1,t) 
= \ppp_xu_{b,\beta}(1,t)$ for all $t>0$.  Therefore, by \eqref{2.8} 
we have
\begin{equation}\label{2.10}
\int^{\infty}_0 (\ppp_x\www{u}_{a,\alpha} - \ppp_x\www{u}_{b,\beta})(1,\tau)
G(t,\tau) d\tau 
= 0, \quad t>0.                                     
\end{equation}
On the other hand, in the definition \eqref{2.5BIS} of $K$, we change the 
variables in the expression of $K$ by $r := \frac{1}{4t}$ and $\tau := \sqrt{z}$, and we 
see 
\[
(K\eta)\left( \frac{1}{4r}\right) 
= \sqrt{\frac{r}{\pi}} \int^{\infty}_0 e^{-rz}\frac{\eta(\sqrt{z})}{\sqrt{z}}
\,dz, \quad r>0,
\]
for $\eta \in L^{\infty}(0,\infty)$.

Thus,  the uniqueness of the Laplace transform implies that $\eta \equiv 0$ 
if $(K\eta)(t) = 0$ for $t>0$ and $\eta \in L^{\infty}(0,\infty)$.
Since $(\ppp_x\www{u}_{a,\alpha} - \ppp_x\www{u}_{b,\beta})(1,\cdot) \in L^{\infty}(0,\infty)$ 
by \eqref{2.4}, it follows from \eqref{2.10} that 
$\ppp_x\www{u}_{a,\alpha}(1,t) = \ppp_x\www{u}_{b,\beta}(1,t)$ for $t>0$.
Therefore,  the proof of Lemma~\ref{Lemma4} is complete.\hfill $\square$ 
%
%

\subsubsection{Second Step: Uniqueness for the linearized inverse hyperbolic 
problem}\label{sec.ufhp}

Here we establish the uniqueness result for the following 
hyperbolic inverse source problem.  We recall that a 
sufficiently small constant $\ep_1>0$ is given.  Moreover 
let $p_0>0$ be arbitrarily given.

Given $p\in C^2([\ep_1,1])$ with $p \ge p_0 > 0$ on $[\ep_1,1]$, 
let $y=y(x,t)$ satisfy 
\begin{equation}\label{3.1} 
\left\{\begin{array}{ll}
 \ppp_t^2y= \ppppp y) + \ppp_x(f(x)\ppp_xR(x,t)), & \ep_1<x<1, \, t>0, 
\\[2mm]
y(1,t) = 0, & t>0, 
\\[2mm]
y(x,0) = \ppp_ty(x,0) = 0, & \ep_1 < x < 1,
\end{array}\right.                               
\end{equation}
where $R(x,t) $ is a given function.

In what follows, we show the uniqueness for an inverse problem of determining of $f(x)$ 
for $\ep_1<x<1$ by boundary data 
$\ppp_xy(1,t)$ for all $t>0$. 

\begin{prop}\label{Proposition1} 
Let $p\in C^2([\ep_1,1])$ with $p>0$ on $[\ep_1,1]$ and let 
$f \in H^1(\ep_1,1)$ satisfy $f(1) = 0$.  We further assume that there 
exists a constant $r_0>0$ such that   
\begin{equation}\label{eq3.2}
\left\{ \begin{array}{l}
 \ppp_tR(x,0) = 0, \quad \vert \ppp_xR(x,0)\vert 
\ge r_0 > 0, \quad \ep_1\le x \le 1,     \\[2mm]
\ppp_t\ppp_xR, \, \ppp_t\ppp_x^2R \in L_{loc}^{\infty}
((\ep_1,1)\times (0,\infty)),                \\[2mm]
\ppp_x^2R(\cdot,0) \in L^{\infty}(\ep_1,1).
\end{array}\right.                          
\end{equation}
If $y \in H^2_{loc}((\ep_1,1) \times (0,\infty))$ satisfies \eqref{3.1},
$\ppp_ty \in H^2_{loc}((\ep_1,1) \times (0,\infty))$ and
\[
\ppp_xy(1,t) = 0, \quad t>0,
\]
then $f(x) = 0$ for $\ep_1 < x < 1$.
\end{prop}

Notice that in this proposition, as is seen by the proof, we need not assume 
$\ppp_xy(1,t) = 0$ 
for all $t>0$, but 
$\ppp_xy(1,t) = 0$ for $0<t<T_*$ with some $T_*>0$ which is determined
by $p(x)$, is sufficient.  However for our puropose
we can assume $\ppp_xy(1,t) = 0$ for all $t>0$.  

The rest of this subsection is devoted to the poof of Proposition~\ref{Proposition1}.

\bigskip\noindent
\textbf{Proof of Proposition~\ref{Proposition1}.}
The methodology for the proof is based on Bukhgeim and Klibanov method from~\cite{BK}, and  also on the arguments by Imanuvilov and Yamamoto from~\cite{IY1} and~\cite{IY2}. 
Recently, Huang, Imanuvilov and Yamamoto in~\cite{HIY} have
simplified those arguments and omits the cut-off procedures.
Here we apply such argument form \cite{HIY}. The key idea is that we take a suitable weight function that originally takes smaller values on the boundary of a domain in $(x,t)$ where data are not given, so that the weight function can well control such unknown data and therefore the cut-off function is not necessary.

The proof will be divided into three steps. 

\medskip 
\noindent 
\textbf{First step of Proposition~\ref{Proposition1}: Carleman estimates.}
 
Henceforth, with arbitrarily chosen constants $M>0$ and $p_0>0$, we assume
\[
p\in C^2([\ep_1,1]), \quad p\ge p_0 > 0 \quad \text{on 
$[\ep_1,1]$}, \quad \Vert p\Vert_{C^2([\ep_1,1])} \le M.
\]
Furthermore, we fix a sufficiently small constant $\gamma>0$ and 
a sufficiently large constant $\la>0$.  
Here and henceforth, for $M$ and $r_0$,
we choose $\theta\in (0,1)$ sufficiently close to $1$
and sufficiently small $\delta>0$ such that 
\[
\ep_1 < \theta + \delta < 1
\]
is always satisfied.

\begin{figure}[h!]
\centering
\includegraphics[width=7cm]{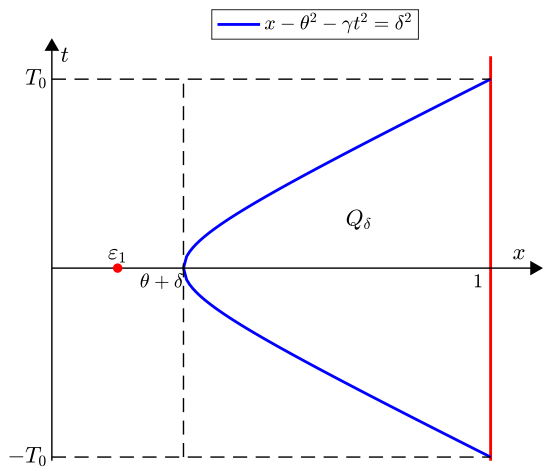}
\caption{Set $Q_\delta$.}
\label{Qdelta}
\end{figure}

We set
\[ 
\psi(x,t) := (x-\theta)^2 - \gamma t^2, \quad
\va(x,t) = e^{\la\psi(x,t)}
\]
and 
\[
Q_\delta := \{(x,t);\, -\infty<t<\infty, \quad \theta+\delta < x < 1, \quad 
\psi(x,t) > \delta^2\}, \quad
T_0:= \sqrt{\frac{(1-\theta)^2 - \delta^2}{\gamma}}.
\]
Then, we note that $Q_\delta$ is bounded by a hyperbolic curve 
$(x-\theta)^2 - \gamma t^2=\delta^2$ and $x=1$, and 
\begin{equation}\label{43}
\mbox{$(x,t)\in Q_\delta$ implies 
$\theta+\delta <x<1$ and $\vert t\vert < T_0, \quad 
Q_\delta \cap \{t=0\} = (\theta+\delta, \,1)$}, 
\end{equation}
see Figure~\ref{Qdelta}.
%

Henceforth $C, C_0>0$, etc. will denote generic constants which are
independent of the parameter $s>0$, but can be dependent on other parameters such as $\la>0$.  
Since other parameters  such as $\la$ are fixed, we do not specify the dependence of 
$C, C_0$, etc. on them.

We show the following result:

\begin{lemma}[Hyperbolic Carleman estimate]\label{Lemma5}
There exist constants $C>0$ and $s_0>0$ depending on $M$ and $r_0$, 
but independent of $s$ and the choices of $p$, such that
\begin{align*}
& \int_{Q_\delta} (s\vert \ppp_xy\vert^2+ s\vert \ppp_ty\vert^2
+ s^3\vert y\vert^2) \weight \,dxdt 
\\[2mm]
& \le C\int_{Q_\delta} \vert F(x,t)\vert^2 \weight \, dxdt
+ C\int_{\ppp Q_\delta} (s\vert \ppp_xy\vert^2 
+ s\vert \ppp_ty\vert^2 + s^3\vert y\vert^2) e^{2s\va} dS 
\end{align*}
for all $s\ge s_0$ and $y \in H^2(Q_\delta)$ satisfying
$$
\ppp_t^2y = \ppppp y) + F \quad \text{in }  Q_\delta
$$
with some $F \in L^2(Q_\delta)$.

\medskip
In particular, if $y(1,t) = \ppp_xy(1,t) = 0$ for $-T_0<t<T_0$, then 
\begin{equation}\label{44}
\int_{Q_\delta} (s\vert \ppp_xy\vert^2+ s\vert \ppp_ty\vert^2
+ s^3\vert y\vert^2) \weight \,dxdt 
\le C\int_{Q_\delta} \vert F(x,t)\vert^2 \weight \, dxdt
+ Cs^3\exex\Vert y\Vert^2_{H^2(Q_\delta)}
\end{equation}
for all $s\ge s_0$ and $y \in H^2(Q_\delta)$.
\end{lemma}

\noindent 
\textbf{Proof of Lemma~\ref{Lemma5}.} This is a classical Carleman estimate:
the proof can be done similarly to Theorem~\ref{th.uniqueness} in Bellassoued and 
Yamamoto \cite{BK} or Imanuvilov \cite{Im}.  Here we omit the details of the 
proof.   
Estimate \eqref{44} can follow because $\psi(x,t) = \delta^2$ for 
$(x,t) \in \ppp Q_\delta \setminus \{ x=1\}$ and 
\begin{align*}
 \Vert \ppp_xy\Vert_{L^2(\ppp Q_\delta)} 
& + \Vert \ppp_ty\Vert_{L^2(\ppp Q_\delta)}
+ \Vert y\Vert_{L^2(\ppp Q_\delta)}
\\[2mm]
& \le  C(\Vert \ppp_xy\Vert_{H^1(Q_\delta)} + \Vert \ppp_ty\Vert_{H^1(Q_\delta)}
+ \Vert y\Vert_{H^1(Q_\delta)})
\le C\Vert y\Vert_{H^2(Q_\delta)}
\end{align*}
by the trace theorem. \hfill $\square$

\begin{remark}\label{remark4.9}
For Carleman estimates for hyperbolic equations, we usually have to assume a 
certain extra condition for the principal part
$p(x)$ e.g., about the wave speed, and the extra condition is restrictive 
in higher spatial dimensions.  However, the one spatial dimension case is 
exceptional and does not require any extra condition thanks to 
the smallness of the domain $Q_\delta$.  
More precisely, in the one dimensional 
case, such an extra condition is described as
\[
\left\vert \frac{p'(x)}{2p(x)}(x-\theta) \right\vert < 1, 
\quad \theta+\delta \le x \le 1
\]
(see e.g., Proposition 2.1 in Imanuvilov and Yamamoto \cite{IY2}).
Since $\displaystyle \max_{\theta+\delta \le x\le 1} \vert x-\theta\vert 
= 1-\theta$ is sufficiently small, we see that 
\begin{equation}\label{441}
\left\vert \frac{p'(x)}{2p(x)}(x-\theta) \right\vert 
\le \frac{M}{2r_0}(1-\theta) < 1.
\end{equation}
Hence, there exists a constant $\mu_*(r_0,M)>0$ such that 
if $0<1-\theta < \mu_*(r_0,M)$, then \eqref{441} and so the 
Carleman estimate \eqref{44} hold.  
\end{remark}

In view of \eqref{441},
the key Carleman \eqref{44} estimate holds near $x=1$ where Cauchy data are 
given, that is, for $\theta+\delta < x < 1$.  Thus, we first apply the Carleman estimate only for 
$\theta + \delta < x < 1$ to prove that $f=0$ in $(\theta+\delta, \, 1)$.
This will be done in second step of the proof of Proposition~\ref{Proposition1}.  Next, we repeat the same argument by 
replacing $x=1$ by $x=\theta$, and this will be done in third step of the proof of Proposition~\ref{Proposition1}.

\bigskip
 
Before proceeding to second step, 
we prove a Carleman estimate for a simple operator 
$\frac{d}{dx}$.  We set 
\[
\va_0(x) = e^{\la (x-\theta)^2}
\]
with fixed large constant $\la>0$.
\begin{lemma}\label{Lemma6}
There exist constants $s_0>0$ and $C>0$ such that 
\[
s\int^1_{\theta+\delta} \vert f(x)\vert^2 e^{2s\va_0(x)} \,dx 
\le C\int_{\theta+\delta}^1 \vert f'\vert^2 e^{2s\va_0} \,dx,
\]
for all $s\ge s_0$ and all $f\in H^1(\theta+\delta,1)$ satisfying
$f(1) = 0$.
\end{lemma}

\medskip
\noindent 
\textbf{Proof of Lemma~\ref{Lemma6}.} 
The proof relies on the estimation of $h:= fe^{s\va_0}$
and the integration by parts, which is a conventional proof of 
the Carleman estimate.

We set 
\[
h(x):= f(x)e^{s\va_0(x)}, \quad 
Lh(x) := e^{s\va_0(x)}\frac{d}{dx}(e^{-s\va_0(x)}h(x))
= h'(x) - s\va_0'(x)h(x).
\]
Then, integrating by parts and taking into account that $f(1)=0$, we have that there exists a positive constant $C$ such that 
\begin{align*}
 \int^1_{\theta+\delta} \vert f'\vert^2 e^{2s\va_0}\,dx
& = \int^1_{\theta+\delta} \vert Lh\vert^2 dx 
= \int^1_{\theta+\delta} \vert h'(x) - s\va_0'(x)h(x)\vert^2 \,dx
\\[2mm]
& \ge -2s\int^1_{\theta+\delta} \va_0'(x)h(x)h'(x)\,dx
= -s\int^1_{\theta+\delta} \va_0'(x)\frac{d}{dx}h^2(x)\, dx
\\[2mm]
& =\Big [s\va_0'(x)h^2(x)\Big]^{x=\theta+\delta}_{x=1}
+ s\int^1_{\theta+\delta} \va_0''(x)h^2(x) \,dx
\\[2mm]
& = 2\la s\delta e^{\la\delta^2}h^2(\theta+\delta)
+ s\int^1_{\theta+\delta} (2\la\va_0(x)+4\la^2(x-\theta)^2\va_0(x))h^2(x) \,dx
\\[2mm]
& \ge s \int^1_{\theta+\delta} 2\la e^{\la(x-\theta)^2}h^2(x) \,dx 
\ge C s\int^1_{\theta+\delta} \vert f\vert^2 e^{2s\va_0} \,dx.
\end{align*}
Thus, the proof of Lemma~\ref{Lemma6} is complete. \hfill $\square$

\bigskip 
\noindent 
\textbf{Second step of the proof of Proposition~\ref{Proposition1}. }
 
We take the even extension in $t$ of the solution $y$ of \eqref{3.1}  
to $(-\infty,0)$, and we denote it by the same letter:
\[
y(x,t) =
\left\{ \begin{array}{ll}
 y(x,t),  & \quad t\ge 0, \\
y(x,-t), &  \quad t < 0.
\end{array}\right.
\]
Also, for $R$, we make the even extension in $t<0$.  Then, using that
$y(\cdot,0) = \ppp_ty(\cdot,0) = 0$ and $\ppp_tR(\cdot,0) = 0$ 
on $[\ep_1,1]$ by \eqref{eq3.2},  
we can verify that 
\begin{equation}\label{45} 
\left\{ \begin{array}{l}
\ppp_ty,\, y \in H^2_{loc}((\ep_1,1) \times (-\infty, \infty)),\\[2mm]
\ppp_tR,\, R \in L^2_{loc}(-\infty,\infty; H^2(\ep_1,1)),           \\[2mm]
\ppp_t\ppp_xR, \, \ppp_t\ppp_x^2R 
\in L^{\infty}_{loc}((\ep_1,1) \times (-\infty, \infty))
\end{array}\right.
\end{equation}
and
\begin{equation}\label{46}
\left\{\begin{array}{ll}
\ppp_t^2y= \ppppp y) + \ppp_x(f(x)\ppp_xR(x,t)),  & 
\ep_1<x<1, \, -\infty < t < \infty, 
\\[2mm]
y(1,t) = \ppp_xy(1,t) = 0, &  -\infty < t < \infty, 
\\[2mm]
y(x,0) = \ppp_ty(x,0) = 0, &  \ep_1 < x < 1.
\end{array}\right.                                   
\end{equation}

In \eqref{46}, we see that $y, \ppp_ty \in H^2((\ep_1,1)\times
(-T,T))$ with arbitrarily given $T>0$.  Setting 
$$
z := \ppp_ty \in H^2((\ep_1,1) \times (-T,T)),
$$
and using \eqref{46}, we have
\begin{equation}\label{47}
\left\{\begin{array}{ll}
\ppp_t^2z= \ppppp z) + \ppp_x(f(x)\ppp_t\ppp_xR(x,t)),  & 
\ep_1<x<1, \, -\infty < t < \infty, 
\\[2mm]
z(1,t) = \ppp_xz(1,t) = 0, &  -\infty < t < \infty, 
\\[2mm]
z(x,0) = 0, \quad \ppp_tz(x,0) = \ppp_x(f(x)\ppp_xR(x,0)), &  \ep_1 < x < 1.
\end{array}\right.                                   
\end{equation}
Notice that the second and the third conditions in \eqref{47} are obtained by 
\eqref{46} taking into account that $\ppp_t^2y \in H^1(-T,T;L^2(\ep_1,1))$ and the Sobolev embedding such as 
$H^1(-T,T;L^2(\ep_1,1)) \subset C^0([-T,T];L^2(\ep_1,1))$.
Here, we also note that 
$\ppp_t\ppp_x^2y, \, \ppp_t\ppp_x^2R \in L^2(-T,T;L^2(\ep_1,1))$, and so
$\ppp_x^2y, \, \ppp_x^2R \in C^0([-T,T];L^2(\ep_1,1))$ by the 
Sobolev embedding.

Setting now
\[
M_0:= \Vert z\Vert_{H^2((\ep_1,1)\times (-T,T))} = \Vert \partial_t y\Vert_{H^2((\ep_1,1)\times (-T,T))}, 
\]
by the regularity assumption \eqref{45} on $y$, we see that $M_0 < \infty$.

Now, we apply the Carleman estimate \eqref{44} to \eqref{47} and we obtain that there exists 
a constant $C>0$, which is independent of $s>0$, such that 
\begin{align}\label{48}
\int_{Q_\delta} (s\vert \ppp_xz\vert^2  & + s\vert \ppp_tz\vert^2
+ s^3\vert z\vert^2) e^{2s\va}\, dxdt \\[2mm]
& 
\le C\int_{Q_\delta} |\ppp_x(f(x)\ppp_x\ppp_tR(x,t)\vert^2 e^{2s\va}\, dxdt
+ Cs^3\exex M_0^2, \nonumber
\end{align}
for all large $s>0$.

We set
\[
w := (\ppp_ty) e^{s\va} = ze^{s\va}. 
\]
Then, direct calculations yield
\[
\ppp_t^2 w  = (\ppp_t^3y) \wwww
+ 2s(\ppp_t^2y)(\ppp_t\va) \wwww
+ (s\ppp_t^2\va + s^2(\ppp_t\va)^2)) \wwww \ppp_ty.
\]
On the other hand, 
\[
p\ppp_x^2w  = p(\ppp_x^2\ppp_ty) \wwww
+ 2ps(\ppp_t\ppp_xy)(\ppp_x\va)\wwww
+ p(s(\ppp_x^2\va) + s^2(\ppp_x\va)^2)) \wwww \ppp_ty,
\]
and also
\[
p'\ppp_xw = p'(\ppp_x\ppp_ty)\wwww + sp'(\ppp_x\va)(\ppp_ty)\wwww.
\]
Therefore, taking the derivative with respect to $t$ in \eqref{3.1} and multiplying by $e^{s\varphi}$, we can write 
\begin{equation}\label{eq.w}
\ppp_t^2w - \ppppp w)
= \wwww \ppp_x(f(x)\ppp_x\ppp_tR(x,t)) + P_s(w),
\end{equation}
where the term $P_s(w)$ satisfies  
\begin{equation}\label{50}
\begin{array}{ll}
\vert P_s(w) \vert  & \le Cs(\vert \ppp_t^2y\vert
+ \vert \ppp_ty\vert + \vert \ppp_x\ppp_ty\vert) e^{s\va}
+ Cs^2\vert \ppp_ty\vert e^{s\va} 
\\[2mm]
 & \le C(s^2\vert w\vert + s\vert \ppp_tw\vert + s\vert \ppp_xw\vert),
\qquad (x,t) \in Q_\delta,
\end{array}
\end{equation}
for all large $s>0$.
For the last inequality, we have used that $w = (\ppp_ty) \wwww$ and 
\[
(\ppp_t^2y) \wwww = \ppp_tw - s(\ppp_t\va)(\ppp_ty)\wwww
= \ppp_tw - s(\ppp_t\va)w,
\]
and
\[
(\ppp_x\ppp_ty)\wwww = \ppp_xw - s(\ppp_x\va)w,
\]
which imply
\[
\vert \ppp_t^2y\vert \wwww \le \vert \ppp_tw\vert 
+ Cs\vert \ppp_ty\vert\wwww
\]
and
\[
\vert \ppp_x\ppp_ty\vert \wwww \le \vert \ppp_xw\vert
+ Cs\vert \ppp_ty\vert\wwww,
\]
respectively.

On the other hand, we it is not difficult to see that we can write \eqref{48}  in terms of $w=ze^{s\va} = (\ppp_ty)e^{s\va}$
\begin{equation}\label{eq51}
\begin{array}{l}
\displaystyle\int_{Q_\delta} (s\vert \ppp_xw\vert^2 + s\vert \ppp_tw\vert^2
+ s^3\vert w\vert^2)\,dxdt
\\[4mm]
\phantom{------}\le \displaystyle C\int_{Q_\delta} (\ppp_x(f(x)\ppp_x\ppp_tR(x,t))\vert^2 e^{2s\va} \, dxdt
+ Cs^3\exex M_0^2,
\end{array}
\end{equation}
for all large $s>0$.

\begin{figure}[h!]
\centering
\includegraphics[width=8cm]{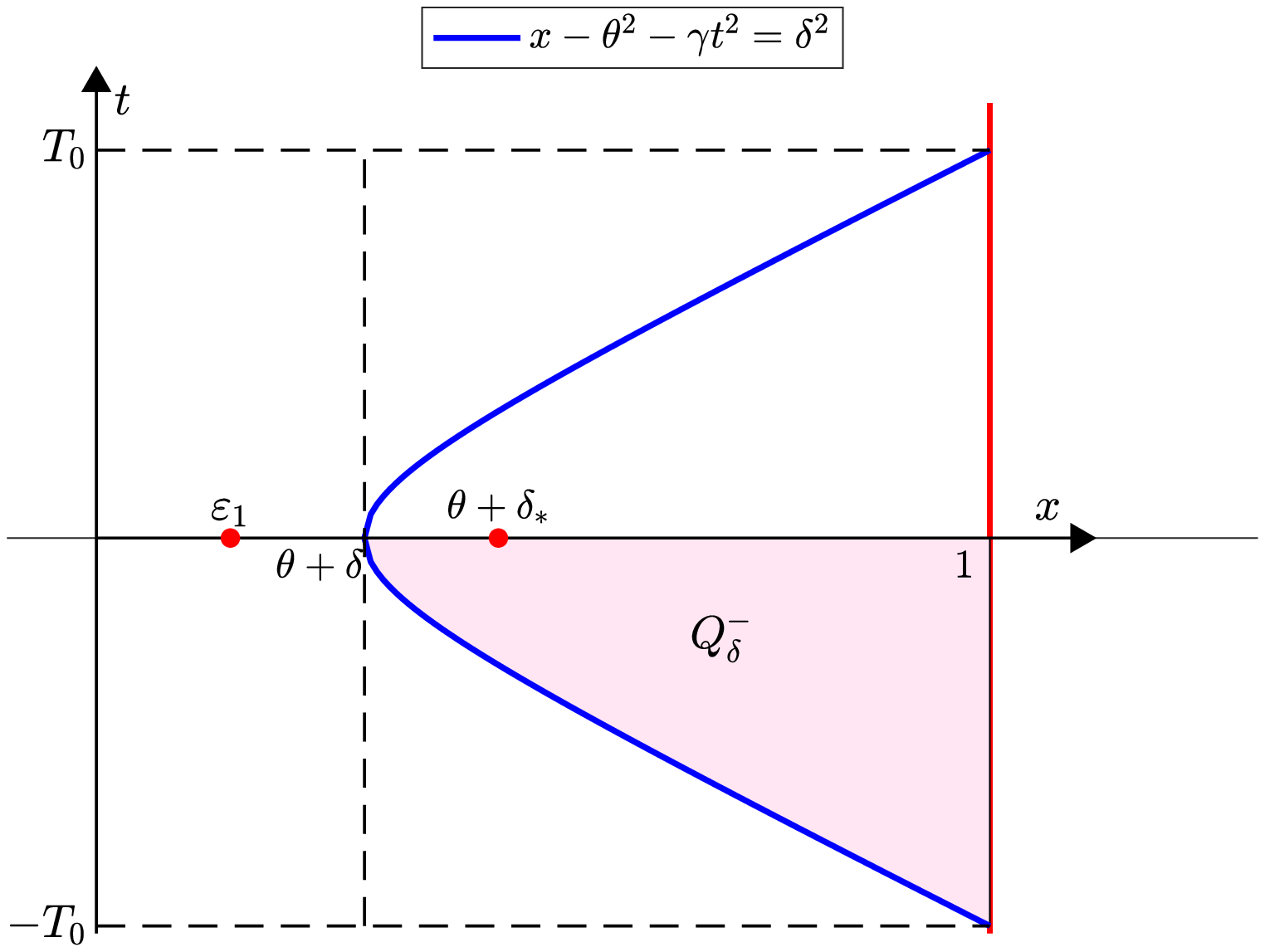}
\caption{Set $Q_\delta^-$.}
\label{Qdeltamenos}
\end{figure}

Multiplying \eqref{eq.w} by $2\ppp_tw$ and integrating over 
$Q_\delta^-:= Q_{\delta} \cap \{ t<0\}$ (see Figure~\ref{Qdeltamenos}), we have
\begin{align}\label{eq57}
B:  &= 2\int_{Q_\delta^-} (\ppp_t^2w)\ppp_tw \,dxdt 
- 2\int_{Q_\delta^-} \ppp_x(p(x)\ppp_xw)\ppp_tw \,dxdt 
\\[2mm]
&= 2\int_{Q_\delta^-} (f'(x)\ppp_t\ppp_xR + f\ppp_t\ppp_x^2R)e^{s\va}
\ppp_tw \,dxdt
+ 2\int_{Q_\delta^-} P_s(w)\ppp_tw \,dxdt. \nonumber
\end{align}

Then, integration by parts in the left hand side of  \eqref{eq57} yields  
\[
B = \int_{Q_\delta^-} \ppp_t(\vert \ppp_tw\vert^2)\, dxdt
+ 2\int_{Q_\delta^-} p(x)(\ppp_xw)(\ppp_t\ppp_xw) \,dxdt
- 2\int_{\ppp Q_\delta^-} p(x)(\ppp_xw)(\ppp_tw) \nu_x \,dS.
\]
Here and henceforth let $\nu = (\nu_x, \nu_t)$ be the unit outward normal 
vector to $\ppp Q_{\delta}^{-}$ in the $(x,t)$-space.  
We have $\nu_x = 0$ and $\nu_t = 1$ on $t=0$.
Since $w = (\partial_t y) e^{s\varphi}$ and $y(1,0) = 0$, we have that
$\ppp_tw(1,t) = \ppp_t((\ppp_ty)e^{s\va})(1,t) = 0$. Moreover, using that
\[
\ppp Q_{\delta}^- = (\ppp Q_{\delta}^- \cap \{ x=1\}) 
\cup (\ppp Q_{\delta}^- \cap \{ x<1\})
= (\ppp Q_{\delta}^- \cap \{ x=1\}) \cup 
(\ppp Q_{\delta}^- \cap \{ \psi(x,t) = \delta^2\}),
\]
we obtain
\begin{align*}
& \left\vert -2\int_{\ppp Q_\delta^-} p(\ppp_xw)(\ppp_tw)\nu_x \,dS
\right\vert 
= \left\vert -2\int_{\ppp Q_\delta^-\cap \{\psi(x,t) = \delta^2\}} 
p(\ppp_xw)(\ppp_tw)\nu_x dS\right\vert              \\[2mm]
\le& \left\vert -2\int_{\ppp Q_\delta^-\cap \{\psi(x,t) = \delta^2\}} 
p(\ppp_x\ppp_ty + s(\ppp_x\va)\ppp_ty)
(\ppp_t^2y + s(\ppp_t\va)\ppp_ty)e^{2se^{\la\delta^2}}\nu_x dS\right\vert 
 \\[2mm]
\le& C\exex s^2 \int_{\ppp Q_\delta^-} 
(\vert \ppp_x\ppp_ty\vert^2 + \vert \ppp_t^2y\vert^2
+ \vert \ppp_ty\vert^2) \,dS
 \\[2mm]
\le& Cs^2\exex (\Vert \ppp_x\ppp_ty\Vert_{H^1(Q_\delta)}^2 
+ \Vert \ppp_t^2y\Vert_{H^1(Q_\delta)}^2 
+ \Vert \ppp_ty\Vert_{H^1(Q_\delta)}^2)
\le Cs^2\exex M_0^2,
\end{align*}
for all large $s>0$.
Hence,
\begin{equation}\label{52}
\left\vert -2\int_{\ppp Q_\delta^-} p(\ppp_xw)(\ppp_tw)\nu_x \,dS
\right\vert \le Cs^2\exex M_0^2.
\end{equation}
In the last inequality, we have used the trace theorem in $Q_{\delta}^-$.

Moreover, using that $\ppp_ty(x,0) = 0$ for $\theta+\delta < x < 1$, we have 
\[
\ppp_xw(x,0) = \ppp_x(e^{s\va}\ppp_ty)\vert_{t=0} = 0, \quad \theta+\delta < x < 1.
\]
Hence, integration by parts yields  
\begin{align*}
 2\int_{Q_\delta^-} p(x)(\ppp_xw)\ppp_t\ppp_xw \,dxdt 
& = \int_{Q_\delta^-} p(x)\ppp_t(\vert \ppp_xw\vert^2) \,dxdt
\\[2mm]
& =  \int^1_{\theta+\delta} p(x)\vert \ppp_xw(x,0)\vert^2 \,dx 
+ \int_{\ppp Q_\delta^-} p(x)\vert \ppp_xw\vert^2 \nu_t \,dS.
\end{align*}
Since $\ppp_ty(x,0)=0$ for $\ep_1<x<1$ implies 
$\ppp_xw(x,0)=0$ for $\theta+\delta<x<1$, similarly to \eqref{52},
by integration by parts, we can estimate the above integral on
$\ppp Q_{\delta}^-$ to see 
\[
\left\vert 2\int_{Q_{\delta}^-} p(\ppp_xw)\ppp_t\ppp_xw \,dx dt
\right\vert
= \left\vert \int_{\ppp Q_{\delta}^-} p\vert \ppp_xw\vert^2 \nu_t \,dS\right\vert
\le Cs^2\exex M_0^2.
\]
Furthermore, we have
\begin{align*}
& \int_{Q_{\delta}^-} \ppp_t(\vert \ppp_tw\vert^2)\, dxdt
= \int^1_{\theta+\delta} \vert \ppp_tw(x,0)\vert^2 \,dx 
+ \int_{\ppp Q_{\delta}^-} \vert \ppp_tw\vert^2 \nu_t \,dS\\
\ge & \int^1_{\theta+\delta} \vert \ppp_tw(x,0)\vert^2 \,dx 
- Cs^2\exex M_0^2.
\end{align*}
Therefore, 
\begin{equation}\label{53}
\vert B\vert 
\ge \int^1_{\theta+\delta} \vert \ppp_tw(x,0)\vert^2 dx 
- Cs^2\exex M_0^2.
\end{equation}
On the other hand, by \eqref{eq3.2} and \eqref{50}, we can estimate the 
right-hand side of (\ref{eq57}) as follows:
\begin{align*}
& \left\vert 2\int_{Q_\delta^-} (f'(x)\ppp_t\ppp_xR 
+ f(x)\ppp_t\ppp_x^2R)e^{s\va}
\ppp_tw dxdt + 2\int_{Q_\delta^-} P_s(w)\ppp_tw dxdt\right\vert
 \\[2mm]
\le& C\int_{Q_\delta^-} 
(\vert f'(x)\vert + \vert f(x)\vert) e^{s\va}\vert \ppp_tw\vert\, dxdt
+ C\int_{Q_\delta^-} (s^2\vert w\vert + s\vert \ppp_xw\vert
+ s\vert \ppp_tw\vert)\vert \ppp_tw\vert \,dxdt.
\end{align*}
Applying the Cauchy-Schwarz inequality to the right-hand side and using 
\[
s^2\vert w \vert \vert\ppp_tw\vert = s^{\frac{3}{2}}\vert w\vert 
s^{\frac{1}{2}}\vert \ppp_tw\vert 
\le s^3\vert w\vert^2 + s\vert \ppp_tw\vert^2,
\]
%
by (\ref{eq51})  and \eqref{eq3.2}  we can obtain
\begin{align*}
& \Big\vert \text{[the right-hand side of (\ref{eq57})]}\Big\vert \\
\le& C\int_{Q_\delta} 
(\vert f'(x)\vert + \vert f(x)\vert) e^{2s\va}\, dxdt
+ C\int_{Q_\delta} \vert \ppp_x(f\ppp_x\ppp_tR)\vert^2 e^{2s\va}\,
dxdt + Cs^3\exex M_0^2\\
\le& C\int_{Q_\delta} (\vert f'(x)\vert + \vert f(x)\vert) e^{2s\va} \, dxdt
+ Cs^3\exex M_0^2.
\end{align*}
Therefore, in view of \eqref{53}, we reach 
\begin{equation}\label{61}
\int_{\theta+\delta}^1 \vert \ppp_tw(x,0)\vert^2\, dx
\le C\int_{Q_\delta} (\vert f'(x)\vert^2 + \vert f(x)\vert^2) 
e^{2s\va(x,t)} dxdt
+ Cs^3\exex M_0^2,
\end{equation}
for all large $s>0$.

Now, we recall that 
$T_0 = \sqrt{\frac{(1-\theta)^2-\delta^2}{\gamma}}$ and that $w=(\partial_t y)e^{s\varphi}$.
By \eqref{3.1} we have
\begin{align*}
\ppp_tw(x,0) & = \ppp_t^2y(x,0)e^{s\va(x,0)}
+ (\ppp_ty)(x,0)(\ppp_t\wwww)(x,0)          \\[2mm]
& = \ppp_t^2y(x,0)e^{s\va(x,0)} 
= \ppp_x(f(x)\ppp_xR(x,0))e^{s\va(x,0)}.
\end{align*}
Hence, \eqref{61} yields
\begin{equation}\label{571}
\begin{array}{l}
\displaystyle\int^1_{\theta+\delta} 
\vert f'(x)\ppp_xR(x,0) + f(x)\ppp_x^2R(x,0)\vert^2 e^{2s\va(x,0)}\,dx 
\\[4mm]
\le C\displaystyle\int_{Q_{\delta}} (\vert f(x)\vert^2 + \vert f'(x)\vert^2) 
e^{2s\va(x,t)}\,dxdt + Cs^3\exex M_0^2,
\end{array} 
\end{equation} 
for all large $s>0$.

Since  
\[
\va(x,t) - \va(x,0) = -e^{\la\vert x-\theta\vert^2}(1-e^{-\la\gamma t^2})
\le -(1-e^{-\la \gamma t^2}),
\]
by \eqref{43} we obtain
\begin{align*}
&\int_{Q_{\delta}} (\vert f(x)\vert^2 + \vert f'(x)\vert^2) 
e^{2s\va(x,t)}\, dxdt
\le C\int^{T_0}_{-T_0}
\left(\int^1_{\theta+\delta} (\vert f\vert^2 + \vert f'\vert^2) 
e^{2s\va(x,0)}\,dx \right) e^{2s(\va(x,t) - \va(x,0))}\, dt      \\[2mm]
& \le C \int^1_{\theta+\delta} \left( \int^{T_0}_{-T_0} 
e^{-2s(1-e^{-\la\gamma t^2})} dt \right) 
(\vert f\vert^2 + \vert f'\vert^2) e^{2s\va(x,0)}\,dx.
\end{align*}
The Lebesgue convergence theorem yields   
\[
\int^{T_0}_{-T_0} e^{-2s(1-e^{-\la\gamma t^2})} dt = o(1)
\]
as $s \to \infty$, and so
\[
\int_{Q_{\delta}} \vert f(x)\vert^2 + \vert f'(x)\vert^2) 
e^{2s\va(x,t)} \,dxdt 
= o(1) \int^1_{\theta+\delta}
(\vert f(x)\vert^2 + \vert f'(x)\vert^2) e^{2s\va(x,0)}\,dx.
\]
Applying this in \eqref{571}, we reach 
\begin{align}\label{3.11}
\int^1_{\theta+\delta} \vert f'(x)\ppp_xR(x,0) + f(x) \ppp_x^2R(x,0)\vert^2
e^{2s\va(x,0)} \, dx
            \nonumber       \\[2mm]
\le o(1) \int^1_{\theta+\delta}
(\vert f(x)\vert^2 + \vert f'(x)\vert^2) e^{2s\va(x,0)}\, dx
+ Cs^3\exex M_0^2                     
\end{align}
for all large $s>0$.

Using now that $\vert \ppp_xR(x,0)\vert \ge r_0 > 0$ for all $x\in [\ep_1,1]$ and 
$\ppp_x^2R(\cdot,0)\in L^{\infty}(\ep_1,1)$ by \eqref{eq3.2}, we see 
\begin{align*}
& \int^1_{\theta+\delta} \vert f'(x)\ppp_xR(x,0) + f(x) 
\ppp_x^2R(x,0)\vert^2 e^{2s\va(x,0)} \,dx \\
\ge & r_0^2\int^1_{\theta+\delta} \vert f'(x)\vert^2 e^{2s\va(x,0)} \,dx
- C\int^1_{\theta+\delta} \vert f(x)\vert^2 e^{2s\va(x,0)} \,dx.
\end{align*}
Hence, \eqref{3.11} yields 
\begin{align*}
& r_0^2\int^1_{\theta+\delta} \vert f'(x)\vert^2 e^{2s\va(x,0)}dx
- C\int^1_{\theta+\delta} \vert f(x)\vert^2 e^{2s\va(x,0)}dx\\
\le& o(1)\int^1_{\theta+\delta} \vert f'(x)\vert^2 e^{2s\va(x,0)}dx
+ o(1)\int^1_{\theta+\delta} \vert f(x)\vert^2 e^{2s\va(x,0)}dx
+ Cs^3\exex M_0^2,
\end{align*}
for all large $s>0$.  
Absorbing the first term on the right-hand side into the left-hand side.
we reach 
\[
\int^1_{\theta+\delta} \vert f'(x)\vert^2 e^{2s\va(x,0)}dx
\le C\int^1_{\theta+\delta} \vert f(x)\vert^2 e^{2s\va(x,0)}dx\\
+ Cs^3\exex M_0^2,
\]
for all large $s>0$. 
Applying Lemma~\ref{Lemma6} to the first term on the right-hand side if this inequality,
we obtain
\[
C_2\int^1_{\theta+\delta} \vert f'(x)\vert^2 e^{2s\va(x,0)}dx
- \frac{C}{s}\int^1_{\theta+\delta} \vert f'(x)\vert^2 e^{2s\va(x,0)}dx
\le Cs^3\exex M_0^2,
\]
for all large $s>0$.  
Choosing $s>0$ sufficiently large, we can absorb the second term into the 
first term on the left-hand side, and we reach
\begin{equation}\label{64}
C_3\int^1_{\theta+\delta} \vert f'(x)\vert^2 e^{2s\va(x,0)}dx
\le Cs^3\exex M_0^2,
\end{equation}
for all large $s>0$.

Next, choose $\delta_* > \delta$ such that $\theta+\delta_* < 1$ (see Figure~\ref{Qdeltamenos}).
Then 
$$
\va(x,0) = e^{\la (x-\theta)^2}
\ge e^{\la\delta_*^2} > e^{\la\delta^2} \quad \mbox{for 
$\theta+\delta_* \le x \le 1$},
$$
and so \eqref{64} yields 
\begin{align*}
&C_3e^{2se^{\la\delta_*^2}}\Vert f'\Vert^2_{L^2(\theta+\delta_*,1)}
\le C_3\int^1_{\theta+\delta_*} \vert f'(x)\vert^2 e^{2s\va(x,0)}dx\\
\le& C_3\int^1_{\theta+\delta} \vert f'(x)\vert^2 e^{2s\va(x,0)}dx
\le Cs^3\exex M_0^2,
\end{align*}
that is,
\[
\Vert f'\Vert^2_{L^2(\theta+\delta_*,1)}
\le C_4s^3M_0^2 
\exp(-2s(e^{\la\delta_*^2} - e^{\la\delta^2})),
\]
for all large $s>0$.  By $e^{\la\delta_*^2} - e^{\la\delta^2} > 0$, letting
$s\to \infty$, we see that $f'=0$ in $(\theta+\delta_*,\, 1)$.
Since $\theta_* > \theta$ is arbitrary, we deduce that $f'=0$ in 
$(\theta+\delta, \, 1)$.  Since $f(1) = 0$, we see that 
$f=0$ in $(\theta+\delta, \, 1)$.  
\\
\vspace{0.2cm}

\bigskip
\noindent
{\bf Third step of the proof of Proposition~\ref{Proposition1}: 
continuation of the uniqueness of $f$ and $y$ toward $x=\ep_1$.}

Substituting $f=0$ in $(\theta+\delta,1)$ in \eqref{46}, we have
\begin{equation}\label{56}
\left\{ \begin{array}{rl}
& \ppp_t^2y = \ppp_x(p(x)\ppp_xy), \quad \theta +\delta < x < 1, \,
-\infty<t<\infty, \\
& y(1,t) = \ppp_xy(1,t) = 0, \quad -\infty<t<\infty.
\end{array}\right.
\end{equation}
We apply \eqref{44} to \eqref{56} and obtain
\[
\int_{Q_\delta} s^3\vert y\vert^2 \weight\, dxdt 
\le Cs^3\exex M_0^2
\]
for all large $s>0$.  Hence, choosing $\delta_* > \delta$ again such that 
$\theta + \delta_* < 1$, and using that $Q_{\delta} \supset Q_{\delta_*}$, we have 
\[
s^3e^{2se^{\la\delta_*^2}}\Vert y\Vert^2_{L^2(Q_{\delta_*})}
\le \int_{Q_{\delta_*}} s^3\vert y\vert^2 \weight \, dxdt 
\le \int_{Q_{\delta}} s^3\vert y\vert^2 \weight \,dxdt 
\le Cs^3\exex M_0^2,
\]
that is,
\[
\Vert y\Vert^2_{L^2(Q_{\delta_*})}
\le C\exp(-2s(e^{\la\delta_*^2} - e^{\la\delta^2}))M_0^2,
\]
for all large $s>0$. Letting $s\to \infty$, we obtain 
$y=0$ in $Q_{\delta_*}$.  Since $\delta_*>\delta$ is arbitrary, we have
$y=0$ in $Q_{\delta}$ as well as $f=0$ in $(\theta+\delta,1)$.
Since also $\delta>0$ is arbitrary, we see that 
$$
\left\{ \begin{array}{rl}
& y = 0 \quad \mbox{in 
$Q_0 := \{ (x,t);\, -\infty<t<\infty, \quad 
\theta < x < 1, \quad \psi(x,t) > 0\}$ }, \\
& f = 0 \quad \mbox{in $(\theta,1)$}.
\end{array}\right.
$$
We note that 
$$
Q_0 = \left\{ (x,t);\, \vert t\vert < \frac{1}{\sqrt{\gamma}}(x-\theta), \,\,
\theta < x < 1\right\}.
$$

For arbitrary $t_0\in \R$, we repeat the same arguments to \eqref{56} with 
translated time variable $t-t_0$, so that 
$$
y=0 \quad \mbox{in 
$\left\{ (x,t);\, \vert t-t_0\vert < \frac{1}{\sqrt{\gamma}}(x-\theta), \,\,
\theta < x < 1\right\}$}.
$$
Since $t_0\in \R$ is arbitrary, we see that 
$$
y = 0 \quad \mbox{in 
$\bigcup_{t_0\in \R} \left\{ (x,t);\, \vert t-t_0\vert 
< \frac{1}{\sqrt{\gamma}}(x-\theta), \,\,
\theta < x < 1\right\} \supset (\theta,1) \times (-\infty,\infty)$}.
$$
Using this with \eqref{3.1}, we derive
$$
\left\{ \begin{array}{l}
\ppp_t^2y = \ppp_x(p(x)\ppp_xy) + \ppp_x(f(x)\ppp_xR(x,t)), 
\quad \ep_1 < x < \theta, \,
-\infty<t<\infty, 
\\[1mm]
 y(\theta,t) = \ppp_xy(\theta,t) = 0, \quad -\infty<t<\infty,
 \\[1mm]
y(x,0) = \ppp_ty(x,0) = 0, \quad \ep_1 < x < \theta.
\end{array}\right.
$$
\bigskip

Now, starting at $x=1$, we advance by the width $\theta$ to obtain
the uniqueness from $x=1$ to $x=\theta$.
In the second step of this proof, we have applied the Carleman estimate in $Q_{\delta}
\subset \{x;\, 1 - (1-\theta) < x < 1\}$ with the weight function
$\psi(x,t) = (x - \theta)^2 - \gamma t^2$.
Now we have to consider a level set in $\{ x;\, \theta - (1-\theta))
< x < \theta \}$ and so we shift the level set $Q_{\delta}$ 
by $x \mapsto x - (\theta-1)$:
$$
Q^1_{\delta} := \{ (x,t);\, -\infty<t<\infty, \quad 
\psi^1(x,t) > \delta^2\},
$$
where 
$$
\psi^1(x,t) := ((x - (\theta-1))-\theta)^2 - \gamma t^2.
$$
Since $1-\theta > 0$ is sufficiently small, by Remark~\ref{remark4.9} we can 
verify the Carleman estimate  from Lemma~\ref{Lemma5} in $Q^1_{\delta}$.
Then we repeat the same argument fro the second step of this proof
by replacing the $x$-interval $(\ep_1, 1)$ by $(\ep_1, \theta)$ and we can 
prove
\[
y(x,t)=0, \quad (x,t) \in (1-2(1-\theta), \, 1-(1-\theta)) \times 
(-\infty, \infty)
\]
and
\[
f(x)=0, \quad x\in (1-2(1-\theta), \, 1-(1-\theta)).
\]
Adjusting the weight function for the Carleman estimate and 
continuing this argument, we reach 
\[
f(x)=0, \quad x\in (1-k(1-\theta), \, 1),
\]
with $k \in \N$ as long as $1-k(1-\theta) > \ep_1$.

We choose $k_0\in \N$ such that 
$$
1-(k_0+1)(1-\theta) < \ep_1 < 1 - k_0(1-\theta).
$$
Then, $(1 - k_0(1-\theta)) - \ep_1 < 1-\theta$ and 
\begin{equation}\label{65}
f=0 \quad \mbox{in $(1 - k_0(1-\theta), \, 1)$}.
\end{equation}
Since we have $1-\theta > 0$ can be arbitrarily small, we can choose 
$\theta >0$ and $k_0\in \N$ such that 
$1-k_0(1-\theta) > \ep_1$ and $1-k_0(1-\theta)$ is arbitrarily close to 
$\ep_1$.   By the continuity 
of $f$, we reach $f=0$ in $(\ep_1,1)$.

Thus, the proof of Proposition~\ref{Proposition1} is complete. \hfill $\square$

\subsubsection{Third Step: Completion of the proof of 
Theorem~\ref{th.uniqueness}}
\label{sec.proofth}

Now we proceed to  finish the proof of  Theorem~\ref{th.uniqueness}. 
Recall that we have assumed \eqref{28}.  Then, Lemma~\ref{Lemma4} yields
\[
\ppp_x\www{u}_{a,\alpha}(1,t) = \ppp_x\www{u}_{b,\beta}(1,t), 
\quad 0<t<\infty.
\]
Here $\www{u}_{a,\alpha}$ and $\www{u}_{b,\beta}$ are the solutions to 
\eqref{2.1} corresponding to $(a,\alpha)$ and $(b,\beta)$, respectively.
Setting 

\[p(x) = x^{\alpha}a(x)\quad y=\www{u}_{a,\alpha} - \www{u}_{b,\beta}, \quad 
f = x^{\alpha}a(x) - x^{\beta}b(x), \quad R(x,t) = \www{u}_{b,\beta}(x,t), 
\]
we have \eqref{3.1} and $\ppp_xy(1,t) = 0$ for $t>0$.

In view of Proposition~\ref{Proposition1}, using the regularity \eqref{2.2}-\eqref{2.4} 
of $y$ and $R$ guaranteed by Lemma~\ref{Lemma1}, we obtain $f=0$ in 
$\ep_1 < x < 1$. 
Thus, the proof of Theorem~\ref{th.uniqueness} is complete. \hfill $\square$ 

\bigskip

\section{Numerical results}\label{sec.numerics1}

In this section, we will consider numerical results for the previous inverse problems. 
We will carry out the reconstruction of the unknown functions $a$ and the power $\alpha$  through the resolution of some appropriate extremal problems. This strategy has been applied in some previous papers for other similar problems, see~\cite{Jone} and~\cite{doubova}. The results of the numerical tests that follow pretend to illustrate the theoretical results in the previous sections.

\subsection{Numerical tests for the linear case I}\label{sec:lincase1}

In this section we present the numerical reconstruction of the constant coefficient $a$ in \eqref{dpe1} by distributed measurements, related to the theoretical result from Theorem~\ref{th:linearcase}. More precisely, we will deal with the following inverse problem: given $t_0\in (0,T)$ as in \eqref{LB},  $u_0=u_0(x)$, $\gamma = \partial_t u(\cdot,t_0) $ and $\beta =u(\cdot,t_0) $, find 
$a\in \mathcal{U}^{a}_{ad}$ such that the corresponding solution $u^a$ of \eqref{dpe1} satisfies the additional conditions: 
\begin{equation*}\label{dpe1_a}
\partial_t u^a(x,t_0) = \gamma(x), \quad  x\partial_x u^a(x,t_0) = \beta(x), \quad  \forall\, x\in (0,\ell),
\end{equation*}
where the admissible set $\mathcal{U}^{a}_{ad}=\{ a: 0<\underline{a}\le a \le \overline{a}  \}$. 

We reformulate the inverse problem as the following extremal problem: find $a\in \mathcal{U}^{a}_{ad}$ such that 
\begin{equation}\label{pb.const11}
J(a) \le  J(a'),\quad \forall\, a'\in \mathcal{U}^{a}_{ad}\,,
\end{equation}
where the functional $J : a\in \mathcal{U}^{a}_{ad} \mapsto \mathds{R}$ is given by

\begin{equation}\label{J.const1}
J(a) = \dfrac{1}{2}\displaystyle\int_0^\ell |\gamma(x) - \partial_t u^a(x,t_0)|^2\, dx 
+  \dfrac{1}{2}\displaystyle\int_0^\ell |\beta(x) - x\partial_x u^a(x,t_0)|^2\, dx,
\end{equation}
with $u^a$ the solution of \eqref{dpe1} corresponding to the unknown $a$. 

\bigskip
We have performed several numerical tests with different values of $a$ to recover: $a>1$, $a=1$ and $a<1$. In them, we fix $T=5$, $\ell=1$ and $u_0(x) = 0.5x^2(1-x)$. In the first three tests we will take, for simplicity, $f(x,t) = 0$. We will denote by $a_d$ the desired value of $a$ that we would like to obtain and by $a_i$ the initial guess for the optimization algorithm. 
Using an optimization algorithm from Optimization Toolbox of MatLab, the \texttt{fmincon} function and \texttt{interior-point} gradient based algorithm, we have computed $a_c$. Our goal is to obtain $a_c$ as close as possible to $a_d$.

\begin{figure}[h!]
\begin{minipage}[t]{0.49\linewidth}
\centering
\includegraphics[width=0.9\linewidth]{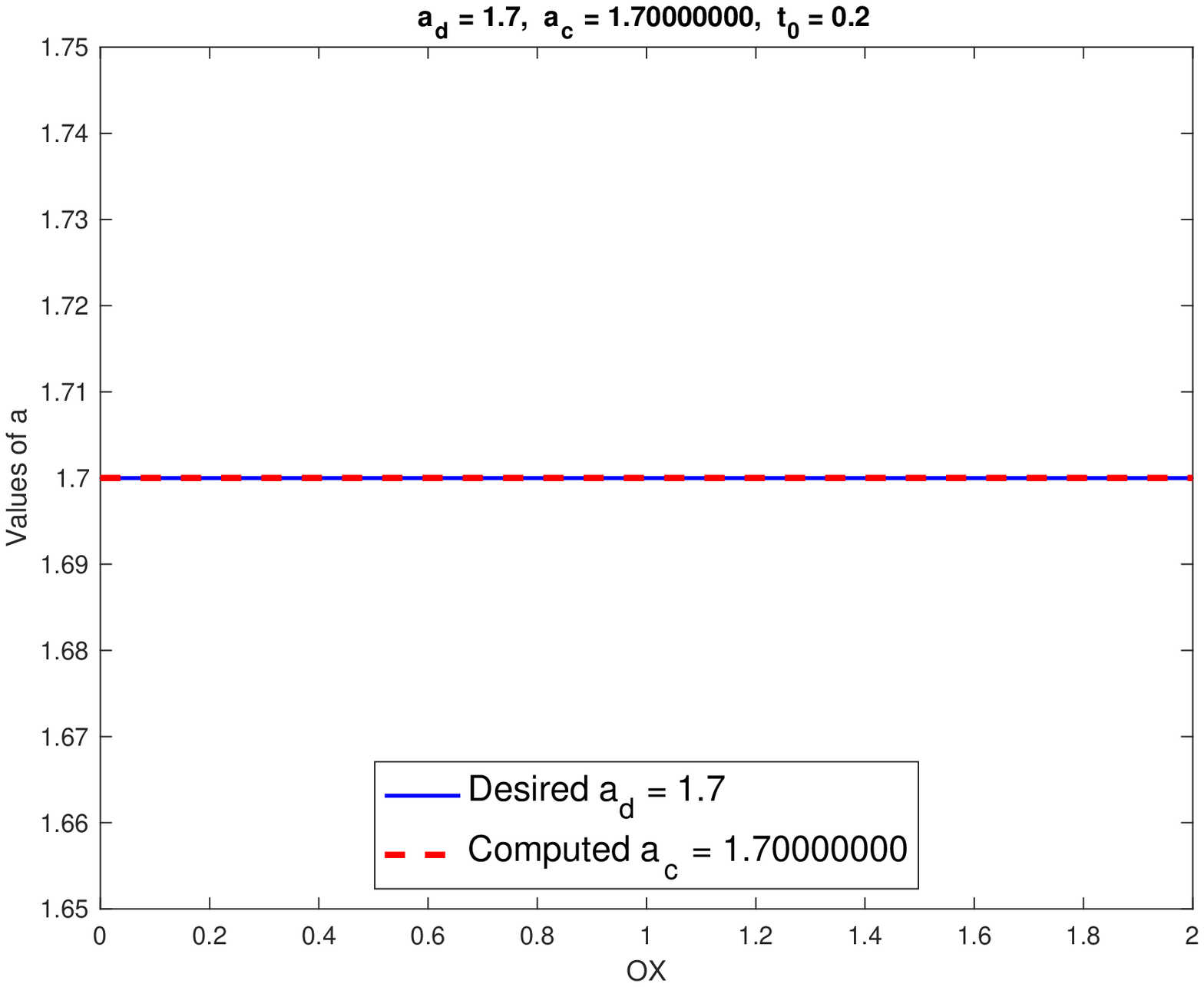}
\caption{Test~\ref{test.const11}, $a_d$ solid line, $a_c$ dashed line.}
\label{alpha171}
\end{minipage}
\hfill
\begin{minipage}[t]{0.49\linewidth}
\centering
\includegraphics[width=0.9\linewidth]{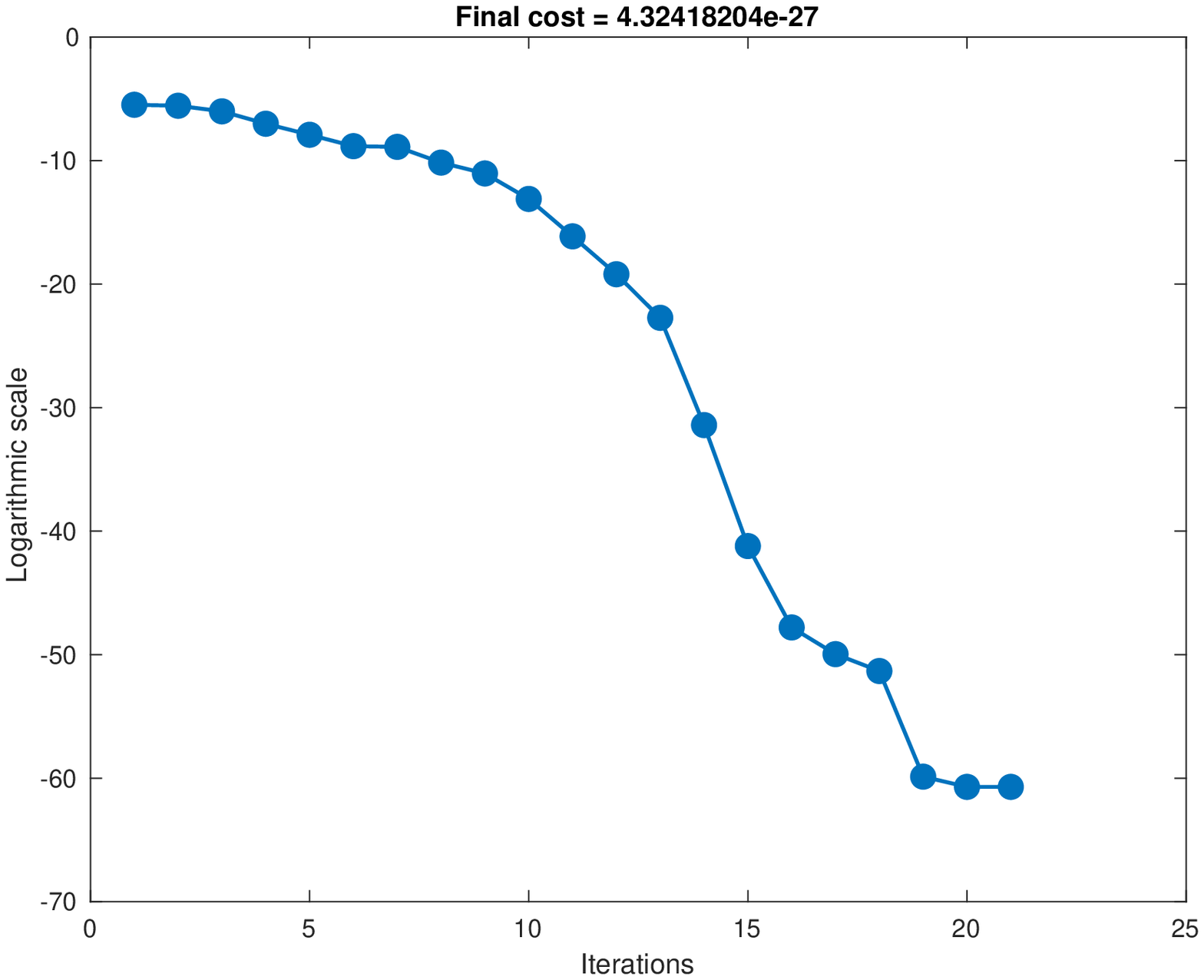}
\caption{Test~\ref{test.const11}, the evolution of the cost.}
\label{fvalalpha171}
\end{minipage}
\end{figure}

\begin{test}\label{test.const11}
Let us take $a_{i} = 0.7$ as an initial guess for the minimization algorithm and $t_0 = 0.2$, that we have computed numerically such that \eqref{LB} holds. The desired value to recover is $a_d=1.7$.  After numerical resolution of the optimization problem \eqref{pb.const11}, we found $a_c = 1.7$. The numerical results can be seen in Figure~\ref{alpha171}, where the desired $a_d$ and computed $a_c$ have been represented. In Figure~\ref{fvalalpha171} we have plotted the evolution of the cost at each iteratios of the optimization algorithm represented by blue points.  In Table~\ref{table11.const} we can see the evolution of the cost when we introduce random noises in the target.
\end{test}

\begin{minipage}[t]{0.49\textwidth}
\centering
\renewcommand{\arraystretch}{1.2}
\medskip
\captionof{table}{Cost at the best stopping and computed $a_c$ with random noises in the target in Test~\ref{test.const11}.}
\medskip
{\begin{tabular}{cccc} \hline
\% noise  & Cost  & Iterations &  $a_c$ \\
\hline
1\%           &   1.e-7       &      17       &   1.70451153  \\ \hline
0.1\%        &    1.e-8      &      16       &  1.70017891\\ \hline
0.01\%      &   1.e-10     &      16      &   1.70017665\\\hline
0.001\%    &   1.e-13     &      15       &   1.70000362\\ \hline
0\%           &   1.e-27     &      20       &  1.7  \\  \hline
\end{tabular}}
\label{table11.const}
\end{minipage}
\hfill
\begin{minipage}[t]{0.49\textwidth}
\vspace{0cm}
\centering
\includegraphics[width = 0.9\linewidth]{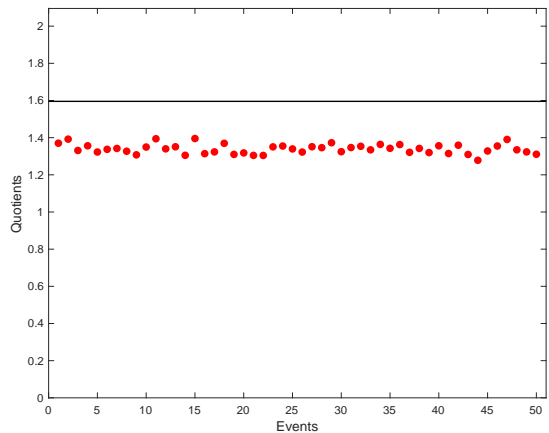}
 \captionof{figure}{Test~\ref{test.const11}, numerical illustration of stability.}
\label{fig.numstab}
\end{minipage}

\bigskip
\begin{remark}[Numerical illustration of stability] In order to illustrate the stability of the considered inverse problem in the Test~~\ref{test.const11}, we have represented in Figure~\ref{fig.numstab} as points the quotients 
\begin{equation}\label{eq:K}
K_j = \dfrac{|a-a^\ast|}{\Vert \gamma(x) - \gamma_j^\ast(x)\Vert + \Vert \beta(x) - \beta_j^\ast(x)\Vert }, \quad j = 1\dots, 50,
\end{equation}
where  $\gamma(x)= \partial_t u^a(x,t_0)$,  $\beta(x)= x\partial_x u^a(x,t_0)$ and for all $j$, $\gamma^\ast_j(x)= \gamma(x)+\varepsilon_j$,  $\beta^\ast_j(x)= \beta(x)+\varepsilon_j$ with the $\varepsilon_j$  random perturbations in $(0, 0.05)$. We can see that all the quotients $K_j$ are bounded.

\end{remark}

\bigskip

\begin{test}\label{test.const22}
Let us now take $a_{i} = 0.2$, $t_0 = 0.2$ and the desired value to recover $a_d=1$. The optimization algorithm gives $a_c = 1$, see Figure~\ref{fvalalpha02lin2}.  In Table~\ref{table2.const2} we can see the evolution of the cost when we introduce random noises in the target.
\end{test}
\medskip

\begin{minipage}[t]{0.49\textwidth}
\centering
\renewcommand{\arraystretch}{1.2}
\medskip
 \captionof{table}{Evolution of the cost and computed $a_c$ with random noises in the target in Test~\ref{test.const22}.}
\medskip
\begin{tabular}{cccc} \hline
\%  noise  & Cost  & Iterations &   $a_c$ \\
\hline
1\%           &   1.e-6       &      12       & 0.99670260  \\ \hline
0.1\%        &    1.e-8      &      13        & 1.00012771\\ \hline
0.01\%      &   1.e-14     &       10        &   0.99995171\\\hline
0.001\%    &   1.e-13     &      14        &  0.99999488\\ \hline
0\%           &   1.e-24     &      14       &  1   \\  \hline
\end{tabular}
\label{table2.const2}
\end{minipage}
\hfill
\begin{minipage}[t]{0.49\textwidth}
\vspace{0cm}
\centering
\includegraphics[width=0.9\textwidth]{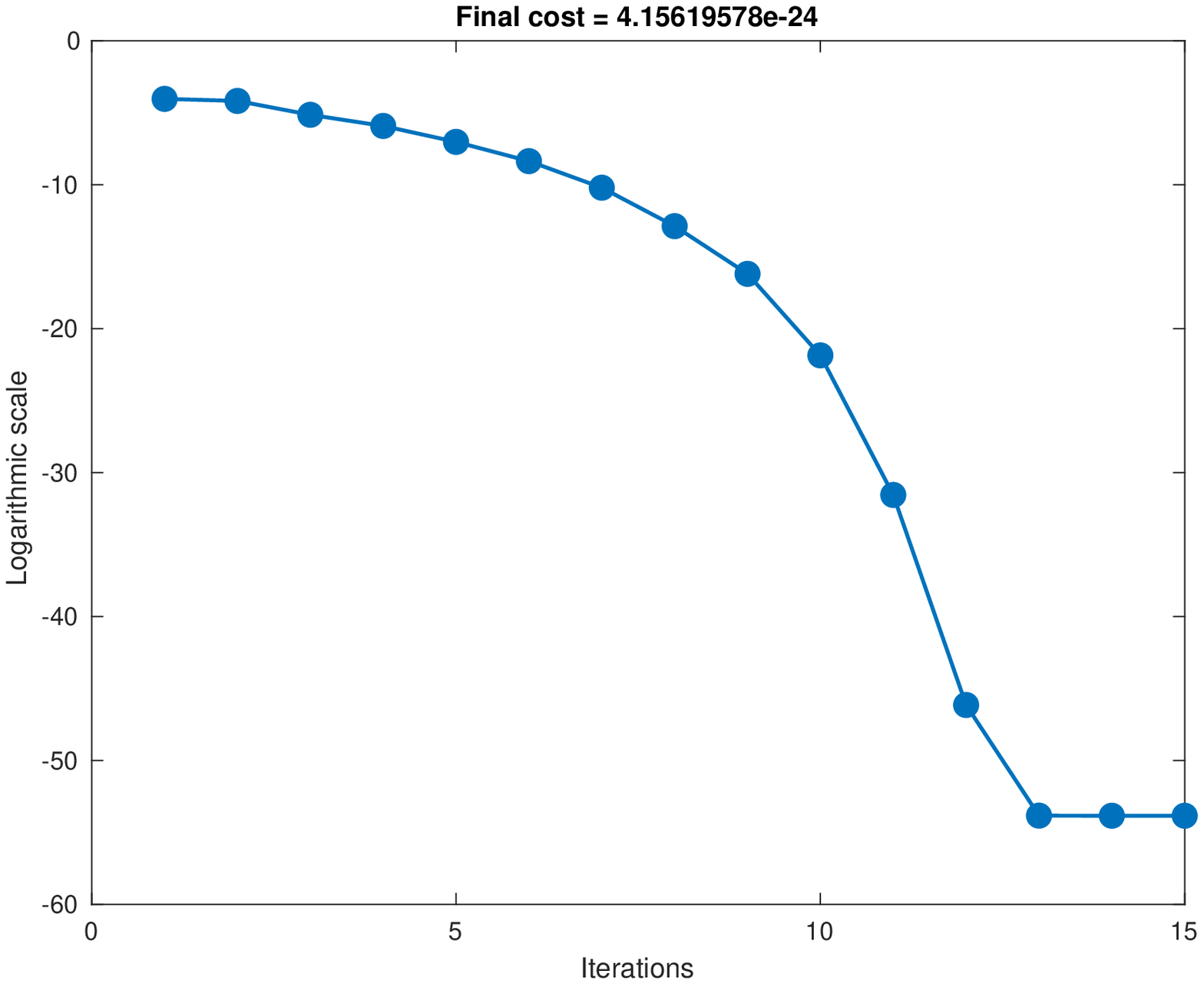}
 \captionof{figure}{Test~\ref{test.const22}, the evolution of the cost.}
\label{fvalalpha02lin2}
\end{minipage}

\bigskip
\begin{test}\label{test.const33}
Let us take $a_{i} = 0.7$ as an initial guess for the minimization algorithm, $t_0 = 0.2$ and $a_d=0.2$. The numerical resolution of the optimization problem gives $a_c=0.2$ (see Figure~\ref{fvalalphalin33} and Table~\ref{table3.const3}). 
\end{test}
\medskip

\begin{minipage}[t]{0.49\linewidth}
\centering
\renewcommand{\arraystretch}{1.2}
\medskip
\captionof{table}{Evolution of the cost and computed $a_c$ with random noises in the target in Test~\ref{test.const33}.}
\medskip
{\begin{tabular}{cccc} \hline
\%  noise  & Cost  & Iterations &   $a_c$ \\
\hline
1\%           &   1.e-7      &    11       & 0.20046140  \\ \hline
0.1\%        &   1.e-13      &    9      &   0.19997676\\ \hline
0.01\%      &   1.e-11       &    16       &   0.19999893\\\hline
0.001\%    &   1.e-13      &    14       &  0.20000028\\ \hline
0\%           &   1.e-27      &   16        &  0.2   \\  \hline
\end{tabular}}
\label{table3.const3}
\end{minipage}
\hfill
\begin{minipage}[t]{0.49\linewidth}
\vspace{0cm}
\centering
\includegraphics[width=0.9\linewidth]{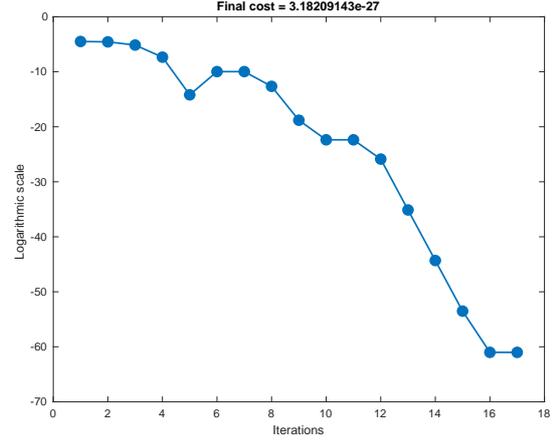}
 \captionof{figure}{Test~\ref{test.const33}, the evolution of the cost.}
\label{fvalalphalin33}
\end{minipage}

\bigskip

\begin{figure}[h!]
\begin{minipage}[t]{0.49\linewidth}
\centering
\includegraphics[width=0.9\linewidth]{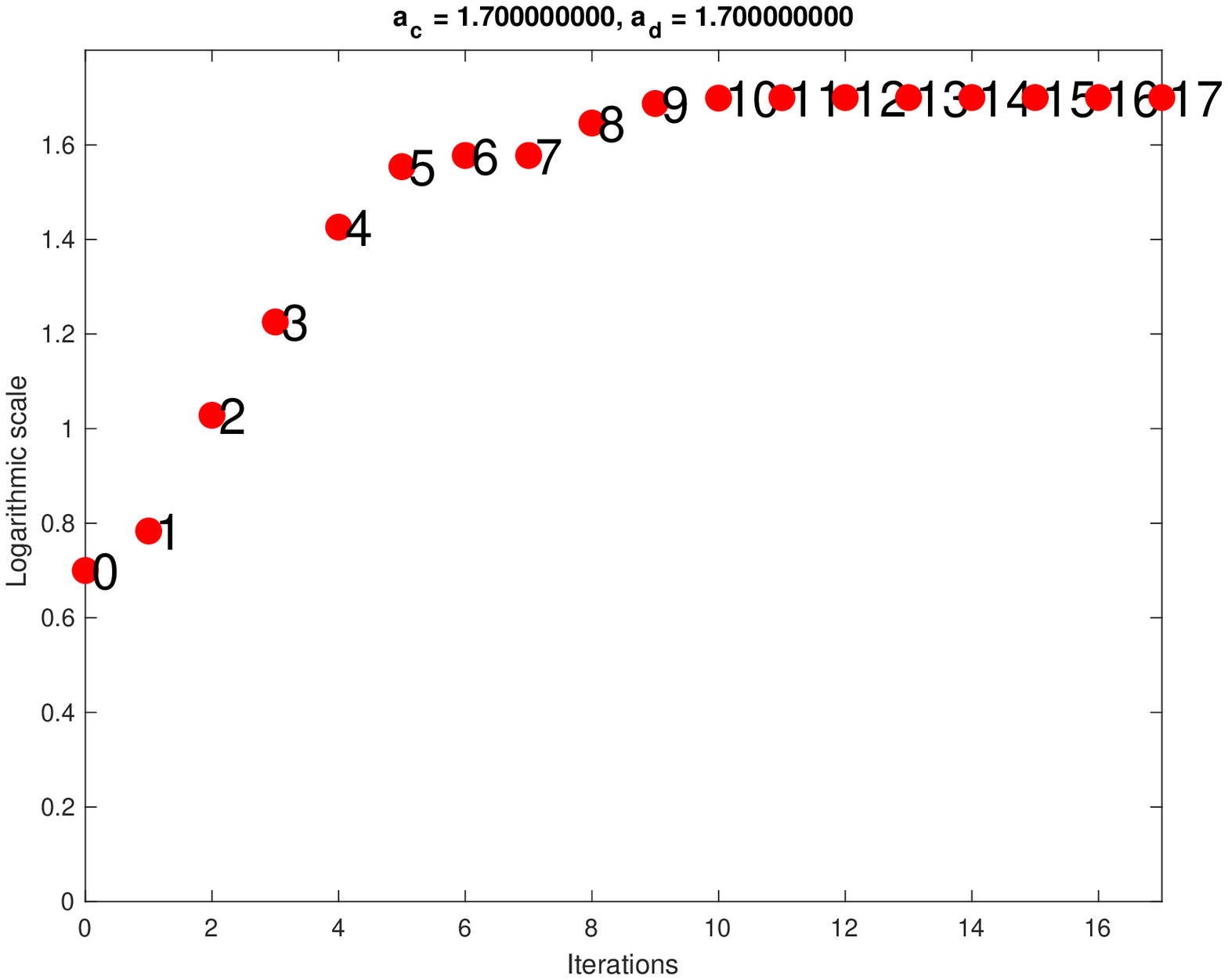}
\caption{Test~\ref{test.const33BIS}, iterations of the algorithm.}
\label{alpha_iter_nh}
\end{minipage}
\hfill
\begin{minipage}[t]{0.49\linewidth}
\centering
\includegraphics[width=0.9\linewidth]{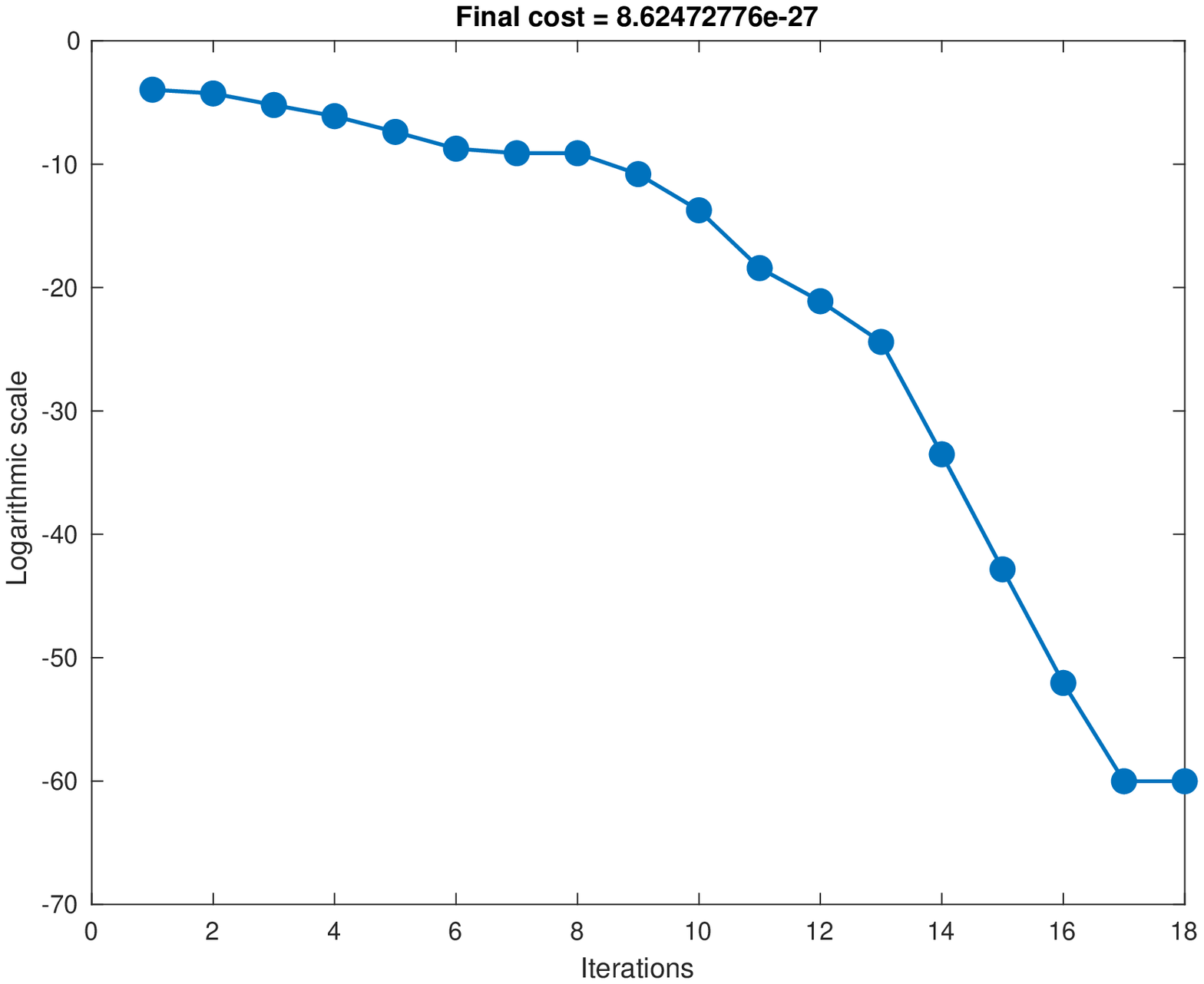}
\caption{Test~\ref{test.const33BIS}, the evolution of the cost.}
\label{fvalalphalin_nonh}
\end{minipage}
\end{figure}

\begin{test}\label{test.const33BIS}
We consider a non-homogeneous problem  \eqref{dpe1}  with $f(x,t) =2xt$.  
Let us take $a_{i} = 0.7$ and  $t_0 = 0.4$. The desired value is $a_d=1.7$. The numerical resolution of the optimization problem gives $a_c=1.7$. Figure~\ref{alpha_iter_nh} shows the evolution of iterations of the optimization algorithm and Figure~\ref{fvalalphalin_nonh}  the evolution of the cost. In Table~\ref{table3.const_nh} we present the results with random noises  in the target and Figure~\ref{fig.sol_nh} representes the solution corresponding to the computed $a_c$.
\end{test}
\medskip

\begin{minipage}[t]{0.49\linewidth}
\centering
\renewcommand{\arraystretch}{1.2}
\captionof{table}{Cost at the best stopping and  $a_c$ with random noises in the target in Test~\ref{test.const33BIS}.}
\medskip
{\begin{tabular}{cccc} \hline
\% noise  & Cost  & Iterations &   $a_c$ \\
\hline
1\%           &   1.e-7       &      14       &   1.69767371  \\ \hline
0.1\%        &    1.e-8      &      15       &  1.69992839\\ \hline
0.01\%      &   1.e-11     &      18     &   1.69997876\\\hline
0.001\%    &   1.e-12     &      15       &   1.69999972\\ \hline
0\%           &   1.e-27     &      17       &  1.7  \\  \hline
\end{tabular}}
\label{table3.const_nh}
\end{minipage}
\hfill
\begin{minipage}[t]{0.49\linewidth}
\vspace{0cm}
\centering
\includegraphics[width=0.9\linewidth]{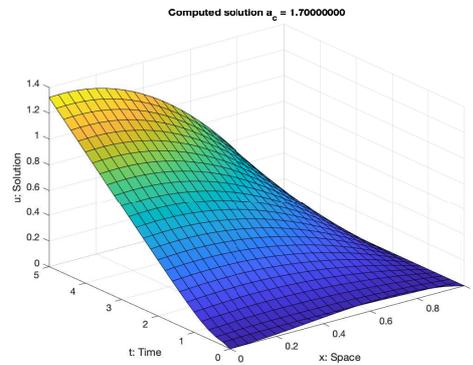}
\captionof{figure}{Test~\ref{test.const33BIS}, computed solution $u=u(x,t)$.}
\label{fig.sol_nh}
\end{minipage}

%

\subsection{Numerical tests for the linear case II}\label{sec.numerics2}

In this section we will present some numerical results for second inverse problem of identification of $a$ from one boundary observation $u_x(1,t)$, mentioned at the beginning of the Section~\ref{sec.lincase}. Arguing as in the previous section, we reformulate the inverse problem as an optimization problem  \eqref{pb.const11}, where the funcional $J : a\in \mathcal{U}^{a}_{ad} \mapsto \mathds{R}$ is now given by
\begin{equation}\label{J.const}
J(a) = \dfrac{1}{2}\displaystyle\int_0^T |\eta(t) - \partial_x u^a(1,t)|^2\, dt, 
\end{equation}
where $\eta =\partial_x u(1,\cdot)$ is given and $u^a$ is the solution of \eqref{dpe1} corresponding to the unknown $a$. 

In what follows we present some numerical tests for this inverse problem. We will consider the same data as in the previous section: $T=5$, $\ell=1$, $u_0(x) = 0.5x^2(1-x)$ and $f(x,t) = 0$. In order to simplify the presentation, we will show just the evolution of the cost and the results with random noises. 

\bigskip
\begin{test}\label{test5}
We take $a_{i} = 0.7$ as an initial guess for the minimization algorithm. The desired value is $a_d=1.7$. The results of the numerical computations are presented in Table~\ref{table.t5} and Figure~\ref{fvalues_t5}. 
\end{test}

\begin{minipage}[t]{0.49\linewidth}
\centering
\renewcommand{\arraystretch}{1.2}
\medskip
\captionof{table}{Evolution of the cost and computed $a_c$ with random noises in the target in Test~\ref{test5}.}
\medskip
{\begin{tabular}{cccc} \hline
\%  noise  & Cost  & Iterations &  $a_c$ \\
\hline
1\%           &   1.e-6       &    13       &  1.69198034   \\ \hline
0.1\%        &    1.e-8    &      13        &  1.70001564 \\ \hline
0.01\%      &   1.e-14     &    10       &   1.70009682\\\hline
0.001\%    &   1.e-18     &    12       &   1.70000143\\ \hline
0\%           &   1.e-25     &    18        &  1.7   \\  \hline
\end{tabular}}
\label{table.t5}
\end{minipage}
\hfill
\begin{minipage}[t]{0.49\linewidth}
\vspace{0cm}
\centering
\includegraphics[width=0.9\linewidth]{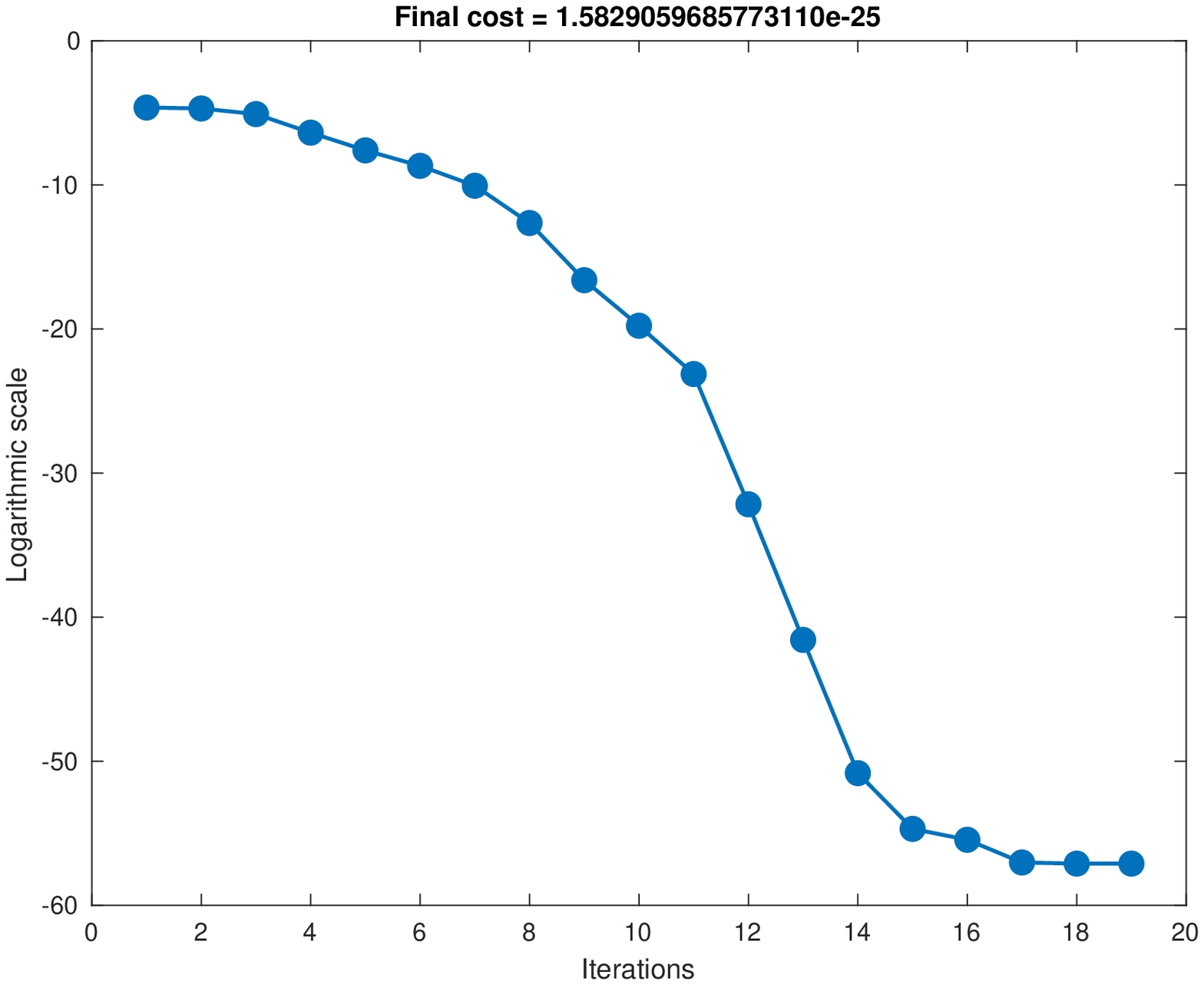}
 \captionof{figure}{Test~\ref{test5}, the evolution of the cost.}
\label{fvalues_t5}
\end{minipage}

\begin{test}\label{test6}
Taking $a_{i} = 0.2$ and  $a_d=1$, we can see the numerical results in Figure~\ref{fvalues_t6} and Table~\ref{table.t6}. 
\end{test}

\begin{minipage}[t]{0.49\linewidth}
\centering
\renewcommand{\arraystretch}{1.2}
\medskip
\captionof{table}{Evolution of the cost and computed $a_c$ with random noises in the target in Test~\ref{test6}.}
\medskip
{\begin{tabular}{cccc} \hline
\% noise  & Cost  & Iterates &   $a_c$ \\
\hline
1\%           &   1.e-6       &    10         &  0.99701925  \\ \hline
0.1\%        &   1.e-8      &    13         &  1.00058988\\ \hline
0.01\%      &   1.e-9       &    16         &   0.99994514\\\hline
0.001\%    &   1.e-13     &    12         &   0.99999587\\ \hline
0\%           &   1.e-27     &    13         &  1   \\  \hline
\end{tabular}}
\label{table.t6}
\end{minipage}
\hfill
\begin{minipage}[t]{0.49\linewidth}
\vspace{0cm}
\centering
\includegraphics[width=0.9\linewidth]{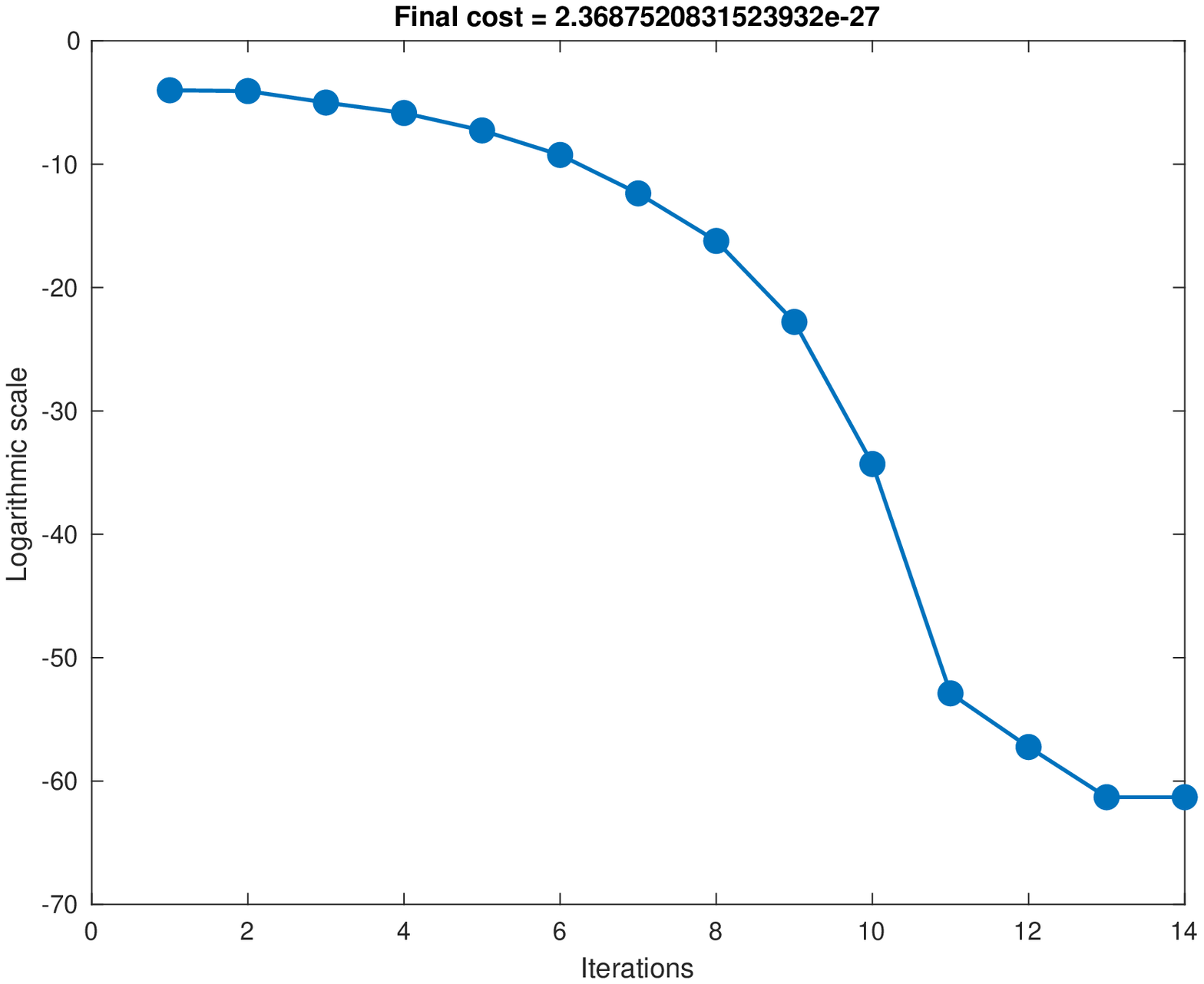}
 \captionof{figure}{Test~\ref{test6}, the evolution of the cost.}
\label{fvalues_t6}
\end{minipage}
\bigskip

\begin{test}\label{test7}
Taking $a_{i} = 0.7$ and $a_d=0.2$, Figure~\ref{fvalues_t7} and Table~\ref{table.t7} show the numerical results. 
\end{test}

\begin{minipage}[t]{0.49\linewidth}
\centering
\renewcommand{\arraystretch}{1.2}
\medskip
\captionof{table}{Evolution of the cost and computed $a_c$ with random noises in the target in Test~\ref{test7}.}
\medskip
{\begin{tabular}{cccc} \hline
\%  noise  & Cost  & Iterations &   $a_c$ \\
\hline
1\%           &   1.e-5      &    15      &   0.20127090   \\ \hline
0.1\%        &    1.e-8     &    15       &  0.20020541\\ \hline
0.01\%      &   1.e-9      &    14       &   0.19997384\\\hline
0.001\%    &   1.e-11     &    14       &   0.20000100\\ \hline
0\%           &   1.e-28     &   18        &  0.2   \\  \hline
\end{tabular}}
\label{table.t7}
\end{minipage}
\hfill
\begin{minipage}[t]{0.49\linewidth}
\vspace{0cm}
\centering
\includegraphics[width=0.9\linewidth]{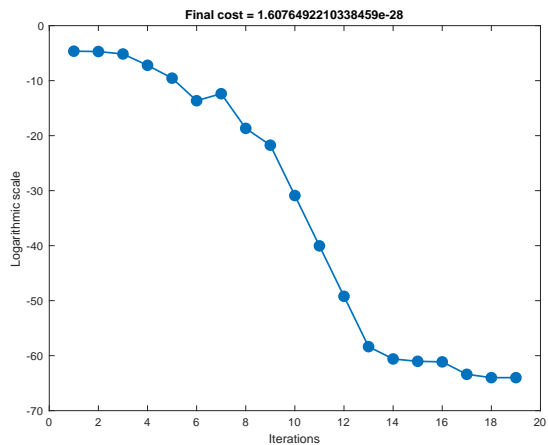}
 \captionof{figure}{Test~\ref{test7}, the evolution of the cost.}
\label{fvalues_t7}
\end{minipage}


\subsection{Numerical tests for the power-like case I}\label{sec:numpower1}

We will present in this section some numerical results of the reconstruction of $\alpha$ for weak and strong degenerate cases \eqref{eq.CP.w} and \eqref{eq.CP} from data that appear in Theorem~\ref{th.power}. More precisely, given $t_0\in (0,T)$ as in \eqref{eq.h4}, $u_0=u_0(x)$ and $\ell<1$, in order to reconstruct $\alpha$ we will solve numerically the following optimization problem: find $\alpha\in (0,2)$ such that 
\begin{equation}\label{cpbalpha}
J(\alpha) \le  J(\alpha'),\quad \forall\, \alpha'\in (0,2),
\end{equation}
where the functional $J : (0,2) \mapsto \mathds{R}$ is given by
\begin{equation*}
J(\alpha) = \dfrac{1}{2}\displaystyle\int_0^1 |\gamma(x) - \partial_t u^\alpha(x,t_0)|^2\, dx 
+  \dfrac{1}{2}\displaystyle\int_0^1 |\beta(x) - x^2\partial_x u^\alpha(x,t_0)|^2\, dx,
\end{equation*}
with $u^\alpha$ the corresponding solution of  \eqref{eq.CP.w} or  \eqref{eq.CP} assoceated to the unknown power $\alpha$. 

\medskip
 In the numerical tests we will take $T=10$, $u_0=0.3x^2(1-x)^2$  and $f(x,t) =0$ and we will use the \texttt{interior-point} MatLab gradient free algorithm from the Optimization Toolbox.

\subsubsection{Weak degeneracy I}

\begin{figure}[h!]
\bigskip
\begin{minipage}[t]{0.49\linewidth}
\centering
\includegraphics[width=0.9\linewidth]{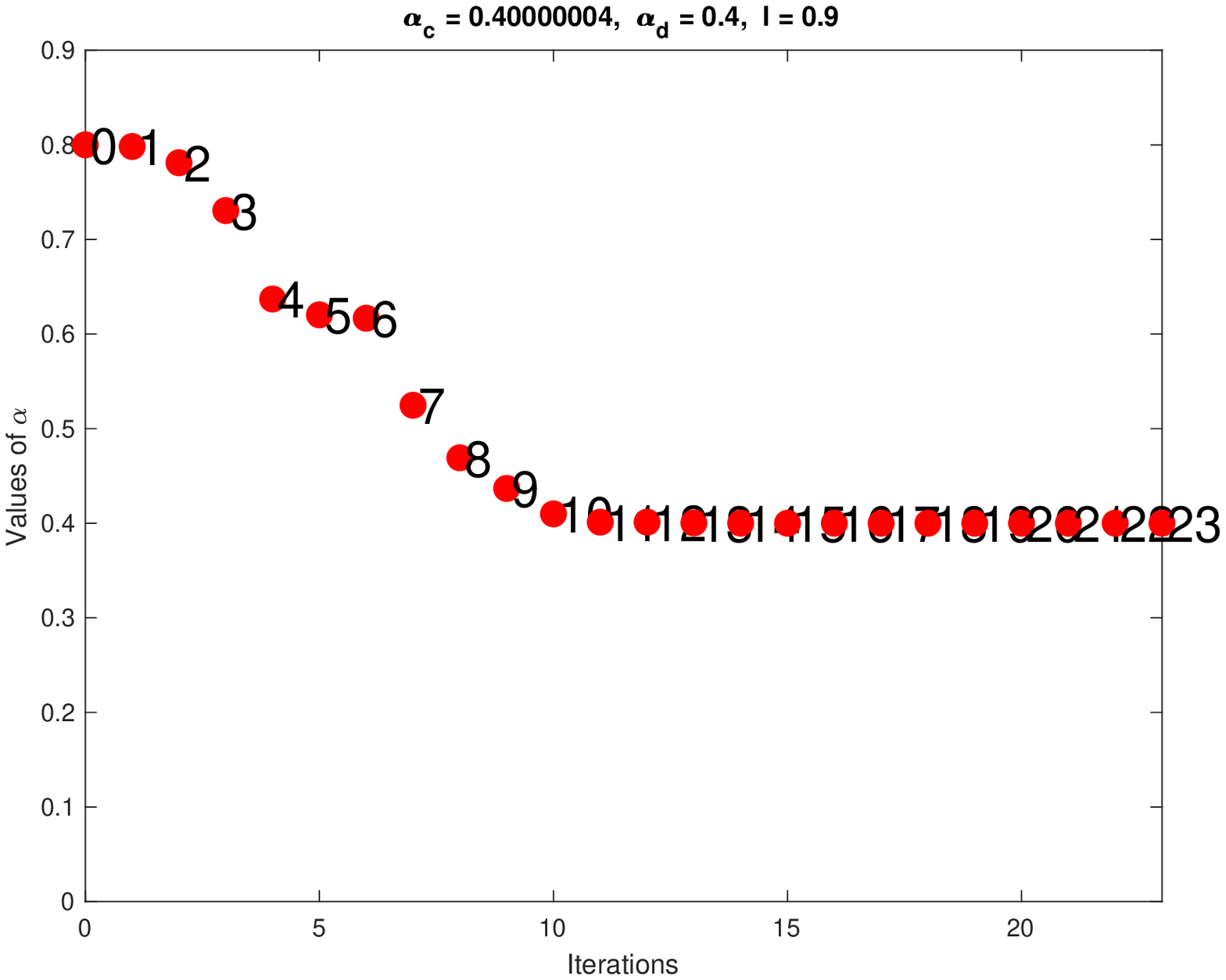}
\caption{Test~\ref{test8}, identification of $\alpha$, $\ell=0.9$.}
\label{alpha07w}
\end{minipage}
\hfill
\begin{minipage}[t]{0.49\linewidth}
\centering
\includegraphics[width=0.9\linewidth]{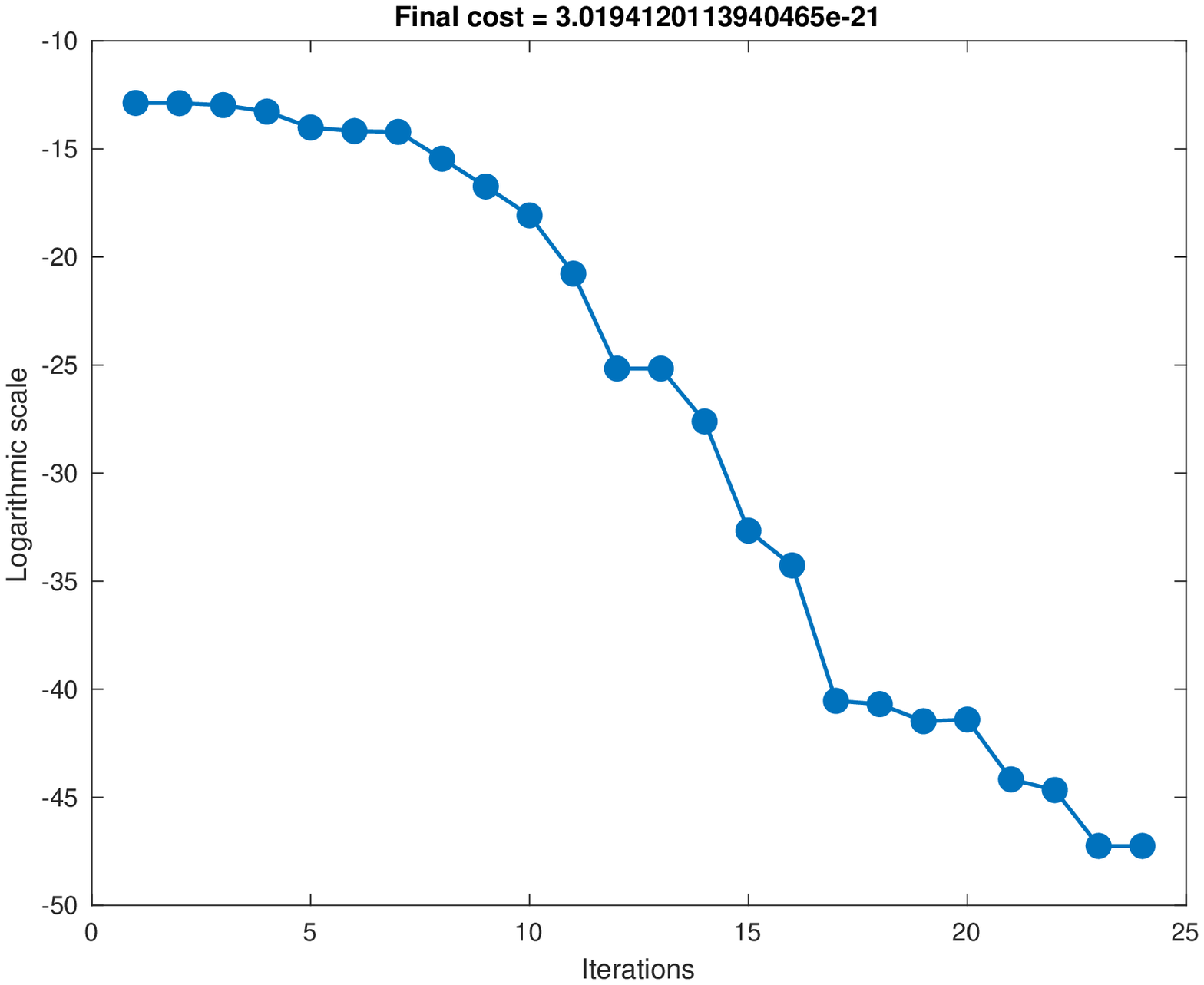}
\caption{Test~\ref{test8}, current function values.}
\label{fvalalpha07w}
\end{minipage}
\end{figure}

\begin{test}\label{test8} Let us  first assume that $0<\alpha < 1$, corresponding to the weak degeneracy case \eqref{eq.CP.w}. 
We take $\alpha_{i} = 0.8$ as an initial guess for the minimization algorithm and the desired value to recover is $\alpha_d=0.4$. In Figure~\ref{alpha07w} we can see the computed values of $\alpha$ at each iteration of the optimization algorithm and  Figure~\ref{fvalalpha07w} shows the evolution of the cost. In the Table~\ref{table.w1}, we present the computed values of $\alpha_c$ depending on the values of $\ell$ close to 1.   
\end{test}

\begin{table}[h!]
\centering
\renewcommand{\arraystretch}{1}
\caption{Computed $\alpha_c$ for $\ell$ close to 1 in Test~\ref{test8}.}
\medskip
{\begin{tabular}{lccl | lccl } \hline
$\ell$  & Cost  & Iterations &   $\alpha_c$ &  $\ell$  & Cost  & Iterations &   $\alpha_c$ \\
\hline
0.9            &   1.e-21        &    23          &    0.40000004 & 
1.00001     &   1.e-19       &    25           &    0.40000010
\\ \hline
0.99           &   1.e-22       &    27           &    0.4   & 
1.0001         &   1.e-23      &    28        &       0.4
\\ \hline
0.999           &   1.e-22      &    32           &    0.39999999 & 
1.001             &   1.e-22      &    31         &    0.39999999
\\ \hline
0.9999           &   1.e-21      &    26         &    0.39999999 & 
1.01               &   1.e-22      &    22        &    0.4 
\\ \hline
0.99999        &   1.e-21      &    30          &     0.4  & 
1.1                &   1.e-12      &    30          &    0.40053247 
\\ \hline
1               &   1.e-19        &    27          &    0.40000008 & 
1.2                  &   1.e-22        &    30          &    0.4 
\\ \hline
\end{tabular}}
\label{table.w1}
\end{table}

\eject

\subsubsection{Strong degeneracy I}

\begin{figure}[h!]
\begin{minipage}[t]{0.49\linewidth}
\centering
\includegraphics[width=0.9\linewidth]{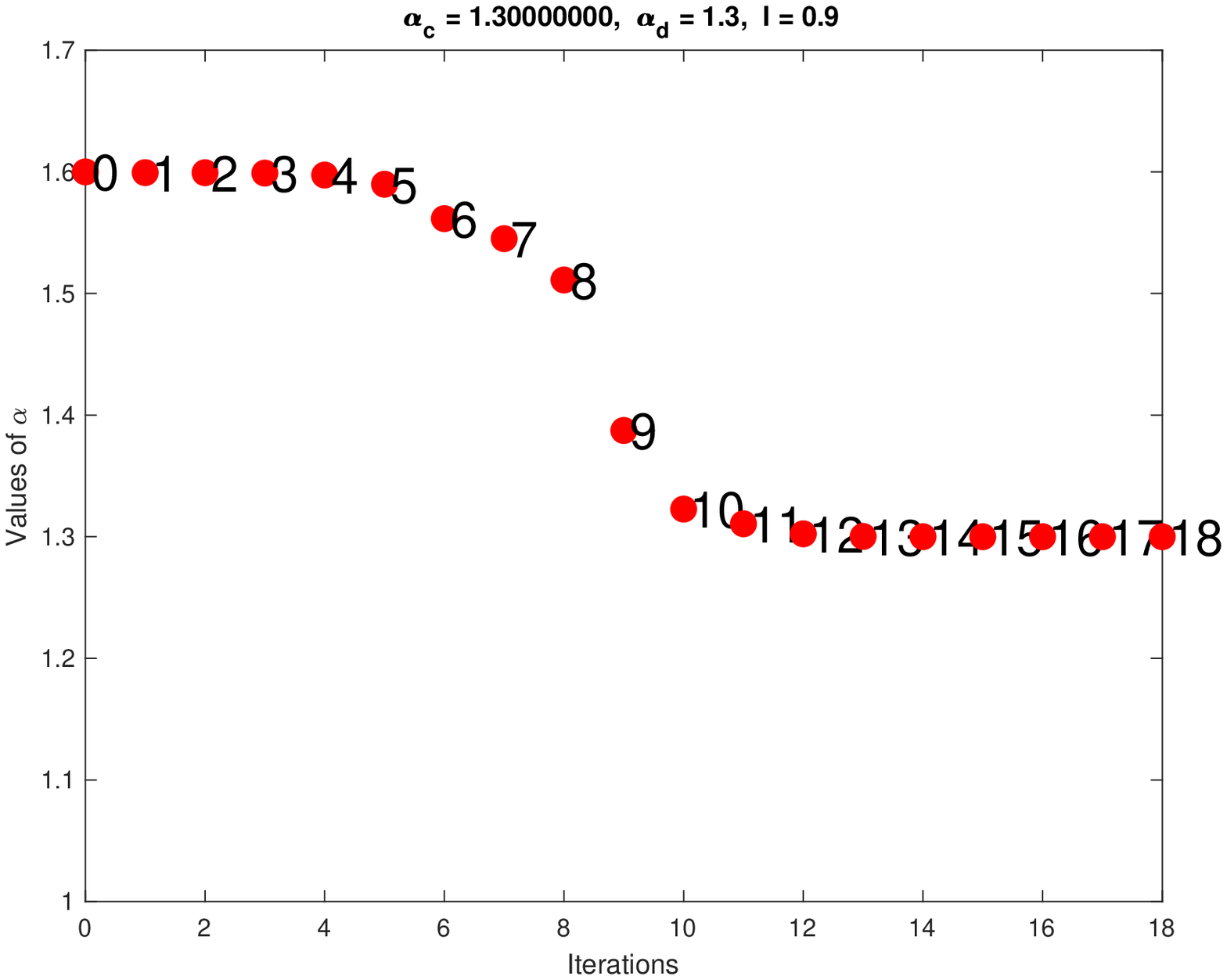}
\caption{Test~\ref{test9}, identification of $\alpha$, $\ell=0.9$.}
\label{alpha13s}
\end{minipage}
\hfill
\begin{minipage}[t]{0.49\linewidth}
\centering
\includegraphics[width=0.9\linewidth]{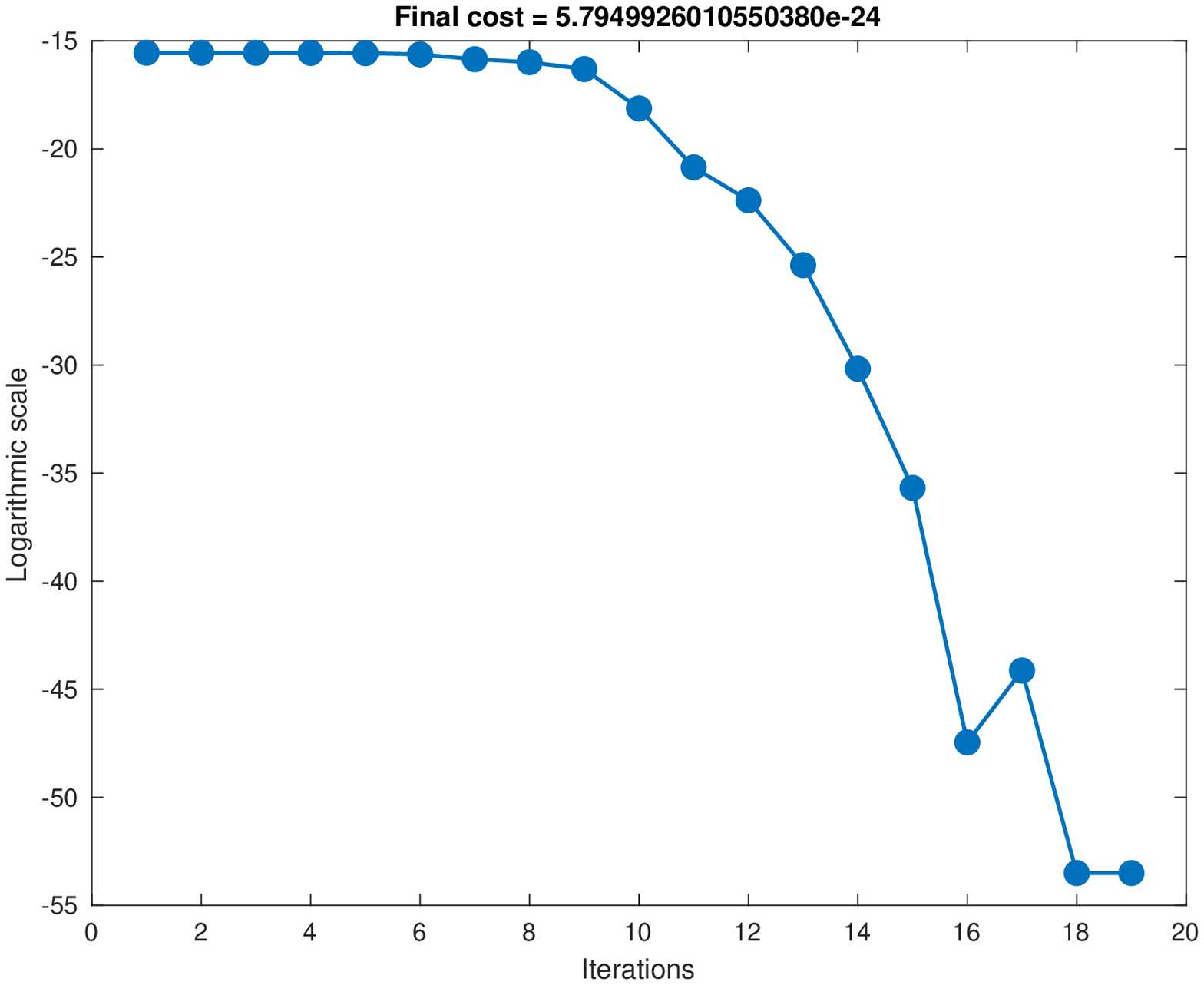}
\caption{Test~\ref{test9}, current function values.}
\label{fvalalpha13s}
\end{minipage}
\end{figure}

\begin{test}\label{test9} Let us consider the strong degeneracy case with $1\leq \alpha <2$, i.e. the problem \eqref{eq.CP}.  We take $\alpha_{i} = 1.6$ and $\alpha_d=1.3$. The numerical results can be seen in the Figures~\ref{alpha13s} and~\ref{fvalalpha13s}. In the Table~\ref{table.sd1}, we present the computed values of $\alpha_c$ depending on the values of $\ell$ close to 1.   
\end{test}

\begin{table}[h!]
\centering
\renewcommand{\arraystretch}{1}
\caption{Computed $\alpha_c$ for $\ell$ close to 1 in Test~\ref{test9}.}
\medskip
{\begin{tabular}{lccl | lccl } \hline
$\ell$  & Cost  & Iterations &   $\alpha_c$ &  $\ell$  & Cost  & Iterations &   $\alpha_c$ \\
\hline
0.9            &   1.e-24        &    18          &    1.3 & 
1.00001     &   1.e-22       &    23          &    1.30000001
\\ \hline
0.99           &   1.e-22       &    24           &    1.30000001   & 
1.0001         &   1.e-23     &    21        &       1.3
\\ \hline
0.999           &   1.e-21      &    19           &    1.30000006 & 
1.001           &   1.e-23       &    25         &    1.3 
\\ \hline
0.9999           &   1.e-23      &    22         &    1.29999999 & 
1.01               &   1.e-24      &    32        &     1.3 
\\ \hline
0.99999        &   1.e-22      &    30          &     1.30000001  & 
1.1                &   1.e-22      &    18          &    1.29999998 
\\ \hline
1                  &   1.e-23        &    27          &    1.30000000 & 
1.2               &   1.e-22        &    20          &    1.3
\\ \hline
\end{tabular}}
\label{table.sd1}
\end{table}

\begin{remark}
From the results of the Tests~\ref{test8} and~\ref{test9}  we can see numerically that the hypothesis  $\ell<1$ from Theorem~\ref{th.power} probably is not really needed for the reconstruction of $\alpha$. But the question whether we can eliminate this condition form this theorem remains open. 
\end{remark}

\subsection{Numerical tests for the power-like case II}\label{sec:powenum2} 

We will present in this section numerical results for the reconstruction of $\alpha$ from one observation of the form $\partial_x u(1,t)$. As before, we will analyze first weak and then strong degenerate cases.  More precisely, given $\ell$, $u_0=u_0(x)$ and $\eta=\partial_x u(1,\cdot)$, we will solve numerically the following optimization problem \eqref{cpbalpha}, but wiht the functional $J : (0,2) \mapsto \mathds{R}$ given by
\begin{equation*}
J(\alpha) = \dfrac{1}{2}\displaystyle\int_0^T |\eta(t) - \partial_x u^\alpha(1,t)|^2\, dt,
\end{equation*}
where $u^\alpha$ is the solution to \eqref{eq.CP.w} or  \eqref{eq.CP}  corresponding to $\alpha$.

\medskip
In that follows, we will take $T=10$, $\ell=1$, $u_0=0.3x^2(1-x)^2$ and $f(x,t)=0$. The computations are performed using 
\texttt{fmincon}  MatLab function form the Optimization Toolbox, combined with \texttt{interior-point} gradient algorithm.

\subsubsection{Weak degeneracy II}

\begin{figure}[h!]
\begin{minipage}[t]{0.49\linewidth}
\noindent
\centering
\includegraphics[width=0.9\linewidth]{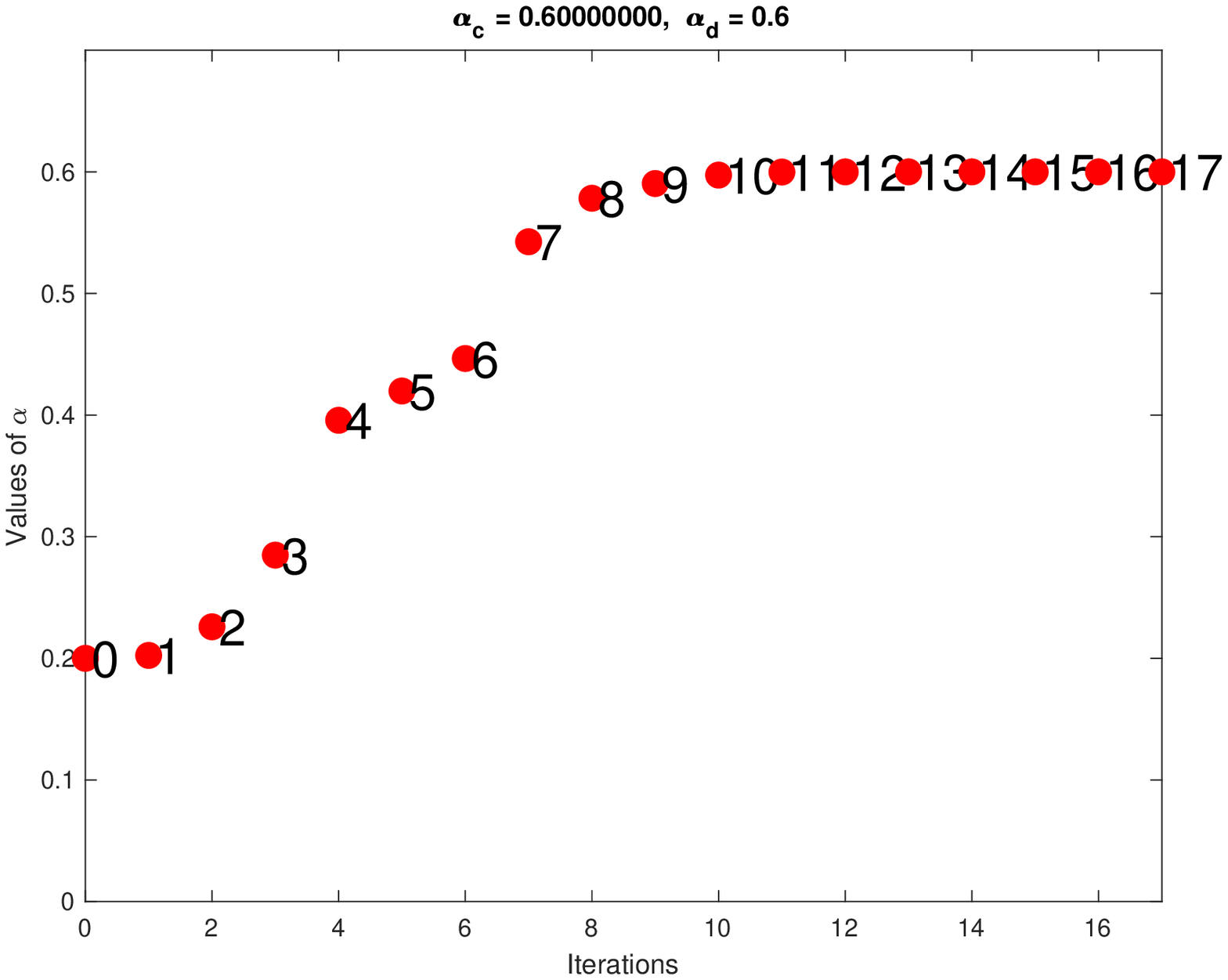}
\caption{Test~\ref{test10}, identification of $\alpha$.}
\label{alpha1B}
\end{minipage}
\hfill
\begin{minipage}[t]{0.49\linewidth}
\noindent
\centering
\includegraphics[width=0.9\linewidth]{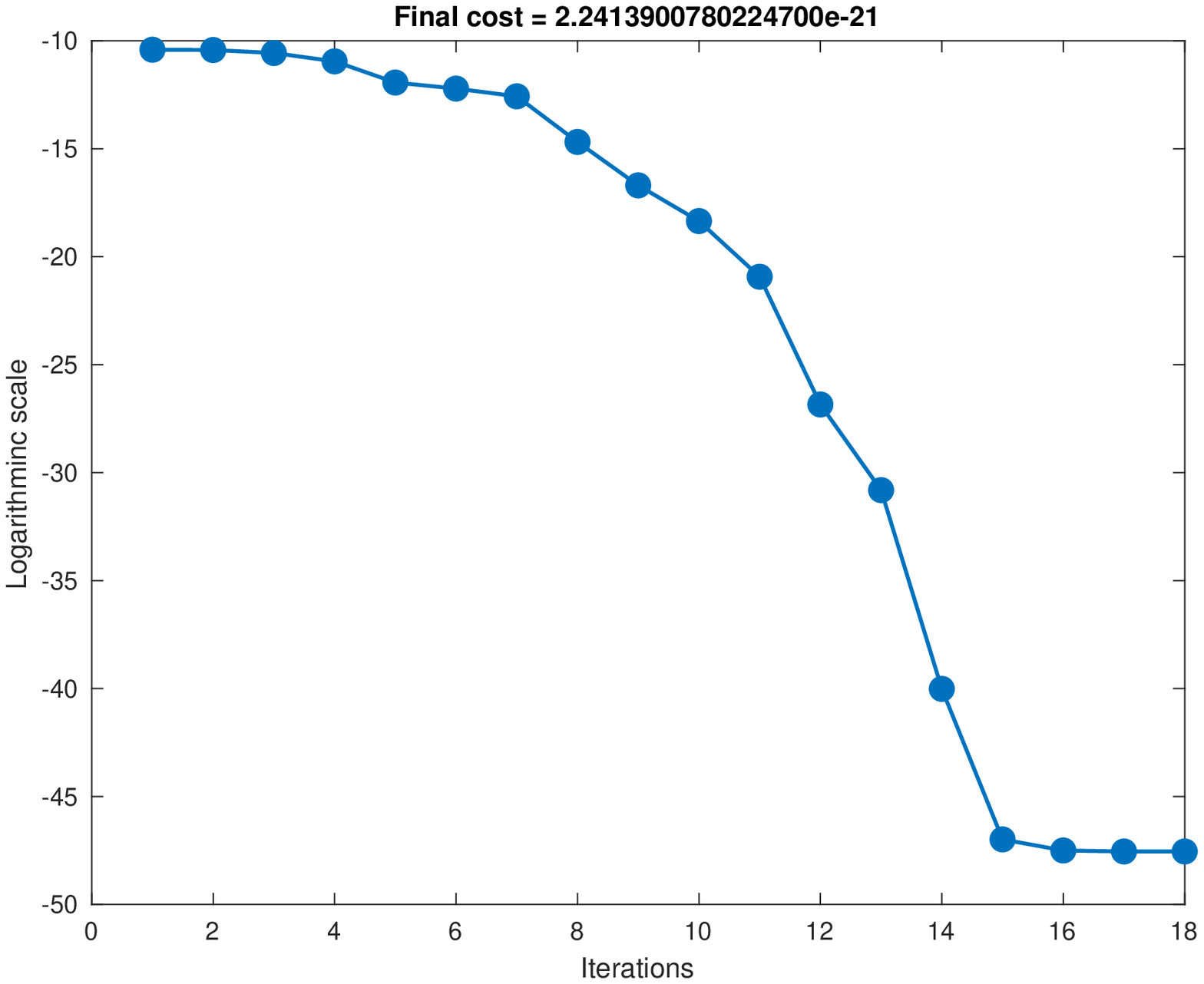}
\caption{Test~\ref{test10}, current function values.}
\label{fvalalpha1B}
\end{minipage}
\end{figure}

\begin{test}\label{test10}
Let us first consider the weak degeneracy case corresponding to the problem \eqref{eq.CP.w} with  $0<\alpha < 1$. 
We take $\alpha_{i} = 0.2$ as an initial guess for the minimization algorithm and the desired value to recover $\alpha_d=0.6$. The numerical results can be seen in Figures~\ref{alpha1B} and~\ref{fvalalpha1B}. Table~\ref{table.t10} shows the computations of $\alpha_c$ and the cost with random noises in the target.  
\end{test}

\bigskip

\begin{minipage}[t]{0.49\linewidth}
\centering
\renewcommand{\arraystretch}{1.2}
\medskip
\captionof{table}{Evolution of the cost and computed $\alpha_c$ with random noises in the target in Test~\ref{test10}.}
\medskip
{\begin{tabular}{cccc} \hline
\% noise  & Cost  & Iterations &  $\alpha_c$ \\
\hline
1\%           &   1.e-8       &    15         &  0.59689497  \\ \hline
0.1\%        &    1.e-10      &    16         &  0.60037182\\ \hline
0.01\%      &   1.e-12       &    16         &   0.59997672\\\hline
0.001\%    &  1.e-14      &    18           &   0.59999807\\ \hline
0\%           &   1.e-21    &    17            &  0.6 \\  \hline
\end{tabular}}
\label{table.t10}
\end{minipage}
\hfill
\begin{minipage}[t]{0.49\linewidth}
\vspace{0cm}
\includegraphics[width=0.9\linewidth]{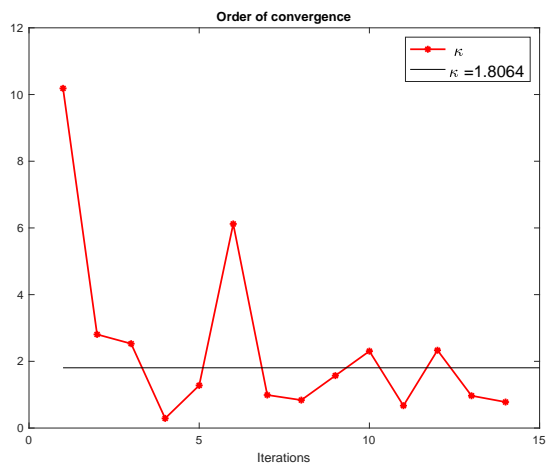}
 \captionof{figure}{Test~\ref{test10}, order of convergence.}
\label{fig.order}
\end{minipage}

\begin{remark}[On the order of convergence] We can  analyze the order of convergence of $\alpha_n$ which we obtain at each iteration $n$ of the optimization algorithm to $\alpha_d$. More precisely, we look for positive real numbers $C$ and $\kappa$  such that 
\[
|e_{n+1}| \approx C|e_n|^\kappa, 
\]
where $e_n = \alpha_n-\alpha_d$. From Figure~\ref{fig.order}, we see that this approximate formula is satisfied for $\kappa = 1.8$ for some finite value of $C$. 

\end{remark}

\subsubsection{Strong degeneracy II}

We will consider here the identification of $\alpha$ in the case $1\leq \alpha <2$ for strong degeneracy case, that is to say for the problem \eqref{eq.CP}. For simplicity of the presentation, the results with random noises in the targets are not presented, because they are similar to those of the previous sections.  

\begin{test}\label{test.id1}
Let us take $\ell = 0.99$ and $\alpha_{i} = 1.6$ as an initial guess for the minimization algorithm. The desired value of $\alpha$ to recover is $\alpha_d=1.3$. The numerical results can be seen in the Figures~\ref{alpha02} and~\ref{fvalalpha02} and the computed value of $\alpha$ is $\alpha_c=1.3$. 
\end{test}

\begin{figure}[h!]
\begin{minipage}[t]{0.49\linewidth}
\noindent
\centering
\includegraphics[width=0.9\linewidth]{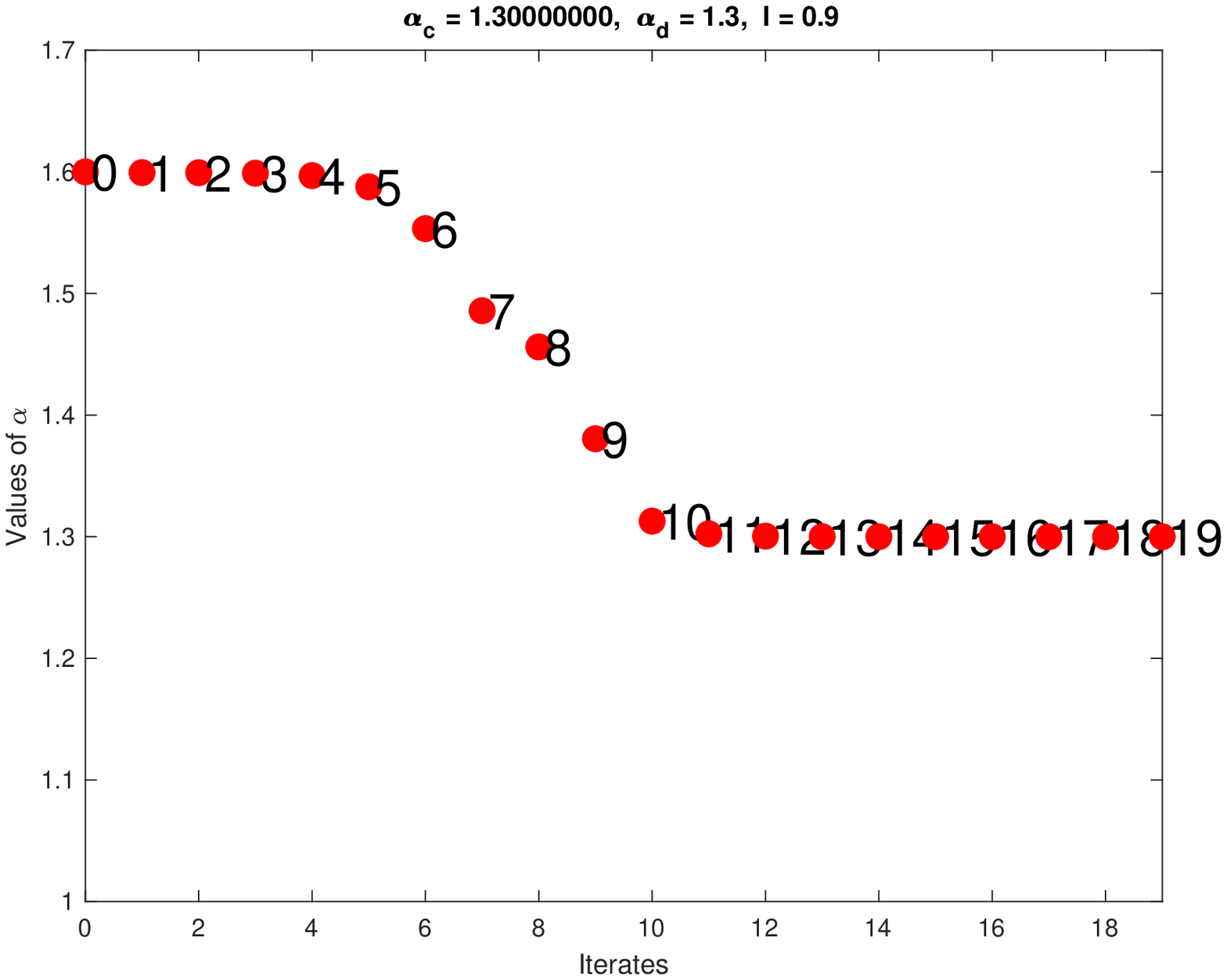}
\caption{Test~\ref{test.id1}, identification of $\alpha$.}
\label{alpha02}
\end{minipage}
\hfill
\begin{minipage}[t]{0.49\linewidth}
\noindent
\centering
\includegraphics[width=0.9\linewidth]{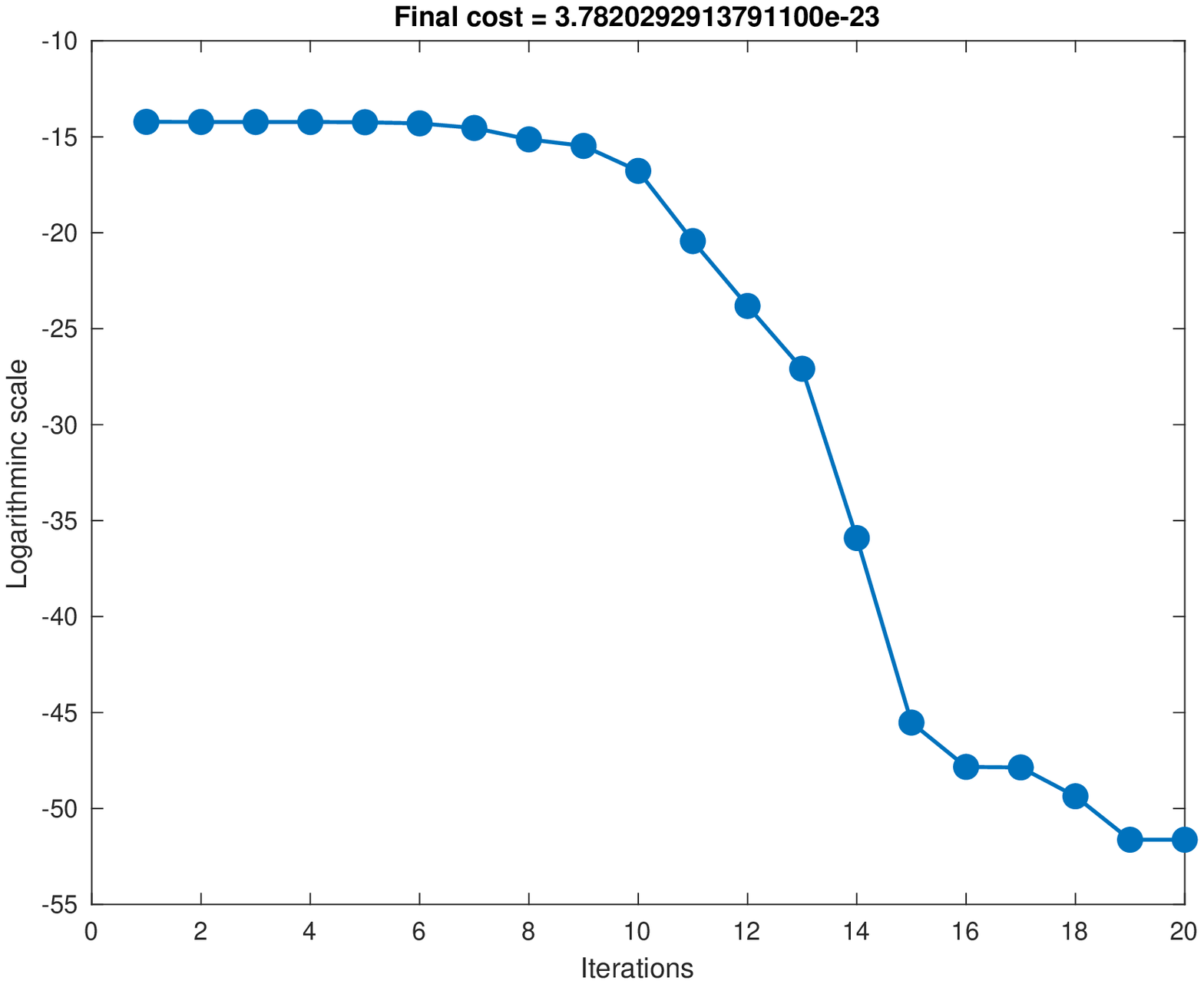}
\caption{Test~\ref{test.id1}, current function values.}
\label{fvalalpha02}
\end{minipage}
\end{figure}


\begin{test}\label{test.id1}
Let us take $\ell = 1$ and $\alpha_{i} = 1.6$ and $\alpha_d=1.3$. The numerical results can be seen in Figures~\ref{alpha02} and~\ref{fvalalpha02} and the computed value is $\alpha_c=1.3$. 
\end{test}

\begin{figure}[h!]
\begin{minipage}[t]{0.49\linewidth}
\noindent
\includegraphics[width=0.9\linewidth]{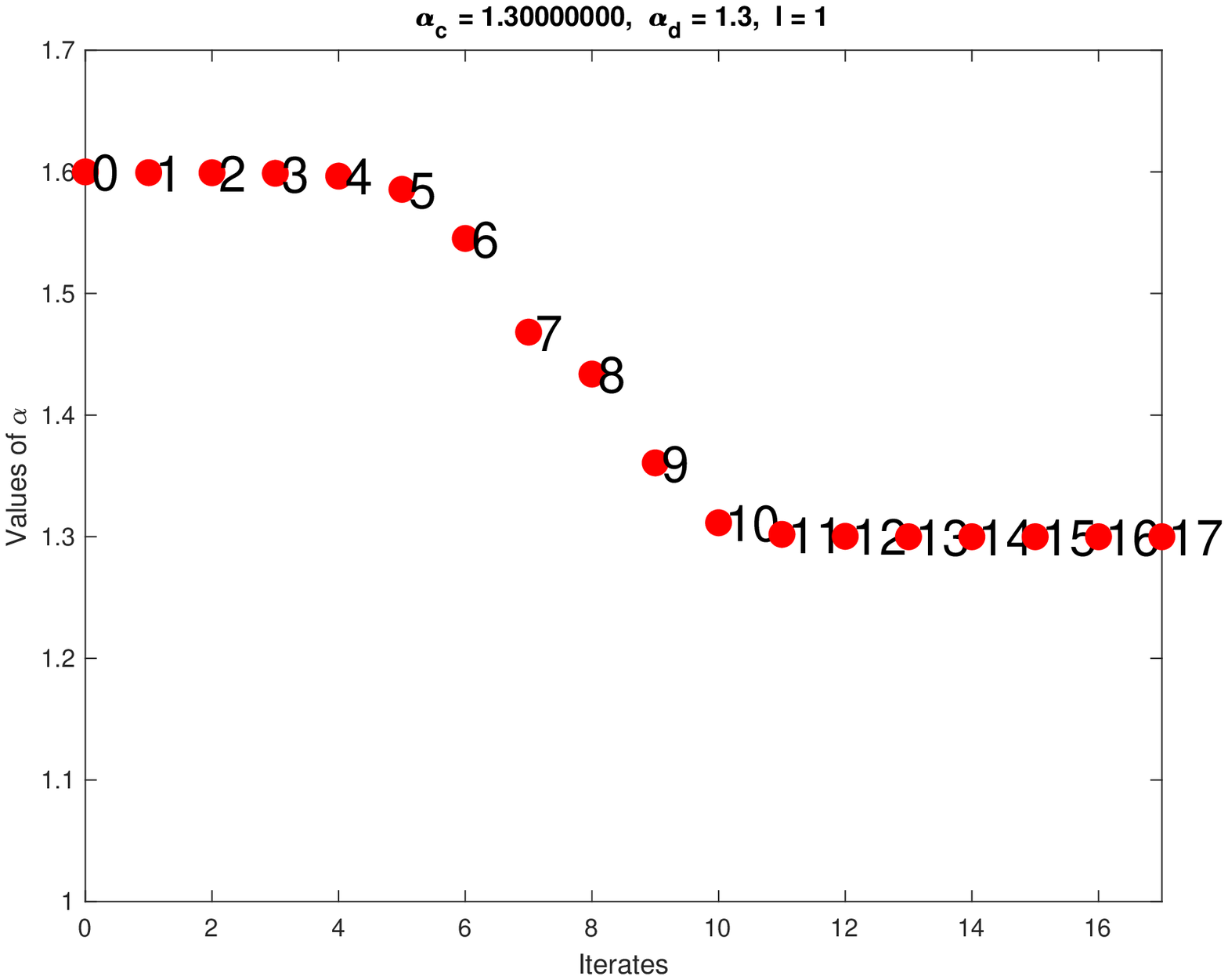}
\caption{Test~\ref{test.id1}, identification of $\alpha$.}
\label{alpha02}
\end{minipage}
\hfill
\begin{minipage}[t]{0.49\linewidth}
\noindent
\includegraphics[width=0.9\linewidth]{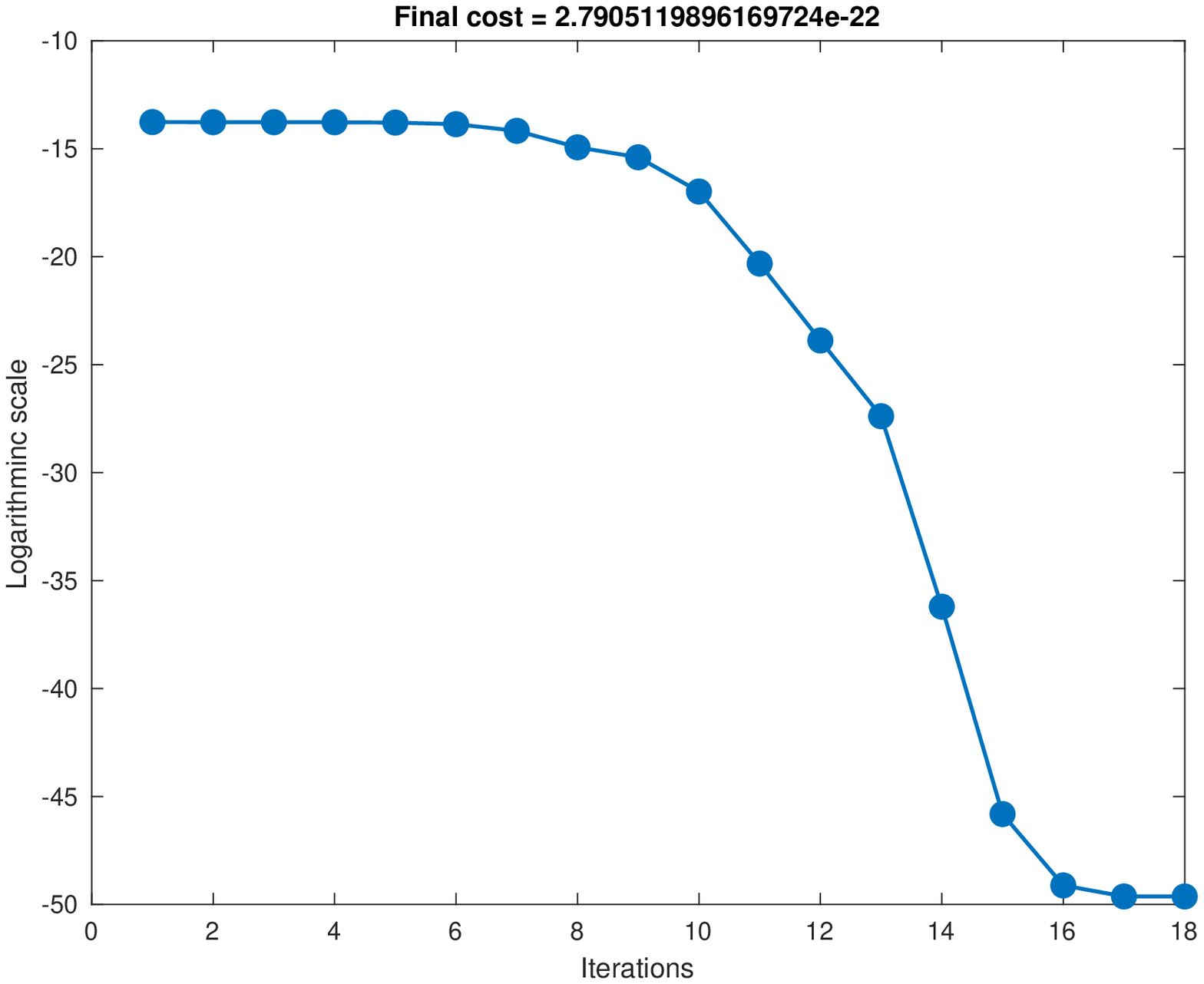}
\caption{Test~\ref{test.id1}, current function values.}
\label{fvalalpha02}
\end{minipage}
\end{figure}


\subsection{Numerical tests for a general case}

In this section we will present numerical results of identification of a function $a=a(x)$ in \eqref{1.1} from $\partial_x u(1,t)$ in order to illustrate the theoretical part form the Section~\ref{sec.gencase}. 

\medskip
Let us take  $\ell =1$, $T = 5$, $\alpha = 0.6$ and $u_0(x) =  0.3x^2(1-x)^2$. In order to reconstruct $a=a(x)$,  again we reformulate the corresponding inverse problem as an optimization problem:  find $a\in \mathcal{A}$ such that 
\begin{equation}\label{optpb1}
J(a) \le  J(a'),\quad \forall\, a'\in\mathcal{A}\,,
\end{equation}
where the functional $J : \mathcal{A} \mapsto \mathds{R}$ is given by
\begin{equation*}\label{J1}
J(a) = \dfrac{1}{2}\displaystyle\int_0^T |\eta(t) - \partial_x u^a(1,t)|^2\, dt,
\end{equation*}
where $u^a$ is the solution of \eqref{1.1} corresponding to $a$ and $\mathcal{A}$ is given by~\eqref{27}. In fact, for the numerical reconstruction we do not need to assume that 
$a$ is known near 1 as in \eqref{27}.
%

\medskip
In order to present some examples, in what follows we will analyze two particular cases of the coefficient $a(x)$: linear and quadratic functions. The numerical experiments will be performed using a gradient based  \texttt{interior-point} algorithm of \texttt{fmincon} MatLab Optimization Toolbox. 

\medskip
Let us notice that we can also reconstruct numerically both $\alpha$ and $a$ from one measurement, but in order to simplify the presentations of the results, we will omit the details.

\begin{description}
\item[Linear function:] Let us consider $a$ of the form $a(x) = bx + c$, where $b$ and $c$ are unknown constants. The optimization problem \eqref{optpb1} is written in this case as follows: find  $b,c\in\mathcal{A}$, with  $\mathcal{A} = \{(b,c) : 0<\underline{b} < b<\overline{b}, \quad  0<\underline{c} < c<\overline{c}\}$ such that 
\begin{equation*}\label{optpb2}
J(b,c) \le J(b',c'), \quad \forall\, (b',c')\in \mathcal{A}.
\end{equation*}
\end{description}

\begin{figure}[h!]
\begin{minipage}[t]{0.49\linewidth}
\noindent
\centering
\includegraphics[width=0.9\linewidth]{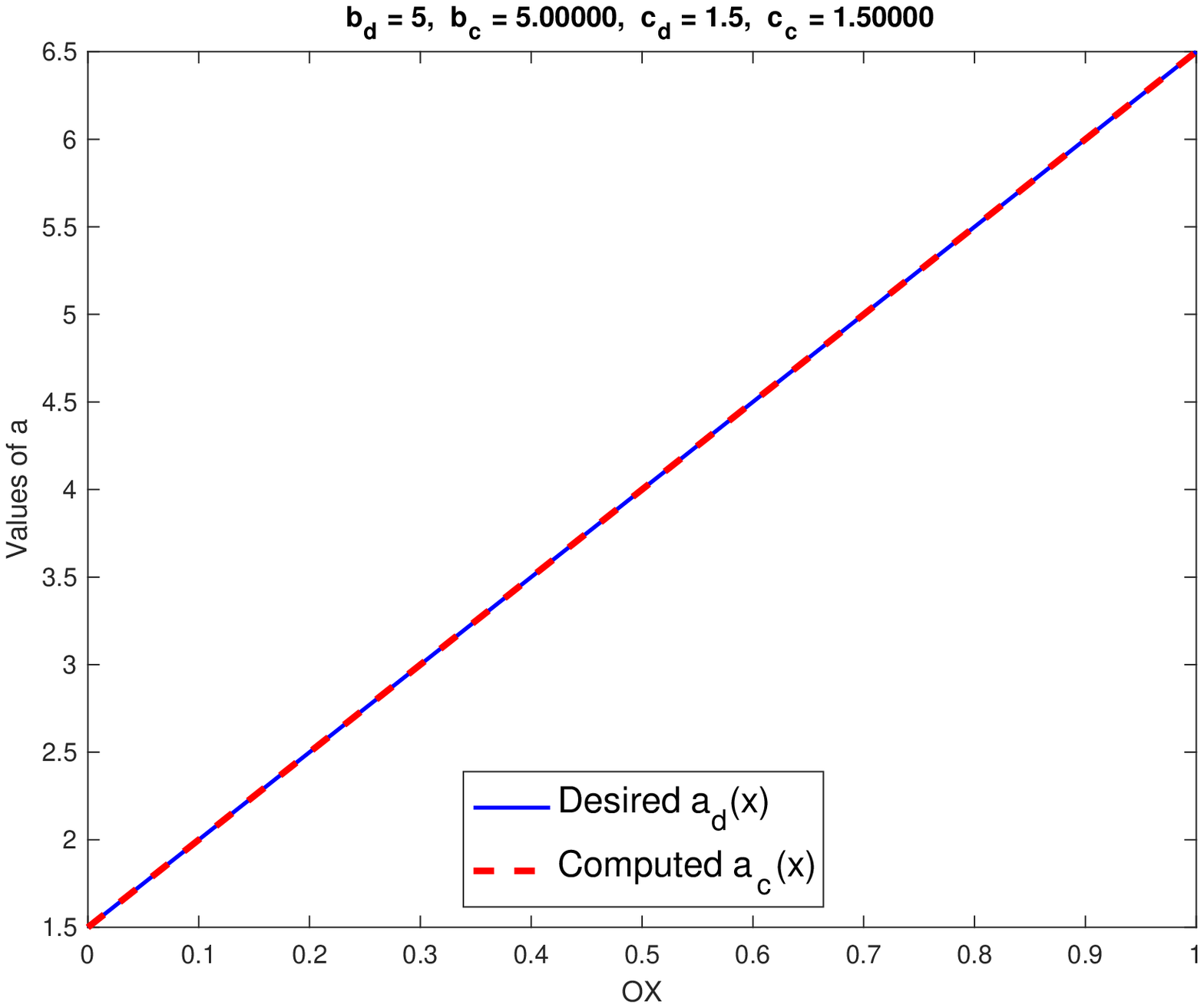}
\caption{Test~\ref{testLin}, $a_d$ solid line, $a_c$ dashed line.}
\label{fig.Test1}
\end{minipage}
\hfill
\begin{minipage}[t]{0.49\linewidth}
\noindent
\centering
\includegraphics[width=0.9\linewidth]{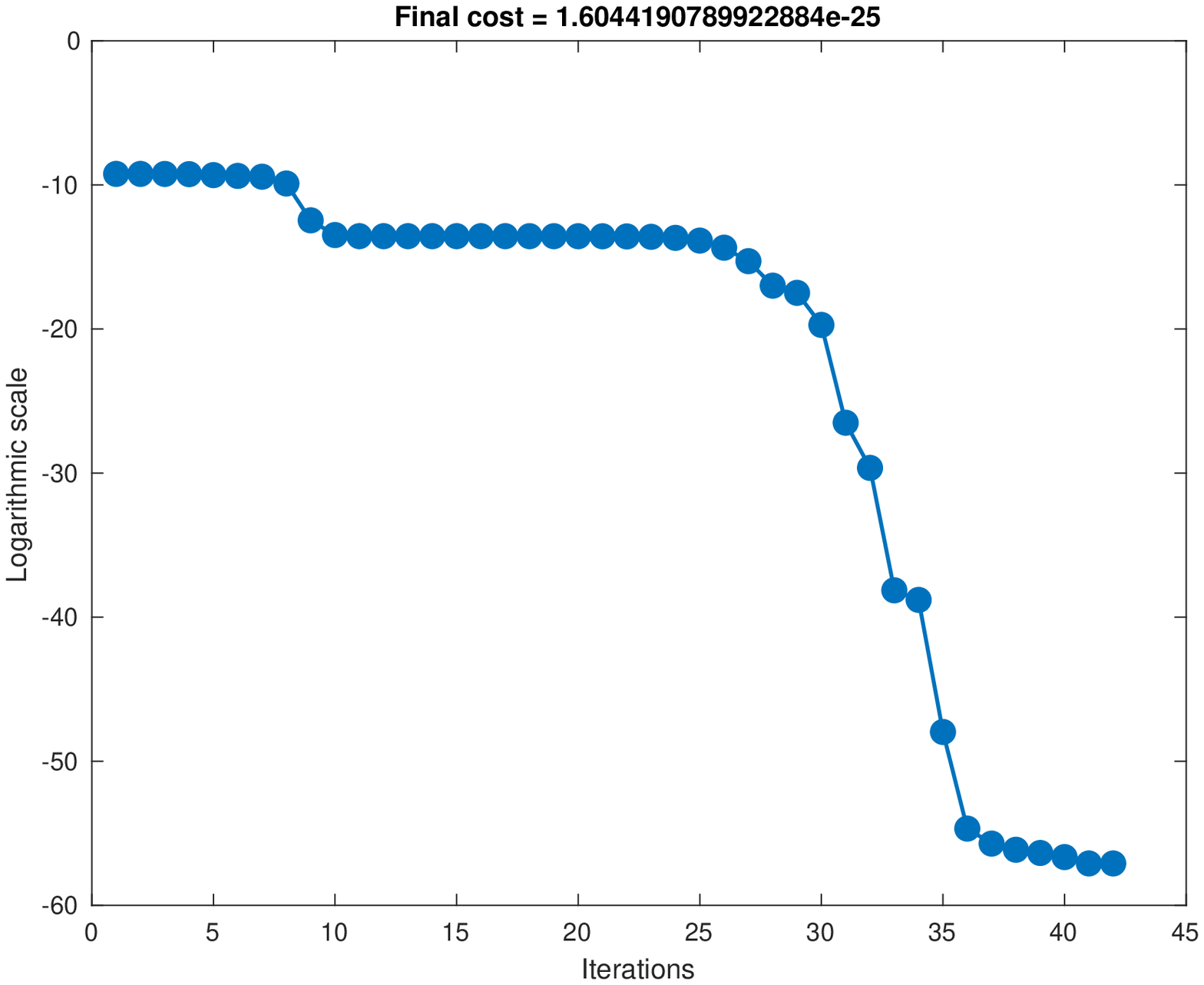}
\caption{Test~\ref{testLin}, evolution of the cost. }
\label{fig.Test2}
\end{minipage}
\end{figure}

\begin{test}\label{testLin}
Let us take $b_i = 1$ and $c_i = 1$ as initial guess for the optimization algorithm and desired values  to recover  $b_d=5$  and $c_d=1.5$. After numerical approximation of the optimization problem above we get the computed values $b_c=5$ and $c_c=1.5$.  The numerical results we have obtained are in Figures~\ref{fig.Test1} and~\ref{fig.Test2}. Table~\ref{table.t11} shows the evolution of the cost during the iterations of optimization algorithm and computed $b_c$  and $c_c$ with random noises in the target. Finally, in Figure~\ref{Sol_Lin} the solution of \eqref{1.1} is represented for the computed values of the coefficients. 
\end{test}

\begin{minipage}[t]{0.49\linewidth}
\centering
\renewcommand{\arraystretch}{1.2}
\medskip
\captionof{table}{Evolution of the cost and computed $b_c$  and $c_c$ with random noises in the target in Test~\ref{testLin}.}
\medskip
{\begin{tabular}{ccccc} \hline
\%  & Cost  & It. &  $b_c$  & $c_c$\\
\hline
1         &   1.e-8         &    29            &  5.02223088   & 1.47406114\\ \hline
0.1        &    1.e-10      &    29           &  5.00636935  & 1.49523326 \\ \hline
0.01      &   1.e-12       &    29           &   4.99989752 & 1.50032299 \\\hline
0.001    &  1.e-14        &    29           &   4.99994089 &  1.50011971  \\ \hline
0           &   1.e-25       &    39           &  5                   &  1.5 \\  \hline
\end{tabular}}
\label{table.t11}
\end{minipage}
\hfill
\hskip 2mm
\begin{minipage}[t]{0.49\linewidth}
\vspace{0cm}
\includegraphics[width=\linewidth]{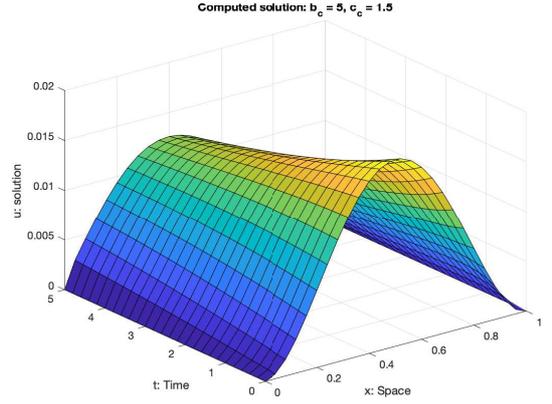}
 \captionof{figure}{Test~\ref{testLin}, computed solution.}
\label{Sol_Lin}
\end{minipage}

\bigskip
\begin{description}
\item[Quadratic function:] Let us now consider $a$ of the form $a(x) = bx^2+ cx +h$. The optimization problem \eqref{optpb1} is written in this case as follows: find  $(b,c,h)\in\mathcal{A}$, with  $\mathcal{A} = \{(b,c,h) : 0<\underline{b} < b<\overline{b}, \quad  0<\underline{c} < c<\overline{c}, \quad 0<\underline{h} < h<\overline{h}\}$ such that 
\begin{equation*}\label{optpb2}
J(b,c,h) \le J(b',c',h'), \quad \forall\, (b',c',h')\in \mathcal{A}.
\end{equation*}
\end{description}

\begin{figure}[h!]
\begin{minipage}[t]{0.49\linewidth}
\noindent
\includegraphics[width=\linewidth]{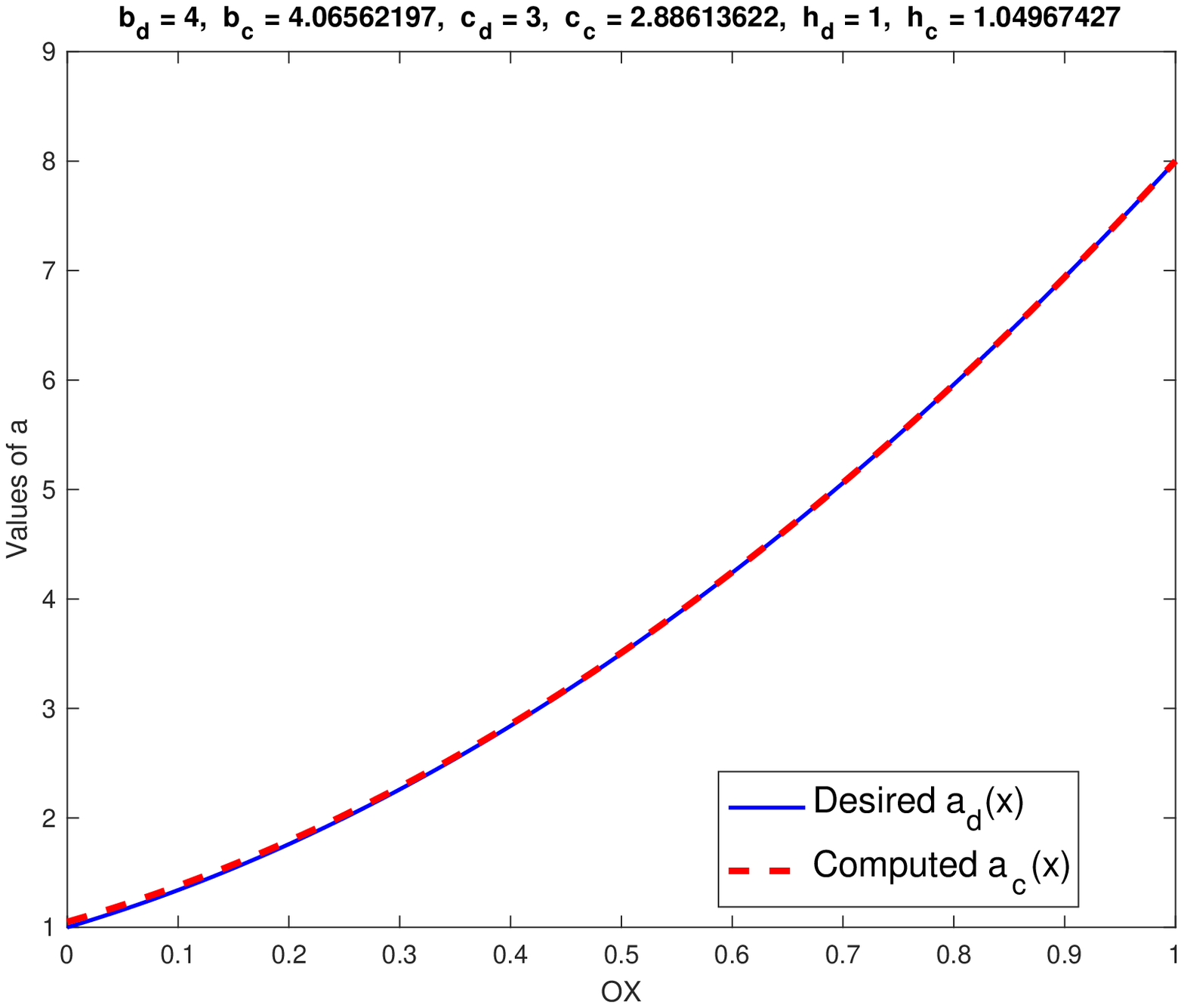}
\caption{Test~\ref{testquad}, $a_d$ solid, $a_c$ dashed lines.}
\label{fig.testquad}
\end{minipage}
\hfill
\begin{minipage}[t]{0.49\linewidth}
\noindent
\includegraphics[width=\linewidth]{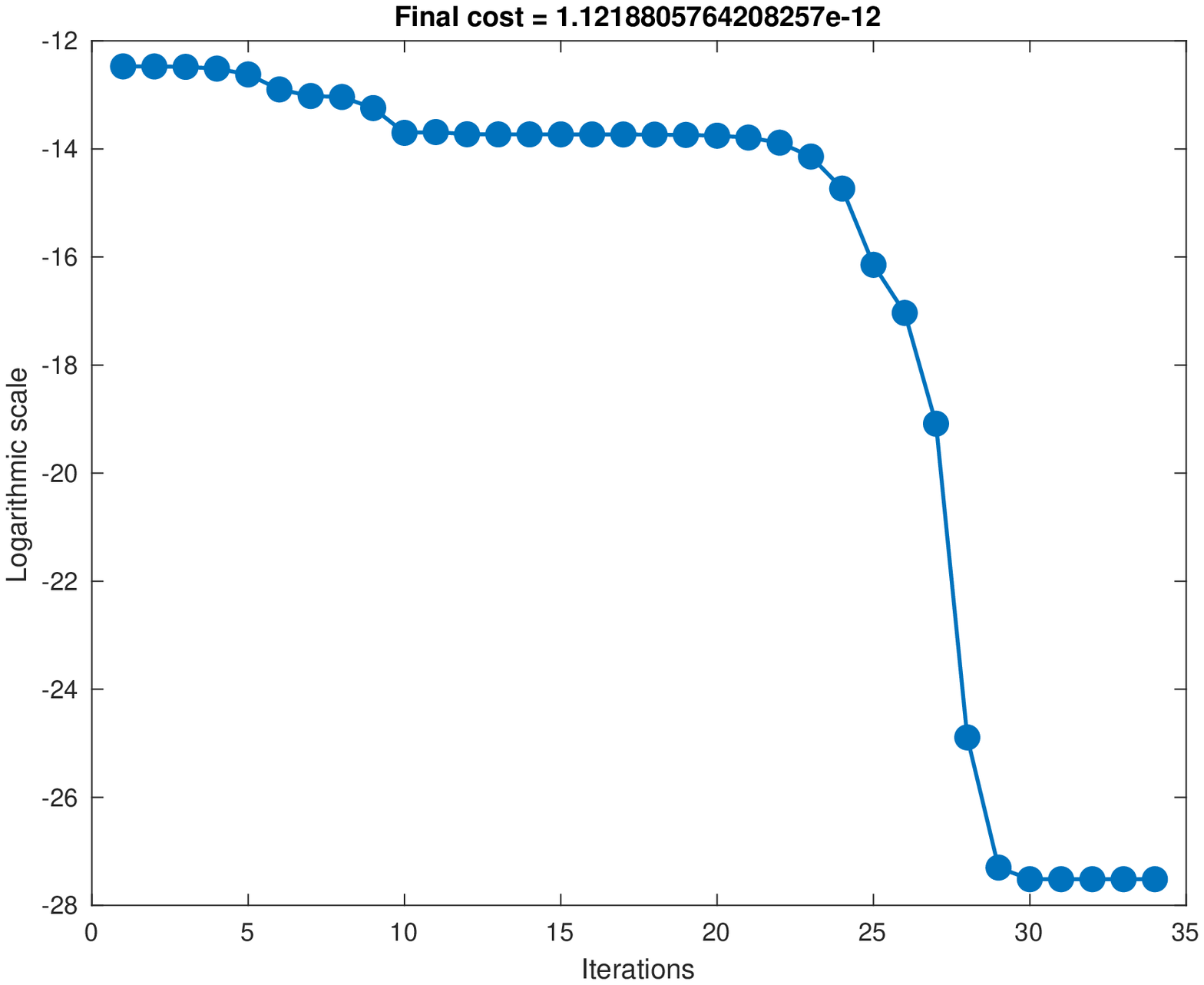}
\caption{Test~\ref{testquad}, evolution of the cost.}
\label{fig.testquad_f}
\end{minipage}
\end{figure}

\begin{test}\label{testquad}
Let us take $b_i = 3.5$, $c_i = 2.5$  and $h_i = 0.5$ as initial guess for the optimization algorithm and  the desired values to recover  $b_d=4$, $c_d=3$ and $h_d=1$.  The numerical results are shown in the Figures~\ref{fig.testquad} and~\ref{fig.testquad_f}, and the computed values are $b_c=4.06562197$, $c_c=2.88613622$ and $h_c=1.04967427$. The numerical results can be seen in Figures~\ref{fig.testquad} and~\ref{fig.testquad_f}. Table~\ref{table.t12} shows the evolution of the cost during the iterations of optimization algorithm and computed $b_c$  and $c_c$ with random noises in the target.


\end{test}

\begin{minipage}[t]{\linewidth}
\centering
\renewcommand{\arraystretch}{1.2}
\medskip
\captionof{table}{Evolution of the cost and computed $b_c$  and $c_c$ with random noises in the target in Test~\ref{testquad}.}
\medskip
{\begin{tabular}{cccccc} \hline
\%  & Cost  & Iterations &  $b_c$  & $c_c$ & $h_c$ \\
\hline
1\%         &   1.e-8         &    41            &  4.05971241   & 2.88655292 & 1.09105166
\\ \hline
0.1\%       &    1.e-10      &    33           &  4.06806662  & 2.88685441 &  1.04274887
\\ \hline
0.01\%      &   1.e-12       &    34           &   4.06565361 & 2.88623199 & 1.05038351
\\\hline
0.001\%    &  1.e-12        &    35           &   4.06567662 &  2.88615326 &  1.04953132
\\ \hline
0\%           &   1.e-12       &    33           &  4.06562197   &  2.88613622  &  1.04967427
\\  \hline
\end{tabular}}
\label{table.t12}
\end{minipage}


\end{document}